\documentclass[11pt]{amsart}

\allowdisplaybreaks[4]
\usepackage{graphicx} 
\usepackage{epsfig}
\usepackage{cite}
\usepackage{xypic}
\usepackage{latexsym,amssymb,amsxtra,amsmath,amsfonts}
\usepackage{verbatim}
\usepackage{color}
\usepackage{kpfonts}

\newtheorem{theorem}{Theorem}
\newtheorem{lemma}[theorem]{Lemma}
\newtheorem{proposition}[theorem]{Proposition}
\newtheorem{corollary}[theorem]{Corollary}

\theoremstyle{definition}

\newtheorem{remark}[theorem]{Remark}

\newcommand{\be}{\begin{equation}}
\newcommand{\ee}{\end{equation}}
\newcommand{\bea}{\begin{eqnarray}}
\newcommand{\eea}{\end{eqnarray}}
\newcommand{\ba}{\begin{array}}
\newcommand{\ea}{\end{array}}
\newcommand{\bean}{\begin{eqnarray*}}
\newcommand{\eean}{\end{eqnarray*}}

\newcommand{\beq}{\begin{equation}}
\newcommand{\eeq}{\end{equation}}
\newcommand{\beqa}{\begin{eqnarray}}
\newcommand{\eeqa}{\end{eqnarray}}
\newcommand{\beaa}{\begin{eqnarray*}}
\newcommand{\ben}{\begin{eqnarray*}}
\newcommand{\eaa}{\end{eqnarray*}}
\newcommand{\een}{\end{eqnarray*}}

\def\be{\beta}

\def\A{\mathcal{A}}
\def\C{\mathcal{C}}
\def\D{\mathcal{D}}

\def\L{\mathcal{L}}
\def\M{\mathcal{M}}
\def\O{\mathcal{O}}
\def\Q{\mathcal{Q}}
\def\R{\mathcal{R}}

\def\E{\mathcal{E}}

\def\ZZ{\mathbb{Z}}

\def\CC{\mathbb{C}}
\def\PP{\mathbb{P}}

\newcommand{\pa}{\partial}

\newcommand{\resi}{{\rm res\/}}

\makeatletter
\newcommand{\doublewidetilde}[1]{{%
  \mathpalette\double@widetilde{#1}%
}}
\newcommand{\double@widetilde}[2]{%
  \sbox\z@{$\m@th#1\widetilde{#2}$}%
  \ht\z@=.9\ht\z@
  \widetilde{\box\z@}%
  }
\makeatother

\begin{document}

\title[Extended HQEs]
{\textbf{ The extended D-Toda hierarchy }}
\author{Jipeng Cheng}
\address{School of Mathematics, China University of Mining and
  technology, Xuzhou, Jiangsu 221116, P.R. China}
\email{chengjp@cumt.edu.cn}

\author{Todor Milanov}
\address{Kavli IPMU (WPI), UTIAS, The University of Tokyo, Kashiwa, Chiba 277-8583, Japan}
\email{todor.milanov@ipmu.jp} 

\begin{abstract}
In a companion paper to this one, we proved that the Gromov--Witten theory of a Fano orbifold line of type $D$ is governed by a system of Hirota Bilinear Equations. The goal of this paper is to prove that every solution to the Hirota Bilinear Equations  determines a solution to a new integrable hierarchy of Lax equations. We suggest the name extended D-Toda hierarchy for this new system of Lax equations, because it should be viewed as the analogue of Carlet's extended bi-graded Toda hierarchy, which is known to govern the Gromov--Witten theory of Fano orbifold lines of type $A$. 
\end{abstract}

\maketitle
\tableofcontents
\section{Introduction}

\subsection{Background and motivation}
Motivated by Gromov--Witten theory, Dubrovin and Zhang have proposed a
general construction which associates an integrable hierarchy to every
semi-simple Frobenius manifold. The definition however is very
complicated and hence the study of these hierarchies is a very
challenging problem. Our strategy is to concentrate on the cases when
the Frobenius manifold corresponds to a semi-simple quantum cohomology of a complex
orbifold $X$, whose coarse moduli space $|X|$ is a projective
variety. We can further seperate these classes of Frobenius manifolds
according to the dimension of $X$. Based on the examples worked out in
the literature, one can speculate that in complex dimension 1, the
corresponding integrable hierarchies can be understood in terms of the
representation theory of generalized Kac--Moody Lie algebras.

Let us discuss the case when $\operatorname{dim}_\CC(X)=1$. The
quantum cohomology of $X$ is semi-simple if and only if 
the coarse moduli space of $X$ is $\PP^1$ (see \cite{Shir}). Let us
divide the orbifold lines into 
three groups depending on whether the orbifold Euler characteristic is
$>0$, $=0$, or $<0$. The orbifolds in these three groups will be called
respectively Fano, elliptic, and hyperbolic orbifold lines. The Fano
case is the easiest and nevertheless it is still unfinished. A Fano orbifold
line has the form $\PP^1_{a_1,a_2,a_3}$, that is, $\PP^1$ with 3
orbifold points with isotropy groups of orders $a_1,a_2$, and $a_3$,
such that $\tfrac{1}{a_1}+\tfrac{1}{a_2}+\tfrac{1}{a_3}>1$. Tripples
$(a_1,a_2,a_3)$ satisfying the above inequality are in one-to-one
correspondence with the Dynkin diagrams of type $ADE$. According to
\cite{MST} the corresponding hierarchy must be an extension of a
certain Kac--Wakimoto hierarchy. In the case $A$, the extension is
known (see \cite{C,MT,CvL}) and it is called the Extended
Bi-graded Toda Hierarchy. Our interest is in the case $D$. We divided
the problem into two parts. In the first part, we find the extension of 
the corresponding Kac--Wakimoto hierarchy in the form of Hirota Bilinear 
Equations. The second part, which is the goal of this paper, is to 
describe the extension in terms of Lax equations.

Let us point out that although the Kac--Wakimoto hierarchies have been
known for a while, it is still an open question to describe the flows
of these hierarchies in terms of Lax equations. There are many cases in
which the Lax equations are known -- usually the answer is a
reduction of some multicomponent KP hierarchy, but there is no general
construction that works for all Kac--Wakimoto hierarchies. In
particular, for the case $D$ in our project, the Lax equations of the
corresponding Kac--Wakimoto hierarchies were unkonwn, so we had to
construct them. We believe that our methods can be generalized and
that one should be able to construct the Lax equations of all
Kac--Wakimoto hierarchies of type $D$. Finally, let us point out that
the symbol of our Lax operator is very similar to the
Landau--Ginzburg potential used in the construction of Frobenius
structures on the orbit spaces of the extended Weyl groups in
\cite{DSZZ}. We can speculate that another possible generalization of our work is to
construct the integrable hierarchies corresponding to the semi-simple Frobenius
manifolds constructed in \cite{DSZZ}.

In the rest of the introduction we will focus on stating our results.
\subsection{Lax operators}\label{sec:Lax_op}
The idea of our construction is partially motivated by Shiota's
approach to the 2-component BKP hierarchy (see \cite{Shi}). 
Let $\R$ be the ring of formal power series in $\epsilon$ whose
coefficients are differential polynomials in the set of
$n+1$ variables 
$\Xi=\{a_1,\dots,a_{n-4}, \alpha,  q_2,
q_3,c_2,c_3\}$ and the functions $e^{\pm \alpha}$. Formally, $\R$ is defined by 
\ben
\R:=\CC[\xi^i\, (i\geq 0,\ \xi\in \Xi), e^{\pm \alpha}][\![\epsilon]\!],
\een
where $\xi^i$ is a formal variable. We identify
$\xi^0:= \xi$ for $\xi\in \Xi$. Let us define the derivation  $\pa_x$ 
\ben
\pa_x(P) =
\sum_{i=0}^\infty \frac{\pa P}{\pa \xi^i} \,
\xi^{i+1}.
\een
Note that $\xi^{i}=\pa_x^i(\xi)$.  The translation operator $\Lambda:=
e^{\epsilon \pa_x}$ acts naturally on the ring $\R$ and we put  $P[m]:=\Lambda^m(P)$.

Suppose that the ring $\R$ is equipped with two commuting
derivations $\pa_2$ and $\pa_3$ both commuting with $\pa_x$.
Given an operator series
\beq\label{series-A}
A=\sum_{j_1,j_2,j_3\in\mathbb{Z}}f_{j_1j_2j_3} \Lambda^{j_1}\partial_2^{j_2}\partial_3^{j_3}，
\eeq
where $f_{j_1,j_2,j_3}\in \R$ and the sum is possibly
infinite, we define the truncations $A_{a,\leq k}$, $A_{a,<k}$,
$A_{a,[k]}$, $A_{a,\geq k}$, and $A_{a,>k}$ by keeping only the terms in
the sum \eqref{series-A} for which $j_a$ is respectively $\leq k$,
$<k$, $=k$, $\geq k$, and $>k$ and truncating the remaining ones,
e.g.,   
\ben
A_{1,\geq k }=\sum_{j_1\geq k}
\sum_{j_2,j_3\in\mathbb{Z}
}f_{j_1j_2j_3} \Lambda^{j_1}\pa_2^{j_2}\pa_3^{j_3}. 
\een
Given another operator series
$B=\sum_{l_1,l_2,l_3\in\mathbb{Z}}g_{l_1l_2l_3} \Lambda^{l_1}\pa_2^{l_2}\pa_3^{l_3}$,
the operator composition of $A$ and $B$, whenever it makes sense, will
be denoted by $AB$ or $A\cdot
B$. If the sum \eqref{series-A} is finite and contains only terms for
which $j_2,j_3\geq 0$, then $A$ is called a {\em differential-difference}
operator. A differential-difference operator $A$ acts naturally on
the space of formal operator series of the type \eqref{series-A}.  We 
denote by $A(B)$ the operator series obtained by applying $A$ to the
coefficients of $B$, that is,  
\ben
A(B):=\sum_{l_1,l_2,l_3\in\mathbb{Z}}A(g_{l_1l_2l_3})\Lambda^{l_1}\pa_2^{l_2}\pa_3^{l_3}.  
\een
Finally, let us introduce the adjoint operation $\#$ 
\begin{eqnarray}
  A^\#=\sum_{j_1,j_2,j_3\in\mathbb{Z}}
  \Lambda^{-j_1}(-\pa_2)^{j_2}(-\pa_3)^{j_3}f_{j_1j_2j_3}， 
\end{eqnarray}
which obeys  $(AB)^\#=B^\#A^\#$ for any two operator series $A$ and
$B$ for which the composition $AB$ makes sense.

Let us denote by
$\E:=\R[\Lambda^{\pm 1},\pa_2,\pa_3 ]$ the ring of differential-difference operators and
define the following 4 operators in $\E$:  
\ben
\L  & := &  \Big(\sum_{i=1}^{n-3} (a_i\Lambda^i - \Lambda^{-i} a_i)\Big)
(\Lambda-\Lambda^{-1})+
\frac{1}{2}\,\pa_2^2 + \frac{1}{2}\,\pa_3^2 
+\frac{1}{4}(c_2-c_3)(\Lambda+\Lambda^{-1}) +\frac{1}{2}(c_2+c_3),\\
H_1 & := & \pa_2\pa_3 + q_1, \\
H_2 & := & (\Lambda-1) \pa_2 -q_2(\Lambda+1), \\
H_3 & := & (\Lambda+1) \pa_3-q_3(\Lambda-1),
\een
where $q_1 := -2(1+\Lambda)^{-1}(q_2 q_3)\in \R$
and $a_{n-3}:=\tfrac{1}{n-2} e^{(n-2)\alpha}$.

Let us introduce the following rings 
\ben
\mathcal{E}_{(\pm)}=\R[\pa_2,\pa_3]((\Lambda^{\mp 1})),\quad
\mathcal{E}_{(2)}=\R[\Lambda,\Lambda^{-1},\pa_3]((\pa_2^{-1})),\quad
\mathcal{E}_{(3)}=\R[\Lambda,\Lambda^{-1},\pa_2]((\pa_3^{-1}))
\een
and
\ben
\mathcal{E}_{(\pm)}^{0}=\R((\Lambda^{\mp 1})),\quad
\mathcal{E}_{(2)}^{0}=\R((\pa_2^{-1})),\quad
\mathcal{E}_{(3)}^{0}=\R((\pa_3^{-1})).
\een
Finally, let us denote by $\A H$ the left ideal in
$\A$ generated by 
$H_1,H_2,H_3$, where $\A$ could be any of the rings $\E_{(\pm)}$,
$\E_{(2)}$, or $\E_{(3)}$. Note that we always have a decomposition
into sum of vector spaces 
\beq\label{proj-A}
\A=\A^0+\A H,\quad \forall\ \A\in \{\E_{(\pm)}, \E_{(2)}, \E_{(3)}\}.
\eeq
We will work out a criteria for the derivations $\pa_2$ and $\pa_3$
that guarantees that \eqref{proj-A} is a direct sum decomposition,
that is,
$\A^0\cap \A H\neq \{0\}$. If \eqref{proj-A} is a direct sum
decomposition, then we denote by $\pi_\alpha$ $(\alpha=\pm,2,3)$ the
corresponding projection $\E_{(\alpha)}\to \E_{(\alpha)}^0$ and we
have the following recursion formulas: 
\begin{eqnarray*}
&&\pi_{\alpha}(\Lambda^{j_1\pm 1}\pa_2^{j_2}\pa_3^{j_3})=\Lambda^{\pm 1}\Big(\pi_\alpha(\Lambda^{j_1}\pa_2^{j_2}\pa_3^{j_3})\Big)\cdot\pi_{\alpha}(\Lambda^{\pm 1}),\\
&&\pi_{\alpha}(\Lambda^{j_1}\pa_2^{j_2+1}\pa_3^{j_3})=\pa_2\Big(\pi_{\alpha}(\Lambda^{j_1}\pa_2^{j_2}\pa_3^{j_3})\Big)+\pi_{\alpha}(\Lambda^{j_1}\pa_2^{j_2}\pa_3^{j_3})\cdot\pi_{\alpha}(\pa_2),\\
&&\pi_{\alpha}(\Lambda^{j_1}\pa_2^{j_2}\pa_3^{j_3+1})=\pa_3\Big(\pi_{a}(\Lambda^{j_1}\pa_2^{j_2}\pa_3^{j_3})\Big)+\pi_{\alpha}(\Lambda^{j_1}\pa_2^{j_2}\pa_3^{j_3})\cdot\pi_{\alpha}(\pa_3),
\end{eqnarray*}
which reduce the computation of $\pi_\alpha$ to the following cases:
\begin{eqnarray*}
&&\pi_{\pm}(\pa_2)=\iota_{\Lambda^{\mp 1}}Q_{2},\quad \pi_{\pm}(\pa_3)=\iota_{\Lambda^{\mp 1}}Q_{3},\\
&&\pi_2(\Lambda)=(\pa_2-q_2)^{-1}\cdot(\pa_2+q_2)=-1+2(\pa_2-q_2)^{-1}\cdot\pa_2,\quad \pi_2(\pa_3)=-\pa_2^{-1}\cdot q_1,\\
&&\pi_3(\Lambda)=-(\pa_3-q_3)^{-1}\cdot(\pa_3+q_3)=1-2(\pa_3-q_3)^{-1}\cdot\pa_3,\quad \pi_3(\pa_2)=-\pa_3^{-1}\cdot q_1,
\end{eqnarray*}
where $Q_2:=(\Lambda-1)^{-1} q_2(\Lambda+1)$, 
$Q_3:=(\Lambda+1)^{-1} q_3(\Lambda-1)$ and $\iota_{\Lambda}$
(resp. $\iota_{\Lambda^{-1}}$) denotes the 
Laurent series expansion at $\Lambda=0$ (resp. $\Lambda=\infty$). 
\begin{proposition}\label{prop:proj-A}
a) The decomposition \eqref{proj-A} is a direct sum of vector spaces
if and only if the derivations $\pa_2$ and $\pa_3$ satisfy the
following 0-curvature condition in $\E_{(\pm)}$
\ben
\pa_2(Q_3) -\pa_3(Q_2) = [Q_2,Q_3].
\een
 
b) There are unique derivations $\pa_2$ and $\pa_3$ such that
\eqref{proj-A} is a direct sum of vector spaces and
\ben
H_a\L\in \E H ,\quad
2\leq a\leq 3.
\een
\end{proposition}
\begin{remark} The 0-curvature condition can be justified as
follows: Note that $\pi_\alpha(\pa_2\pa_3)$ ($\alpha=\pm, 2,3 $) can
be computed using the recursion formulas from above in two different
ways: reducing the power of $\pa_2$ or reducing the power of $\pa_3$.
The two computation will agree if and only if
\ben
\pa_2\Big(\pi_{\alpha}(\pa_3)\Big)-
\pa_3\Big(\pi_{\alpha}(\pa_2)\Big)+
[\pi_\alpha(\pa_3),\pi_\alpha(\pa_2)]=0,
\een
which is equivalent to the 0-curvature condition. Similarly, the
projection $\pi_\alpha(\Lambda \pa_b)$ ($\alpha=\pm, 2,3,\ b=2,3$) can
be computed in two different ways, which will agree if and only if 
\ben
\pa_b\Big(\pi_{\alpha}(\Lambda)\Big)+
\pi_\alpha(\Lambda)\cdot\pi_a(\pa_b)=
\Lambda\Big(\pi_\alpha(\pa_b)\Big)\cdot\pi_\alpha(\Lambda).
\een
This relation is also equivalent to the 0-curvature condition.\qed
\end{remark}

\subsection{Dressing operators}\label{sec:dress_op}
There are unique operators $L_1=b_{1,0}\Lambda + \sum_{i=1}^\infty b_{1,i}
\Lambda^{1-i}$, $L_a= \pa_a +\sum_{i=1}^\infty b_{a,i} \pa_a^{1-i}$ ($a=2,3$)
with coefficients in $\R$ such that 
\ben
b_{1,0}=\exp \Big(\frac{(\Lambda-1)(n-2)}{\Lambda^{n-2}-1}(\alpha)\Big),\quad
b_{2,1}=b_{3,1}=0
\een
and 
\beq\label{projs-L}
\pi_+(\L)= \frac{1}{n-2} L_1^{n-2},\quad 
\pi_2(\L)= \frac{1}{2} L_2^2,\quad
\pi_3(\L) = \frac{1}{2} L_3^2.
\eeq
Let us construct 3 differential ring extensions $\R_i$ ($1\leq i\leq
3$) of $\R$. Put
\ben
\R_1:=\CC[\xi^j (\xi\in \Xi, j\geq0),e^{\pm \alpha}, e^{\pm \phi} ,
\psi_{1,1},\psi_{1,2},\dots][\![\epsilon]\!],
\een
and 
\ben
\R_a=\CC[\xi^j (\xi\in \Xi, j\geq0), e^{\pm \alpha}, \psi_{a,1},\psi_{a,2},\dots][\![\epsilon]\!],
\quad  2\leq a\leq 3.
\een
The derivation $\epsilon \pa_x$ is extended uniquely to a derivation of $\R_1$ in such a way
that $\epsilon\pa_x(\phi)=\alpha$ and $L_1= S_1 \Lambda S_1^{-1}$, where
\ben
S_1 = \psi_{1,0} + \sum_{i=1}^\infty \psi_{1,i} \Lambda^{-i},\quad
\psi_{1,0} = e^{\frac{(n-2)\epsilon \pa_x}{1-\Lambda^{n-2}} (\phi)}. 
\een
More explicitly, by comparing the coefficients
in front of $\Lambda^{1-i}$ we get that the identity $L_1 S_1 =
S_1\Lambda$ will be satisfied if we define 
\beq\label{x-der-S1}
\epsilon \pa_x \Big( \frac{\psi_{1,i}}{\psi_{1,0}}\Big) =
  \frac{\epsilon \pa_x}{1-e^{\epsilon \pa_x}}
  \left(\sum_{s=1}^i b_{1,s} \frac{\psi_{1,i-s}[1-s]}{\psi_{1,0}}
  \right).
\eeq
Similarly, there exists a unique extension of the derivation $\pa_a$ $(a=2,3)$ to a derivation of
$\R_a$, such that, $L_a T_a=T_a \pa_a$, where
\ben
T_a=1+\psi_{a,1}\pa_a^{-1}+\psi_{a,2}\pa_a^{-2} +\cdots.
\een
Such an extension exists and the derivation $\pa_a$ $(a=2,3)$ is uniquely determined. Indeed,
substituting $L_a=\pa_a+\sum_{i=1}^\infty b_{a,i} \pa_a^{1-i}$ and the
above expansion of $T_a$ in $L_a T_a=T_a \pa_a$, comparing the
coefficients in front of $\pa_a^{-k}$ for $k\geq 1$, and using that
$b_{a,1}=0$ we get  
\ben
\pa_a(\psi_{a,k}) +
\sum_{i=2}^{k+1}
\sum_{s=0}^{k+1-i}
{1-i\choose s} b_{a,i} \pa_a^s(\psi_{a,k+1-i-s}) = 0,\quad
k\geq 1.
\een
This formula alows us to define recursively $\pa_a(\psi_{a,k})$ for
all $k\geq 1$. We are going to prove that the derivations $\epsilon
\pa_x$, $\pa_2$, and $\pa_3$ can be extended uniquely to pairwise
commuting derivations of $\R_i$ ($1\leq i\leq 3$), such that, the
following conjugation relations hold:
for the operator $S_1$
\begin{align}
  \label{conj-S1}
  S_1 \Lambda S_1^{-1} =L_1,\quad
  S_1 \pa_a S_1^{-1}=\pa_a-Q_a\ (a=2,3),
\end{align}
for the operator $T_2$
\begin{align}
  \label{conj-T2}
  T_2 (\Lambda-1) T_2^{-1}=(\pa_2+q_2)^{-1} H_2 ,\quad
  T_2 \partial_2 T_2^{-1}=L_2, \quad 
  T_2 \pa_3 T_2^{-1}=\pa_2^{-1} H_1,
\end{align}
and for the operator $T_3$
\begin{align}
  \label{conj-T3}
  T_3 (\Lambda+1) T_3^{-1}= (\pa_3+q_3)^{-1} H_3 ,\quad
  T_3 \pa_2 T_3^{-1} = \pa_3^{-1} H_1 ,\quad
 T_3 \partial_3 T_3^{-1}=L_3.
\end{align}
We need to modify slightly the definition of $T_2$ and $T_3$. It turns
out that the following proposition holds:
\begin{proposition}\label{prop:S_a}
  There exists an operator $S_a\in 1+ \R_a[\![\pa_a^{-1}]\!]$
  ($a=2,3$), such that, $S^{-1}_a T_a$ commutes with the derivations
  $\epsilon\pa_x,\pa_2$, and $\pa_3$ and $S_a^\#=\pa_a S_a^{-1} \pa_a^{-1}$. 
\end{proposition}
Note that the conjugation formulas \eqref{conj-T2} and \eqref{conj-T3}
remain valid if we replace $T_a$ by $S_a$. The operators $S_i$ ($1\leq
i\leq 3$) will be called {\em dressing operators}. 

\subsection{Lax equations}\label{sec:Lax_eqns}
Let $L_i$ ($1\leq i\leq 3$) be the operator series defined by
\eqref{projs-L}. Let us define
\ben
B_{1,k}:= \sum_{m=0}^\infty \Big(
\left(L_1^k\, \Lambda^{-2m-1}\right)_{1,\geq 0} +
\left(\Lambda^{2m+1} (L_1^\#)^k\right)_{1,<0}
\Big)\, (\Lambda-\Lambda^{-1}),\quad k\geq 1,
\een
and 
\ben
B_{a,2l+1}:=(L_2^{2l+1})_{2,\geq 0},\quad a=2,3,\quad l\geq 0.
\een
Using the dressing operators we define
\ben
\log L_1 := S_1 \epsilon\pa_x S_1^{-1}=\epsilon\pa_x - \ell_1,\quad \ell_1\in \R[\![\Lambda^{-1}]\!],
\een
\ben
\log \Big( (\pa_2+q_2)^{-1} H_2 +1\Big)
:= S_2 \epsilon\pa_x S_2^{-1}
=\epsilon \pa_x - \ell_2,\quad
\ell_2\in \R[\![\pa_2^{-1} ]\!] \pa_2^{-1},
\een
and
\ben
\log \Big( (\pa_3+q_3)^{-1} H_3 -1\Big)
:= S_3 \epsilon\pa_x S_3^{-1}
=\epsilon \pa_x - \ell_3,\quad
\ell_3\in \R[\![\pa_3^{-1} ]\!] \pa_3^{-1},
\een
where $\ell_i:=\epsilon\pa_x(S_i)\cdot S_i^{-1}$ and the fact that the
coefficients of $\ell_i$ belong to $\R$ will be established later on. 
Put
\ben
A_{1,k}^+:=\frac{1}{(n-2)^k k!} L_1^{(n-2)k} (\epsilon \pa_x -\ell_1-h_k),
\een
\ben
A_{1,k}^-:= \iota_\Lambda (\Lambda-\Lambda^{-1})^{-1} (A_{1,k}^+)^\#
(\Lambda-\Lambda^{-1}),
\een
and
\ben
A_{a,k}:=\frac{1}{2^k k!} L_a^{2k} (\epsilon \pa_x-\ell_a),\quad a=2,3,
\een
where $k\geq 1$ and $h_k=\tfrac{1}{n-2}\left(1+\tfrac{1}{2}+\cdots
  +\tfrac{1}{k}\right)$. Let us define
$B_{0,k}= B_{0,k,1}+B_{0,k,2}+B_{0,k,3}$, where 
\begin{align}
   B_{0,k,1}& =\sum_{m=0}^{+\infty}
   \left(
   \Big(A_{1,k}^+\cdot\Lambda^{-2m-1}\Big)_{1,\geq 0}+
   \Big(A_{1,k}^-\cdot\Lambda^{2m+1}\Big)_{1,<0}
   \right)\cdot(\Lambda-\Lambda^{-1}),\nonumber\\
     B_{0,k,2} &=\Big(A_{2,k}\Big)_{2,>0}+
     \frac{1}{2}\Big(A_{2,k}\Big)_{2,[0]}\cdot (1+\Lambda^{-1}),\nonumber\\
     B_{0,k,3}& =\Big(A_{3,k}\Big)_{3,>0}+
     \frac{1}{2}\Big(A_{3,k}\Big)_{3,[0]}\cdot (1-\Lambda^{-1}). \nonumber
\end{align}
Our first result can be stated as follows.
\begin{theorem}\label{thm:Lax}
Let $(i,k)$ be an arbitrary pair of non-negative integers, such that $0\leq i\leq 3$,
$k\geq 1$, and $k=2l+1$ is odd for $i=2,3$. 

a) There exists a unique derivation $\pa_{i,k}:\R\to \R$ commuting with
$\pa_x$ such that
\begin{align}\nonumber
  \pa_{i,k}(\L) - [B_{i,k},\L] & \in \E H, \\
  \nonumber
\pa_{i,k} (H_a) -[B_{i,k},H_a] & \in  \E H ,\quad a=2,3.
\end{align}

b) The set of derivations $\{\pa_{i,k}\}$ pairwise commute. 
\end{theorem} 

The integrable hierarchy defined by the infinite set of commuting
derivations in Theorem \ref{thm:Lax} will be called the {\em Extended
  D-Toda Hierarchy}. Note that $\pi_+(\pa_{i,k}(M))=\pa_{i,k}(M)$ for
$M=\L,H_a(1\leq a\leq 3)$. Therefore, the conditions in part a) of
Theorem \ref{thm:Lax} yield the following formula:
\beq\label{Lax-pi_+}
\pa_{i,k} (M) = \pi_+[B^+_{i,k},M],\quad M\in \{\L,H_1,H_2,H_3\},
\eeq
where $B^+_{i,k}:=\pi_+(B_{i,k})$. We will prove also that
formula \eqref{Lax-pi_+} provides an equivalent 
formulation of the Extended D-Toda Hierarchy (see Lemma
\ref{le:der_R}).  

The Extended D-Toda hierarchy can be formulated also in terms of
the dressing operators. In other words, the derivations $\pa_{i,k}$
can be extended to $\R_a$ ($1\leq a\leq 3$) and the extended
derivations are pairwise commuting. The formulas for the derivatives
of the dressing operators can be found in Section
\ref{sec:dress_evol}. 
  
\subsection{Hirota bilinear equations}\label{sec:HBEs}
Let $\mathbf{t} =
(\mathbf{t}_0,\mathbf{t}_1, \mathbf{t}_2, \mathbf{t}_3)  $ be 4
sequences of formal variables of the following form:
\ben
\mathbf{t}_i=(t_{i,k})_{k\geq 1} \quad (i=0,1),\quad
\mathbf{t}_a=(t_{a,{2k-1}})_{k\geq 1},\quad (a=2,3).
\een
The 2 variables $y_a:=t_{a,1}$ ($2\leq a\leq 3$) will play a special
role. Let $\O_\epsilon:= \O(\CC)[\![\epsilon]\!]$ be the ring of
formal power series in $\epsilon$ whose coefficients are holomorphic
functions on $\CC$. Let $x$ be the standard coordinate function on
$\CC$.
The ring $\O_\epsilon[\![\mathbf{t}]\!]$ is equipped with a
translation operator 
\ben
\Lambda:\O_\epsilon[\![\mathbf{t}]\!]\to
\O_\epsilon[\![\mathbf{t}]\!],\quad 
\Lambda(f)(x,\mathbf{t}):=f(x+\epsilon,\mathbf{t})
\een
and an infinite set of pairwise commuting differentiations
$\tfrac{\partial}{\partial x}$ and $\tfrac{\partial}{\partial
  t_{i,k}}$. Note that a ring homomorphism 
\ben
\phi:\R \to \O_\epsilon[\![\mathbf{t}]\!],\quad
\mbox{such that}\quad
\phi\circ \partial_x = \frac{\pa}{\pa x}\circ \phi
\een
is uniquely determined by a set of $n+1$ functions in $\O_\epsilon[\![\mathbf{t}]\!]$
\ben
\phi(\xi)=
\xi(x,\mathbf{t}) \quad \xi\in \Xi.
\een
We say that the set of functions $\xi(x,\mathbf{t})$ ($\xi\in \Xi$) is
a solution to the Extended D-Toda Hierarchy if
the corresponding ring homomorphism satisfies
\ben
\phi\circ \pa_a= \tfrac{\pa}{\pa y_a}\circ \phi\ (a=2,3),\quad
\phi\circ \pa_{i,k} = \tfrac{\pa }{\pa t_{i,k}}\circ \phi.
\een
The second goal of our paper is to prove that solutions of the above
type can be constructed from a system of Hirota Bilinear Equations
(HBEs), which we describe now.

Suppose that 
$\tau(x,\mathbf{t}) \in \O_\epsilon[\![\mathbf{t}]\!] $ is an arbitary
invertible formal power series, i.e., $\tau(x,0)\in
\O(\CC)[\![\epsilon]\!]$ is a formal power series in $\epsilon$
whose leading order term is a holomorphic function on $\CC$ that has 
no zeros. 
Let $\D_\epsilon:=\D(\CC)[\![\epsilon]\!]$ be the ring of
formal power series in $\epsilon$ whose coefficients belong to the
ring $\D(\CC)$ of  holomorphic differential operators on $\CC$. We
extend the anti-involution $\#$ to operator series of the form \eqref{series-A}
with coefficients $f_{j_1,j_2,j_3}\in\D_\epsilon[\![\mathbf{t}]\!]$, so that
$(\partial_x)^\#=-\partial_x$. 
Let us introduce the formal power series 
\ben
\Psi^+_1(x,\mathbf{t},z)=\psi^+_1(x,\mathbf{t},z)e^{\xi_1(\mathbf{t},z)}z^{x/\epsilon-\frac{1}{2}}
\quad\mbox{and}\quad
\Psi^-_1(x,\mathbf{t},z)=z^{-x/\epsilon-\frac{1}{2}} e^{-\xi_1(\mathbf{t},z)}\psi^-_1(x,\mathbf{t},z)
\een
taking values in respectively 
$\D_\epsilon(\!(z^{-1})\!) [\![\mathbf{t}]\!] z^{x/\epsilon-\frac{1}{2}}$ and 
$z^{-x/\epsilon-\frac{1}{2}}\D_\epsilon(\!(z^{-1})\!) [\![\mathbf{t}]\!]$, 
where
\ben
&&
\xi_1(\mathbf{t},z)=
\sum_{k\geq 1}\left(t_{1,k}z^k+t_{0,k} 
\,(\epsilon\pa_x-h_k)\, 
\frac{z^{(n-2)k}}{(n-2)^k k!} \right),\\
&&
\psi^{\pm}_1(x,\mathbf{t},z)=
\frac{e^{\mp \sum_{k\geq 1} \frac{z^{-k}}{k} \partial_{t_{1,k}}}\, 
\tau(x\mp  \epsilon,\mathbf{t})}{\tau(x,\mathbf{t})}=
\sum_{k=0}^\infty
\psi^{\pm}_{1,k}(x,\mathbf{t})z^{-k},
\een
where $h_k=\frac{1}{n-2}\left(1+\frac{1}{2}+\cdots+\frac{1}{k}\right)$.
Similarly, let us define for $a=2,3$ the formal series
\ben
\Psi^+_a(x,\mathbf{t},z)=\psi^+_a(x,\mathbf{t},z)e^{\xi_a(\mathbf{t},z)}
\quad\mbox{and}\quad  
\Psi^-_a(x,\mathbf{t},z)=e^{-\xi_a(\mathbf{t},z)}\psi^-_a(x,\mathbf{t},z)
\een
taking values in $\D_\epsilon(\!(z^{-1})\!) [\![\mathbf{t}]\!] $, where 
\ben
&&
\xi_a(\mathbf{t},z)=\sum_{l\geq 1}\left(
t_{a,2l-1}z^{2l-1}+t_{0,l}\,  \epsilon\pa_x\, \frac{z^{2l}}{2^l l!}  \right),\\  
&&
\psi^{\pm}_a(x,\mathbf{t},z)=
\frac{ e^{\mp 2\sum_{l\geq 1}
\frac{z^{-2l+1}}{2l-1}\, \partial_{t_{a,l}}  }
\tau(x,\mathbf{t})}{\tau(x,\mathbf{t})}=\sum_{k=0}^\infty \psi_{a,k}^{\pm}(x,\mathbf{t})z^{-k}.
\een
The Hirota Bilinear Equations of the Extended D-Toda hierarchy
are given by the following system of quadratic equations 
\beqa\label{hbe-psi}
&&
{\rm  Res}_{z=0}\frac{z^{(n-2)k}}{(n-2)^kk!}\frac{dz}{z}
\left(\Psi_1^+(x,\mathbf{t},z)\Psi_1^-(x+m\epsilon,\mathbf{t}',z)+
  \left(\Psi_1^+(x+m\epsilon,\mathbf{t}',z)\Psi_1^-(x,\mathbf{t},z)\right)^\#\right)
\\
\nonumber
&=&
{\rm Res}_{z=0}\frac{z^{2k}}{2^k k!}\frac{dz}{2z}
\left(\Psi_2^+(x,\mathbf{t},z)\Psi_2^-(x+m\epsilon,\mathbf{t}',z)-
  (-1)^{m}\Psi_3^+(x,\mathbf{t},z)\Psi_3^-(x+m\epsilon,\mathbf{t}',z)\right),
\eeqa
where $k\geq 0$ and $m$ are arbitary integers and 
$\Psi^{\pm}_i(x+m\epsilon,\mathbf{t}',z):=\Lambda^m(\Psi^{\pm}_i(x,\mathbf{t}',z))$. If
the above equations are satisfied then we say that 
$\tau$ is a {\em tau-function} of the Extended D-Toda hierarchy, the
formal series 
$\Psi^{\pm}_1$ and $\Psi_a:=\Psi_a^+$ ($a=2,3$) will be 
called {\em wave functions}, and the operator series
\begin{align}
  \nonumber
  S_1^+(x,\mathbf{t},\Lambda)= & \sum_{j=0}^\infty
                                 \psi_{1,j}^+(x,\mathbf{t})\Lambda^{-j}\\
  \nonumber
  S_1^-(x,\mathbf{t},\Lambda)= & \sum_{j=0}^\infty\Lambda^{-j}
                                 \psi_{1,j}^-(x,\mathbf{t})\\
  \nonumber
  S_a(x,\mathbf{t},\pa_a)= & \sum_{j=0}^\infty
                             \psi_{a,j}^+(x,\mathbf{t})\pa_a^{-j}\quad (a=2,3)
\end{align}
will be called {\em wave operators}. Let us introduce also the
following auxiliarly Lax operators: 
\ben
&&
L_1^+(x,\mathbf{t},\Lambda):=
S_1^+(x,\mathbf{t},\Lambda)\cdot\Lambda\cdot S_1^+(x,\mathbf{t},\Lambda)^{-1}=:
u_{1,0}^+(x,\mathbf{t})\Lambda+\sum_{j=1}^\infty u_{1,j}^+(m,\mathbf{t})\Lambda^{1-j},\\
&&
L_1^-(x,\mathbf{t},\Lambda):=
S_1^-(x,\mathbf{t},\Lambda)^\#\cdot\Lambda^{-1}\cdot
\left(S_1^-(x,\mathbf{t},\Lambda)^\#\right)^{-1}=:
u_{1,0}^-(x,\mathbf{t})\Lambda^{-1}+\sum_{j=1}^\infty u^-_{1,j}(x,\mathbf{t})\Lambda^{j-1},\\
&&
L_a(x,\mathbf{t},\pa_a):=
S_a(x,\mathbf{t},\pa_a)\cdot\pa_a^{-1}\cdot
S_a(x,\mathbf{t},\pa_a)^{-1}=:
\pa_a+\sum_{j=1}^\infty u_{a,j}(x,\mathbf{t})\pa_a^{-j}\quad (a=2,3).
\een
To avoid cumbersome notation we put $L_1(x,\mathbf{t},\Lambda) :=
L_1^+(x,\mathbf{t},\Lambda)$ and $u_{1,j} = u^+_{1,j}$.
\begin{theorem}\label{thm:HBE-sol}
  If $\tau(x,\mathbf{t})$ is a tau-function of the Extended D-Toda
  hierarchy, then there is a uniquely determined solution of the
  Extended D-Toda hierarchy, such that,
  \ben
  c_a(x,\mathbf{t}) = \partial^2_{t_{a,1}} \log
  \tau(x,\mathbf{t}),\quad 
  q_a(x,\mathbf{t}) = \partial_{t_{a,1}} \log
    \frac{\tau(x,\mathbf{t})}{\tau(x+\epsilon,\mathbf{t}) }
 \quad (a=2,3),
  \een
  and the Lax operator
  \ben
\mathcal{L} & = & \frac{\pa_2^2}{2}+\frac{\pa_3^2}{2}+\left(
  \left( \tfrac{L_1^{n-2}}{n-2} \, \sum_{m=0}^\infty
  \Lambda^{-2m-1} \right)_{1,>0} -
  \left(\sum_{m=0}^\infty
  \Lambda^{2m+1}\,  \tfrac{(L_1^\#)^{n-2} }{n-2}\right)_{1,<0} \right)\,
(\Lambda-\Lambda^{-1})+ \\
&&
+\frac{1}{4}(c_2-c_3)(\Lambda+\Lambda^{-1}) +\frac{1}{2}(c_2+c_3).
\een
\end{theorem}
The proof of Theorem \ref{thm:HBE-sol} will be given in Sections
\ref{sec:HBE-lax-1} and \ref{sec:HBE-lax-2}. In fact, Sections
\ref{sec:HBE-lax-1} and \ref{sec:HBE-lax-2} can be read
independently. Our strategy is to construct Lax operators and derive
their evolution equations directly from the HBEs. In fact, the system
of Lax equations \eqref{Lax-pi_+} was discovered exactly in this way.
 
The inverse of Theorem \ref{thm:HBE-sol} is also true, i.e., we have the
following theorem:
\begin{theorem}\label{thm:inverse HBE-sol}
Any solution to the Extended D-Toda hierarchy in
$\O_\epsilon[\![\mathbf{t}]\!]$ is obtained from a tau-function.
\end{theorem}
The proof of Theorem \ref{thm:inverse HBE-sol} has two steps. 
The first one is to prove that for any solution of the Lax equations
\eqref{Lax-pi_+}  the corresponding dressing operators satisfy the
bilinear relations (\ref{hbe-psi}). Note that the wave function can be
expressed only in terms of the dressing operators, so we can interpret
\eqref{hbe-psi} as a system of bilinear relations for the dressing
operators. Using Taylor's series expansion, it is easy to check that
\eqref{hbe-psi} is equivalent to  (\ref{ext-id}), so our proof of 
 Proposition \ref{prop:ext_rel} (see Section \ref{sec:ext_rel})
 contains also the first step in the proof of the existence of a
 tau-function. The second step is to show the exsitence of the
 tau-function.  The argument is essentially the same as in \cite{Mi3}
 (see also \cite{CvL}) and we leave the details to the interested reader.   

\medskip

{\bf Acknowledgements.}
T.M. would like to thank Mikhail Kapranov for a very useful discussion
on localization of difference operators, which helped us to define
the ring of rational difference operators in Section \ref{sec:rdo}. We
would like to thank also Atsushi Takahashi for letting us know the
reference \cite{Shir} about quantum cohomology of orbifold lines. 
The work of T.M. is partially supported by JSPS Grant-In-Aid (Kiban C)
17K05193 and by the World Premier International Research Center
Initiative (WPI Initiative),  MEXT, Japan.

\section{Projections}
The first goal of this section is to prove Proposition
\ref{prop:proj-A}.  The second goal is to establish several key properties of the rings
$\E_{(\alpha)}$ $(\alpha\in \{\pm,2,3\})$ and the corresponding
projections $\pi_\alpha:\E_{(\alpha)}\to \E_{(\alpha)}^0$. We prove
two propositions which will be used in an essential way in the rest of
the paper. While the second proposition (see Proposition
\ref{prop:proj_pi}) is straightforward to prove, the first one (see
Proposition \ref{prop:ker_pi_a}) will require some
preparation. Namely, we will have to introduce rings of rational
difference operators.

\subsection{Proof of Proposition \ref{prop:proj-A}}
\proof
a) The neccessity of the condition is clear because 
\ben
\pa_2(Q_3) -\pa_3(Q_2) - [Q_2,Q_3] = [\pa_2-Q_2,\pa_3-Q_3] \in
\E_{(+)}^0\cap \E_{(+)}H.
\een
Suppose that the 0-curvature condition is satisfied. Let us prove that
\eqref{proj-A} is a direct sum decomposition for $\A=\E_{(+)}$. The
argument for $\E_{(-)}$ is identical. Note that $H_1$ already belongs
to the ideal in $\E_{(+)}$ generated by $H_2$ and $H_3$. We have to
prove that if $P:=a (\pa_2-Q_2)+ b (\pa_3-Q_3)\in \E_{(+)}^0$ then $P=0$. Note that
$a$ and $b$ must have the form
\ben
a=\sum_{i=0}^m a_i(\Lambda,\pa_3) \pa_2^i,\quad
b=\sum_{j=0}^{m+1} b_j(\Lambda,\pa_3)\pa_2^j,
\een
where the coefficients $a_i,b_j\in \R[\pa_3](\!(\Lambda^{-1})\!)$. We
argue by induction on $m$. If $m=0$ then by comparing the coefficients
in front of $\pa_2$ we get $a_0=-b_1 (\pa_3-Q_3)$. Using the
0-curvature condition we get $b_1 Q_2 (\pa_3-Q_3) \in
\E_{(+)}^0$. This implies that $b_1=0$, otherwise we write
$b_1=\sum_{i=0}^k b_{1,i}\pa_3^i$ with $b_{1,k}\neq 0$. The
coefficient in front of the highest power of $\pa_3$ is $b_{1,k} Q_2$
and it must vanish -- contradiction because $\E_{(+)}^0$ does not have
0-divisors. If $m>0$, then by comparing the coefficeints in front of
$\pa_2^{m+1}$ we get $a_m=-b_{m+1}(\pa_3-Q_3)$. Using the 0-curvature
condition we get
\ben
a_m \pa_2^m (\pa_2-Q_2) + b_{m+1} \pa_2^{m+1}  (\pa_3-Q_3) =
b_{m+1} \Big(
[Q_3,\pa_2^m] (\pa_2-Q_2) + \pa_2^m \cdot Q_2\cdot (\pa_3-Q_3)\Big). 
\een
Therefore we can write $P$ as a linear combination of $ \pa_2-Q_2$ and
$\pa_3-Q_3$ whose coefficients are differential operators of $\pa_2$
of orders respectively $<m$ and $<m+1$. Recalling the inductive
assumption we get $P=0$.

Let us move to the case $\A=\E_{(2)}$. After a direct computation we
get that the 0-curvature condition is equivalent to the following
identities
\beq\label{pa_a(q_b)1}
\pa_2(q_3)=\pa_3(q_2) = \frac{e^{\epsilon \pa_x}-1}{e^{\epsilon
    \pa_x}+1} (q_2q_3)
\eeq
Using \eqref{pa_a(q_b)1} we get the following relation:
\ben
\pa_2\cdot  H_3 + q_3 H_2 = (\Lambda+1)\cdot H_1.
\een
Therefore, the left ideal $\E_{(2)}H$ is generated by $H_1$ and
$H_2$. Put
\ben
\widetilde{H}_1:=\pa_2^{-1}H_1=\pa_3+\pa_2^{-1} q_1
\een
and
\ben
\widetilde{H}_2:=(\pa_2+q_2)^{-1}H_2 = (1+\pa_2^{-1} q_2)^{-1} (1-\pa_2^{-1}
q_2)\Lambda-1.
\een
We have to prove that if
$P=a \widetilde{H}_1 + b \widetilde{H}_2\in \E_{(2)}^0$ then $P=0$.
Similarly to the above case, since the coefficients in front of the
highest power of  $\pa_3$ in $P$ vanish, we get that $a$ and $b$ must
have the form   
$a=\sum_{i=0}^m a_i\pa_3^i$ and
$b=\sum_{j=0}^{m+1} b_j \pa_3^j$, where $a_i, b_j\in
\R[\Lambda^{\pm 1}](\!(\pa_2^{-1})\!)$.
The coefficients in front of $\pa_3^{m+1}$ in $P$ must vanish
$\Rightarrow$ $a_m=-b_{m+1} \widetilde{H}_2$. We are going to prove that
$[\widetilde{H}_1, \widetilde{H}_2]=0.$ Assuming this fact we get
\ben
a_m \pa_3^m \widetilde{H}_1 + b_{m+1} \pa_3^{m+1} \widetilde{H}_2 =
b_{m+1}\Big(
[\pa_3^m, \widetilde{H}_2] \, \widetilde{H}_1 +
\pa_3^m\cdot (-\pa_2^{-1} q_1)\cdot \widetilde{H}_2
\Big). 
\een
Therefore we can complete the proof by induction on $m$. Let us prove
that $[\widetilde{H}_1, \widetilde{H}_2]=0$. We claim that the
vanishing of this commutator follows from the identities
\beq
\label{q-der}
\pa_2(q_3)=\pa_3(q_2) = q_1 + q_2 q_3 = -q_2 q_3 -q_1[1],
\eeq
where the last identity is just the definition of $q_1$, while the
rest of the identities, being equivalent to \eqref{pa_a(q_b)1}, are
consequence of the 0-curvature condition. It is convenient to put
$M=(1+\pa_2^{-1} q_2)^{-1} (1-\pa_2^{-1} q_2)$. The vanishing of
the commutator $[\widetilde{H}_1, \widetilde{H}_2]$ is equivalent to
\beq\label{pa_3(M)}
\pa_3(M) = -\pa_2^{-1} q_1  M + M \pa_2^{-1} q_1[1].
\eeq
We have
\ben
\pa_3(M) = -(1+\pa_2^{-1} q_2)^{-1} \pa_2^{-1} \Big(\pa_3(q_2)
M + \pa_3(q_2)\Big).
\een
Substituting the above formula in \eqref{pa_3(M)}, collecting the
terms that have $M$ as the rightmost factor, and multiplying both
sides of the equation from the left by  $\pa_2 (1+\pa_2^{-1} q_2)$, we
get that \eqref{pa_3(M)}  is equivalent to 
\ben
(-\pa_3(q_2) + q_1 + q_2 \pa_2^{-1} q_1) M -\pa_3(q_2)=
\pa_2 (1+\pa_2^{-1} q_2) M \pa_2^{-1} q_1[1].
\een
Using formula \eqref{q-der} we transform the above equation into
\ben
- q_2 \pa_2^{-1} q_3  (\pa_2-q_2) -\pa_3(q_2)=
(\pa_2-q_2) \pa_2^{-1} q_1[1].
\een
This however follows easily from \eqref{q-der}.

b) Let us first prove the uniqueness. We will argue that the
derivations $\pa_2$ and $\pa_3$ are uniquely determined from the
0-curvature condition and the 4 projection constraints 
\ben
\pi_+(H_2 \L)=\pi_+(H_3\L)=0,\quad \pi_2(H_2 \L)=0,\quad \pi_3(H_3 \L)=0.
\een
We already know that the 0-curvature condition is equivalent to
formulas \eqref{q-der}, that is, $\pa_a(q_b)$ for $a\neq b$ is uniquely
fixed. Let us prove that
\beq\label{pa_a(q_a)1}
\pa_a(q_a) = \frac{1}{2}\,\left(1-e^{\epsilon\pa_x}\right)(c_a),\quad a=2,3.
\eeq
The proof of the two formulas is  identical, so let us consider only
the case $\pa_2(q_2)$. Note that
$\pi_2(\L)  =\tfrac{\pa_2^2}{2} + c_2 + O(\pa_2^{-1})$ and that
\ben
(\pa_2-q_2)^{-1} H_2 = \Lambda - 1 -2\sum_{i=1}^\infty (\pa_2^{-1} q_2)^i.
\een
Comparing the coefficients in front of $\pa_2^0$ in
\ben
\pi_2((\pa_2-q_2)^{-1} H_2\L) = 0
\een
we get
\ben
(\Lambda-1)(c_2) + 2 \pa_2(q_2) = 0.
\een
The formula for $\pa_2(q_2)$ follows. 

Note that $\pi_+(H_a \L)=0$ is equivalent to
\beq\label{pa_2(L)1}
\pa_a(\L)=\pi_+([Q_a,\L]).
\eeq
The derivatives $\pa_2(a_i)$, $\pa_3(a_i)$ ($1\leq i\leq n-4$) and 
$\pa_a(c_b)$ ($2\leq a,b\leq 3$)  can be
determined by comparing the coefficients in front of $\Lambda^i$ for
$0\leq i\leq n-3$. Indeed, let us consider only the case $\pa_2(a_i)$
and $\pa_2(c_b)$, because the other case is analogous. The operator
$\L$ has the form 
\ben
&&
a_{n-3} \Lambda^{n-2} + a_{n-4} \Lambda^{n-3} + \sum_{k=2}^{n-4}
(a_{k-1} - a_{k+1}) \Lambda^{k} +
\Big(-a_2+\frac{1}{4}(c_2-c_3)\Big)\Lambda +\\
&&
\Big(-a_1-a_1[-1] + \frac{1}{2} (c_2+c_3)\Big) \Lambda^0+ \cdots,
\een
where the dots stand for terms involving only negative powers of
$\Lambda$. Let us split $[Q_2,\L]$ into sum
of two commutators $[Q_2,\L-\tfrac{1}{2}(\pa_2^2+\pa_3^2)]$ and
\ben
\frac{1}{2}
[Q_2, \pa_2^2+\pa_3^2] = -\frac{1}{2}(\pa_2^2+\pa_3^2)(Q_2) -
\pa_2(Q_2)\pa_2 -\pa_3(Q_2) \pa_3.
\een
The first commutator is already in $\E^{0}_{(+)}$ and it is a Laurent
series in $\Lambda^{-1}$ whose coefficients are differential
polynomials involving only $\pa_x$-derivatives. The projection $\pi_+$
of the second commutator is
\ben
-\frac{1}{2}(\pa_2^2+\pa_3^2)(Q_2) -
\pa_2(Q_2)Q_2 -\pa_3(Q_2) Q_3 .
\een
A straightforward computation, using formulas \eqref{pa_a(q_b)1} and
\eqref{pa_a(q_a)1}, shows that the above expression has 
leading order term of the type
\ben
&&\Big(\frac{1}{4}\pa_2 (c_2-c_2[-1]) -\frac{1}{2}\frac{1-e^{-\epsilon\pa_x}}{1+e^{\epsilon\pa_x}}\Big(\frac{e^{\epsilon\pa_x}-1}{e^{\epsilon\pa_x}+1}(q_2q_3)\cdot q_3+\frac{1}{2}(c_3-c_3[1])q_2\Big)+\\
&&+\frac{1}{2}(c_2-c_2[-1])q_2[-1]+
\frac{1}{2}(c_3-c_3[-1])q_3[-1]\Big) \Lambda^0 + O(\Lambda^{-1}). 
\een
Comparing the coefficients in front of $\Lambda^k$ in \eqref{pa_2(L)1} for
$1\leq k\leq n-2$ we get that $\pa_2(a_i)$ $(1\leq i\leq n-3)$ and
$\pa_2(c_2-c_3)$ can be expressed as 
differential polynomials that involve only $\pa_x$-derivatives.  
Comparing the coefficients in front of $\Lambda^0$ in \eqref{pa_2(L)1} we get that
\ben
\frac{1}{2}\pa_2(c_2+c_3) - \frac{1}{4}\pa_2 (c_2-c_2[-1])
= \frac{1}{8}(3+e^{-\epsilon\pa_x}) \pa_2(c_2+c_3)-\frac{1}{8}(1-e^{-\epsilon\pa_x})\pa_2(c_2-c_3)
\een
can be expressed in terms of differential polynomials that involve
only $\pa_x$-derivatives. The operator $3+e^{-\epsilon\pa_x}$ is
invertible, so the derivative $\pa_2(c_2+c_3)$ is also a differential
polynomial involving only $\pa_x$-derivatives. Finally, by comparing
the coefficients in front of $\Lambda^{n-2}$ we get 
\ben
\pa_2(a_{n-3}) = a_{n-3} (1-e^{(n-2)\epsilon \pa_x})(q_2[-1])=a_{n-3}
( q_2[-1]-q_2[n-3]) .
\een
Since $a_{n-3}=\tfrac{1}{n-2} e^{(n-2)\alpha}$ the above equation implies
$\pa_2(\alpha) = \tfrac{1}{n-2}(q_2[-1]-q_2[n-3])$. 

Let us define $\pa_2$ and $\pa_3$ as above, i.e., 
\beq\label{der-qc}
\pa_2(q_3)=\pa_3(q_2) =q_1+q_2q_3,\quad
\pa_a (q_a) = \frac{1}{2}(c_a-c_a[1]),\ (a=2,3),
\eeq
and
\ben
\Big(\pa_a(\L)-\pi_+([Q_a,\L])\Big)_{1,\geq 0} = 0.
\een
We have to prove that $H_a \L\in \E H$ for $a=2,3$. We will give the
argument for $a=2$ only, because the case $a=3$ is analogous. Let us first
prove that $\pa_2(\L)-\pi_+([Q_2,\L]=0$. Using formulas \eqref{der-qc}
it is straightforward to verify that
\beq\label{Lax-decomp}
\L=\A + \iota_{\Lambda^{\pm 1}}\C + \iota_{\Lambda^{\pm 1}}\Q,
\eeq
where
\ben
\A=\sum_{i=1}^3 (a_i\Lambda^i-\Lambda^{-i} a_i)(\Lambda-\Lambda^{-1}),
\een
\ben
\C=
\frac{1}{4} (\Lambda-1)^{-1} (\Lambda c_2-c_2\Lambda^{-1}) (\Lambda+1) -
\frac{1}{4} (\Lambda+1)^{-1} (\Lambda c_3-c_3\Lambda^{-1}) (\Lambda-1)
+ \frac{1}{2}  (Q_2^2 + Q_3^2),
\een
and
\ben
\Q= \frac{1}{2} \Big( (\pa_2+Q_2)(\pa_2-Q_2) + (\pa_3+Q_3)(\pa_3-Q_3)
\Big).
\een
Note that the above operators have the following symmetries
\beq\label{Lax-symm}
M^\# =
(\Lambda-\Lambda^{-1})\cdot M\cdot (\Lambda-\Lambda^{-1})^{-1},
\quad M=\A,\C,
\eeq
and
\ben
\Q^\#- (\Lambda-\Lambda^{-1})\cdot \Q\cdot (\Lambda-\Lambda^{-1})^{-1}
= (\Lambda-\Lambda^{-1})\cdot
(\pa_2(Q_2) + \pa_3(Q_3) )\cdot
(\Lambda-\Lambda^{-1})^{-1}.
\een
Using the 0-curvature condition we also have
\ben
\pi_+([Q_2,(\pa_a+Q_a)(\pa_a-Q_a)] =
-Q_a \pa_2(Q_a) -\pa_2(Q_a) Q_a  -\pa_2 \pa_a(Q_a).
\een
Therefore,
\ben
(\pi_+[Q_2,\Q])^\# -  (\Lambda-\Lambda^{-1})[Q_2,\Q]
(\Lambda-\Lambda^{-1})^{-1} =
(\Lambda-\Lambda^{-1})(\pa_2^2(Q_2) + \pa_2\pa_3(Q_3))
 (\Lambda-\Lambda^{-1})^{-1}. 
\een
Using the above symmetries it is strightforward to check that
\ben
\left(  (\Lambda-\Lambda^{-1}) \Big(
  \pa_2(\L)-\pi_+([Q_2,\L]) \Big)\right)^\# +
(\Lambda-\Lambda^{-1}) \Big(
  \pa_2(\L)-\pi_+([Q_2,\L]) \Big)=0.
\een
By definition the operator $(\Lambda-\Lambda^{-1}) \Big(
  \pa_2(\L)-\pi_+([Q_2,\L]) \Big)$ has the form $\sum_{i=0}^\infty
  b_i\Lambda^{-i}$, so the above symmetry implies that $b_i=0$ for all
  $i$.

Note that $\pa_2-Q_2, \Q\in \E_{(+)} H$. Therefore,
\ben
[\pa_2-Q_2,\A+\C] = 
\pi_+( [\pa_2-Q_2,\L]) = \pa_2(\L)-\pi_+([Q_2,\L]) = 0.
\een
Since $H_2=(\Lambda-1) (\pa_2-Q_2)$ the vanishing of the above
commutator  implies that
\ben
H_2\L = (\Lambda-1)\A(\Lambda-1)^{-1}  H_2 + (\Lambda-1)\C
(\Lambda-1)^{-1} H_2 + H_2 \Q. 
\een
The first term on the RHS is already in $\E H$, so we need to verify
that the remaining two terms add up to some element in $\E H$. A long but
straightforward computation using the explicit formulas for $\C,
\Q$, and $H_a$ and formulas \eqref{der-qc} gives that $(\Lambda-1)\C
(\Lambda-1)^{-1} H_2 + H_2 \Q$ is given by the following formula
\ben
\left(
  \frac{1}{2}(\pa_2^2+\pa_3^2) +
  \frac{1}{4}
  \Big(\Lambda (c_2-c_3) + (c_2-c_3)\Lambda^{-1}
  +c_2+c_3+c_2[1]+c_3[1]
  \Big)\right) H_2 + \pa_3(q_2) H_3.
\een
Let us give the explicit formulas for the other case. We have 
\ben
H_3\L = (\Lambda+1)\A(\Lambda+1)^{-1}  H_3 + (\Lambda+1)\C
(\Lambda+1)^{-1} H_3 + H_3 \Q 
\een
and $(\Lambda+1)\C (\Lambda+1)^{-1} H_3 + H_3 \Q $ is given by
\ben
\left(
  \frac{1}{2}(\pa_2^2+\pa_3^2) +
  \frac{1}{4}
  \Big(\Lambda (c_2-c_3) + (c_2-c_3)\Lambda^{-1}
  +c_2+c_3+c_2[1]+c_3[1]
  \Big)\right) H_3 + \pa_2(q_3) H_2.\qed
\een
Let us point out that if the ring $\R$ is equipped with the derivations
$\pa_2$ and $\pa_3$ defined by Proposition \ref{prop:proj-A}, b), then
we have
\ben
H_1=\tfrac{1}{2}(\pa_2-q_2) H_3 - \tfrac{1}{2}(\pa_3-q_3) H_2. 
\een
In particular, $H_1\L\in \E H$, that is, the condition in Proposition
\ref{prop:proj-A}, b) holds also for $a=1$.

\subsection{Rational difference operators}\label{sec:rdo}
Let $\CC(\Lambda)$ be the field of rational functions in $\Lambda$. If
$A$ is any $\CC$-algebra, then we denote by $A(\Lambda) :=A\otimes_\CC
\CC(\Lambda)$. Using elementary fraction decomposition, we get that
the $A$-module $A(\Lambda)$ can be decomposed into a direct sum of $A$-modules as
follows:
\ben
A(\Lambda)= \bigoplus_{n\geq 0} A\otimes\Lambda^n\bigoplus
\bigoplus_{a\in \CC, n\geq 1} A\otimes (\Lambda-a)^{-n}.
\een
We apply this construction for $A:=\CC[\xi^i\ |\ i\geq 0, \xi\in
\Xi]$. The key observation is the following: if $\Lambda$ commutes
with the elements of $A$, then clearly $A(\Lambda)$ is a ring. This is
not the case in our setting, but the commutator $[\Lambda,a]$ is
proportional to $\epsilon$, so the failiure of commutativity may be
offset by allowing infinite power series in $\epsilon$. This idea 
can be realized as follows: Put 
\ben
\E_{\rm rat} := \CC[\xi^i\, (i\geq 0,\ \xi\in \Xi), e^{\pm \alpha}](\Lambda)[\![\epsilon]\!].
\een
We claim that $\E_{\rm rat}$ has a natural ring structure. It is
sufficient to define $(\Lambda-a)^{-1} \cdot P$, where $P\in A$ is a
differential polynomial. In order to justify the definition, note that
we have the following formula:
\beq\label{rat-product}
(\Lambda-a)^{-1} P =\sum_{n=0}^\infty (-\operatorname{ad}_\Lambda)^n(P)
\cdot (\Lambda-a)^{-n-1},
\eeq
where both sides make sense if we expand them in the powers of
$\Lambda^{-1}$, that is, the above equality makes sense in the ring
$A(\!(\Lambda^{-1})\!)[\![\epsilon]\!]$. For the proof, note that we have
\beq\label{rat-prod-step}
(\Lambda-a)^{-1} P = 
P\cdot (\Lambda-a)^{-1} +(\Lambda-a)^{-1} (-\operatorname{ad}_\Lambda(P))\cdot (\Lambda-a)^{-1}.
\eeq
Replacing in this formula $P$ by $(-\operatorname{ad}_\Lambda(P))$, we
get
\ben
(\Lambda-a)^{-1} (-\operatorname{ad}_\Lambda(P)) =
(-\operatorname{ad}_\Lambda(P)) (\Lambda-a)^{-1} +
(\Lambda-a)^{-1} (-\operatorname{ad}_\Lambda)^2(P) (\Lambda-a)^{-1}.
\een
Therefore,
\ben
(\Lambda-a)^{-1} P=P (\Lambda-a)^{-1} +
(-\operatorname{ad}_\Lambda)(P) (\Lambda-a)^{-2} +
(\Lambda-a)^{-1} (-\operatorname{ad}_\Lambda)^2(P) (\Lambda-a)^{-2}.
\een
Clearly, continuing this process by applying formula
\eqref{rat-prod-step} with $P$ replaced by
$(-\operatorname{ad}_\Lambda)^m(P)$, we will get formula
\eqref{rat-product}. On the other hand, if $P\in \E_{\rm rat}$, then
$(-\operatorname{ad}_\Lambda)^m(P)\in \epsilon^m\E_{\rm
  rat}$. Therefore, the RHS of \eqref{rat-product} defines an element
of $\E_{\rm rat}$. In other words, the Laurent series expansion operation
\ben
\iota_{\Lambda^{-1}}: \E_{\rm rat} \to A(\!(\Lambda^{-1})\!)[\![\epsilon]\!]
\een
provides an $A$-module embedding of $\E_{\rm rat}$ and the image is a
subring of $A(\!(\Lambda^{-1})\!)[\![\epsilon]\!]$. In particular, the
$A$-module $\E_{\rm rat}$ has a unique associative product, such that
formula \eqref{rat-product} holds. The elements of the ring $\E_{\rm
  rat}$ will be called {\em rational difference operators}. 

For every $a\in \CC$ let us denote by
\ben
\iota_{\Lambda-a}: \E_{\rm rat} \to A(\!(\Lambda-a)\!)[\![\epsilon]\!]
\een
the Laurent series expansion operations. Let us recall also the rings
$\E_{(\pm)}^0=\R(\!(\Lambda^{\mp 1})\!)$. Note that we have an
embedding 
\ben
\R(\!(\Lambda^{\mp 1})\!)=A[\![\epsilon]\!](\!(\Lambda^{\mp 1})\!)
\subset
A(\!(\Lambda^{\mp 1})\!)[\![\epsilon]\!],
\een
which allows us to think of $\R(\!(\Lambda^{\mp 1})\!)$ as the
elements in $A(\!(\Lambda^{\mp 1})\!)[\![\epsilon]\!]$ that have a
finite order pole at $\Lambda=\infty$ or $0$. Finally, note that we
can identify 
\ben
\R[\![\Lambda-a]\!]=A[\![\epsilon]\!][\![\Lambda-a]\!]= A[\![\Lambda-a]\!][\![\epsilon]\!]
\een
with the elements of $A(\!(\Lambda-a)\!)[\![\epsilon]\!]$ that are
regular at $\Lambda=a$, that is, the series that do not involve
negative powers of $\Lambda-a$. There is a unique associative ring
structure on $\R[\![\Lambda-a]\!]$, such that the multiplication satisfies
\ben
(\Lambda-a) \cdot P  = P[1]\cdot (\Lambda-a) + a \Delta(P),\quad P\in \R,
\een
where $\Delta(P):=P[1]-P$ is the forward difference operator. 

The following properties are straightforward to check:
\begin{itemize}
\item[(i)]
  For every $a\in \CC$,
  there exists a unique associative ring structure on
  $A(\!(\Lambda-a)\!)[\![\epsilon]\!]$, such that 
  formula \eqref{rat-product} holds for $P\in
  A(\!(\Lambda-a)\!)[\![\epsilon]\!]$. 
\item[(ii)] The maps $\iota_{\Lambda^{\pm 1}}$ and $\iota_{\Lambda-a}$
  are injective ring homomorphisms.
\item[(iii)]
  The ring $\R(\!(\Lambda^{\mp 1})\!)$ is a subring of
  $A(\!(\Lambda^{\mp 1})\!)[\![\epsilon]\!]$. 
\item[(iv)] The ring
  $\R[\![\Lambda-a]\!] =A[\![\Lambda-a]\!][\![\epsilon]\!]$
  is a subring of
    $A(\!(\Lambda-a)\!)[\![\epsilon]\!]$.
\end{itemize}

\subsection{The kernel of the projections $\pi_\alpha$}

\begin{lemma}\label{le:inj}
a) Suppose that $P=\sum_{i=0}^m a_i \pa_2^i+ \sum_{i=0}^nb_i \pa_3^i$ is a differential
operator with coefficients $a_i,b_i\in \R$. If $\pi_{\pm }(P)=0$ then
$P=0$.

b) Suppose that $P=\sum_{i=0}^m a_i \pa_2^i+ \sum_{i=-n_1}^{n_2}b_i \Lambda^i$ is a differential
operator with coefficients $a_i,b_i\in \R$. If $\pi_3(P)=0$ then
$P=0$.

c) Suppose that $P=\sum_{i=0}^m a_i \pa_3^i+ \sum_{i=-n_1}^{n_2}b_i \Lambda^i$ is a differential
operator with coefficients $a_i,b_i\in \R$. If $\pi_2(P)=0$ then
$P=0$.

\end{lemma}
\proof
Let us give the argument only for part a) for the case when
$\pi_+(P)=0$. The remaining statements are proved in the same way.

Let us consider first the case when all $b_i=0$. We claim that
$a_0=0$. Indeed, using the rlation $\pa_2=(\Lambda-1)^{-1} q_2
(\Lambda+1) + (\Lambda-1)^{-1} H_2$, we get that $P=G_0+\sum_{i=1}^mG_i H_2^i$,
where the coefficients $G_i\in \E_{\rm rat}$ are rational difference
operators that have a finite order pole at $\Lambda=\infty$ and are
regular at $\Lambda=a$ for all $a\neq 1$. By definition,
$\iota_{\Lambda^{-1}}(G_0)=\pi_+(P)$, which is given to be
$0$. Therefore, $G_0=0$ and we get $P=GH_2$ for some $G\in \E_{\rm
  rat}[\pa_2]$. Since $P$ is independent of $\Lambda$, we get
$P= \iota_{\Lambda+1}(G)H_2$ in the ring
$\R[\pa_2][\![\Lambda+1]\!]$. Comparing the coefficients in front of
$\pa_2^0 (\Lambda+1)^0$ we get that $a_0=0$.

Let $m=\operatorname{max}\{i: a_i\neq 0\}$ be the order of the
differential operator $P$. We argue by induction on $m$ that 
$a_i= 0$ for all $i$. Suppose that $\operatorname{ker}(\pi_+)$ does
not contain differential operators of order $\leq m-1$. Note that the
operator $\widetilde{P}:=(a_m \pa_3 - \pa_3(a_m)) 
\cdot P$ is also in the kernel of $\pi_+$. However,
\ben
\widetilde{P} = \sum_{i=1}^{m-1} (a_m \pa_3(a_i)-\pa_3(a_m) a_i)
\pa_2^i + \sum_{i=1}^{m-1} a_m a_i \pa_2^{i-1} (-q_1+H_1)
\een
and if we set $H_1=0$ in the above formula then we still have an
operator whose projection $\pi_+$ is 0 and whose order is at most
$m-1$.  The inductive assumption implies that
\ben
\widetilde{P} = \sum_{i=1}^{m-1} (a_m \pa_3(a_i)-\pa_3(a_m) a_i)
\pa_2^i + \sum_{i=1}^{m-1} a_m a_i \pa_2^{i-1} (-q_1) = 0.
\een
Comparing the coefficients in front of $\pa_2^j$ for $j=0,1,\dots,
m-2$ we get that $a_i=0$ for all $i=1,2,\dots, m-1$. This proves that
$\widetilde{P}=0$ and since the ring $\R[\pa_2]$ does not have zero divisors
we get $P=0$.

Similar argument works if $a_i=0$ for all $i$. Let
$m:=\operatorname{max}\{i: a_i\neq 0\}$
and
$n:=\operatorname{max}\{i: b_i\neq 0\}$. Then it remains only to
consider the case $m>0$ and $n>0$. We argue by induction on $n$ that
such a case is impossible. Note
that the operator $(b_n\pa_2-\pa_2(b_n)) \cdot P$ has the form
\ben
\sum_{i=0}^m (b_n\pa_2-\pa_2(b_n)) a_i \pa_2^i + 
\sum_{i=0}^{n-1} (b_n \pa_2(b_i)-\pa_2(b_n) b_i) \pa_3^i +
b_nb_0 \pa_2+
\sum_{i=1}^n b_n b_i \pa_3^{i-1} (-q_1+ H_1).
\een 
If we set $H_1=0$ in the above expression we will get an operator
whose projection $\pi_+$ is also $0$. The coefficient in front of
$\pa_2^{m+1}$ is $a_m b_n\neq 0$ while the highest possible power of
$\pa_3$ is $n-1$ -- contradiction with our inductive assumption. 
\qed

\begin{proposition}\label{prop:ker_pi_a}
The following equality of left ideals holds $\E_{(\alpha)}H\cap \E= \E
H$ for all $ \alpha\in \{\pm,2,3\}$.
\end{proposition}
\proof
The arguments in all 4 cases are similar, so let us consider only the
case $\alpha=+$. Suppose $P\in \E_{(+)}H\cap \E$. We may assume that
$P$ is polynomial in $\Lambda$. Indeed, if not then we can always find
a positive integer $n$ such that $\Lambda^n P\in
\R[\pa_2,\pa_3,\Lambda]$. The operator $\Lambda^n P\in \E_{(+)}H\cap
\E$, so if we knew that $\Lambda^n P \in \E H$, then $P\in \E H$. 

Suppose that $P$ is polynomial in $\Lambda$. We will reduce the proof
to the case when $P=P_1+P_2+P_3$, where $P_1\in \R[\Lambda]$, $P_a\in
\R[\pa_a]$ $(a=2,3)$. Let us  write 
\ben
P=\sum_{i,j,k} p_{i,j} (\Lambda)\pa_2^i \pa_3^j,\quad
p_{i,j}(\Lambda)\in \R[\Lambda],
\een
where the sum in $i$ and $j$ is finite. Using that $\pa_2 \pa_3= H_1-q_1$ we can write $P$ as a sum of an
element in $\E H$ and $P_2+P_3$, where $P_a\in \R[\Lambda,\pa_a]$
$(a=2,3)$.  Let us write $P_2$ (resp. $P_3$) as a sum of monomials of
the type $(\Lambda-1)^i \pa_2^j$ (resp. $(\Lambda+1)^i
\pa_3^j$). Using the relation $(\Lambda-1)\pa_2 = H_2 + q_2
(\Lambda+1)$ (resp.  $(\Lambda+1)\pa_3 = H_3 + q_3
(\Lambda-1)$) we can reduce $P_a$ to a sum of a polynomial in
$\R[\Lambda]$ and a polynomial in $\R[\pa_a]$. 

Suppose now that 
$P=P_1+P_2+P_3$, where $P_1\in \R[\Lambda]\Lambda$, $P_a\in  \R[\pa_a]$ $(a=2,3)$.  
We are given that $0=\pi_+(P)=P_1 + \pi_+(P_2) + \pi_+(P_3)$. The
projections $\pi_+(P_a)$ ($a=2,3$) have only non-positive powers of
$\Lambda$. Therefore, by comparing the coefficients in front of the
postive powers of $\Lambda$, we get $P_1=0$. Therefore,
$\pi_+(P_2+P_3)=0$. Recalling Lemma \ref{le:inj}, a) we get
$P_2+P_3=0$.
\qed

\subsection{Residue formulas for the projections}

\begin{proposition}\label{prop:proj_pi}
  a) If $k>0$ then the following formulas hold
  \ben
  \pi_\pm(\pa_2^k) & = & \tfrac{1}{2}\iota_{\Lambda^{\mp 1}}
  \Big( \pa_2^k H_2^{-1}(\pa_2+q_2) \Big)_{2,[0]} (\Lambda-\Lambda^{-1}),
  \\
  \pi_3(\pa_2^k) & = & 
  \Big( \pa_2^k H_1^{-1}\, \pa_2\,  \Big)_{2,[0]} \pa_3.
  \een
  b) If $k>0$ then the following formulas hold
  \ben
  \pi_\pm(\pa_3^k) & = & \tfrac{1}{2}\iota_{\Lambda^{\mp 1}}
  \Big( \pa_2^k H_3^{-1}(\pa_3+q_3) \Big)_{3,[0]} (\Lambda-\Lambda^{-1}),
  \\
  \pi_2(\pa_3^k) & = & 
  \Big( \pa_3^k H_1^{-1}\, \pa_3\,  \Big)_{3,[0]} \pa_2.
  \een
c) If $k>0$ then the following formulas hold
\ben
\pi_2(\Lambda^{\pm k})  & = &  (-1)^k \pm 2\Big(\Lambda^{\pm k}
\iota_{\Lambda^{\mp 1}} \left(
H_2^{-1}(\Lambda^{-1}+1)^{-1} \right)\Big)_{1,[0]} \cdot \pa_2,\\
\pi_3(\Lambda^{\pm k})  & = &  1 \pm 2\Big(\Lambda^{\pm k}
\iota_{\Lambda^{\mp 1}} \left(
H_3^{-1}(\Lambda^{-1}-1)^{-1} \right)\Big)_{1,[0]} \cdot \pa_3.
\een
\end{proposition}
\proof
Let us prove the first formula in a)
\begin{eqnarray}
\pi_{+}(\pa_a^k)=\frac{1}{2}\,\,\iota_{\Lambda^{-1}} \left(
\pa_a^k\, H_2^{-1} (\pa_a+q_a)
\right)_{2,[0]}\cdot (\Lambda-\Lambda^{-1}).\label{pipluspaa}
\end{eqnarray}
The remaining formulas are proved with the same method. Put
\beq\label{h2-expansion}
H_2^{-1}\cdot(\pa_2+q_2)=(\Lambda-1)^{-1}+\sum_{j=1}^{+\infty}\pa_2^{-j}\cdot v_j.
\eeq
Using the formula
\ben
f\cdot \pa_2^{-l}=
\sum_{j=0}^\infty\binom{l+j-1}{j}\pa_2^{-l-j}\cdot \pa_2^j(f)
\een
and comparing the coefficients on the right of $\pa_2^{-l}$ in
\ben
(\pa_2+q_2)=((\Lambda-1)\pa_2-q_2(\Lambda+1)) \, \Big(
(\Lambda-1)^{-1}+\sum_{j=1}^{+\infty}\pa_2^{-j}\cdot v_j
\Big) 
\een
we get 
\beq\label{vl-recursion}
v_1=2 Q_2(\Lambda-\Lambda^{-1}),\
v_{l+1}=\sum_{j=0}^{l-1}\binom{l-1}{j}\pa_2^j(Q_2)\cdot v_{l-j},\quad
l=1,2,3,\cdots, 
\eeq
where recall that $Q_2=(\Lambda-1)^{-1} q_2 (\Lambda+1)$. 
Furthermore, note that 
\beq\label{vj-relation}
\sum_{l=1}^\infty\pa_2^{-l}\cdot v_{l+1}=Q_2\cdot\sum_{l=1}^\infty \pa_2^{-l}\cdot v_l.
\eeq
Indeed, the formula follows by comparing the pseudo-differential part,
i.e., the terms involving only negative powers of $\pa_2$ in
\ben
\Big(\pa_2-Q_2\Big)\cdot \Big((\Lambda-1)^{-1}+
\sum_{j=1}^\infty \pa_2^{-j}\cdot v_j\Big)=
(\Lambda-1)^{-1}\cdot (\pa_2+q_2).
\een
Recalling the definition of $\pi_+$ we get
\ben
\pi_+(\pa_2)=Q_2,\quad \pi_+(\pa_2^{l+1})=\pa_2(\pi_+(\pi_2^l))+\pi_+(\pa_2^l)\cdot Q_2,\quad l=1,2,3,\cdots.
\een
It is straightforward to check the relation (\ref{pipluspaa}) is
correct for $l=1.$ Let us denote by $A_l$ the RHS of
(\ref{pipluspaa}). We prove (\ref{pipluspaa}) by induction on
$k$. Clearly, we need only to show that
\begin{eqnarray}
A_{l+1}=\pa_2(A_l)+A_l\cdot Q_2.\label{recursional}
\end{eqnarray}
Substituting the expansion \eqref{h2-expansion} in $A_l$ we get 
\begin{eqnarray}
A_l=\frac{1}{2}\sum_{j=1}^l\pa_2^{-j+l}(v_j)\cdot (\Lambda-\Lambda^{-1}).
\end{eqnarray}
After inserting the above formula for $A_l$ in the recursion relation
(\ref{recursional}), we get that (\ref{recursional}) is equivalent to
\begin{eqnarray}
v_{l+1}=-\sum_{j=1}^l\pa_2^{-j+l}(v_j)\cdot Q_2^{\#},
\end{eqnarray}
which according to formula \eqref{vl-recursion} is equivalent to
\begin{eqnarray}
  -\sum_{j=0}^l\pa_2^{-j+l}(v_j)\cdot Q_2^{\#}=
  \sum_{j=0}^{l-1}\binom{l-1}{j}\pa_2^j(Q_2)\cdot v_{l-j}.\label{sumvjq2}
\end{eqnarray}
Furthermore, let us rewrite (\ref{sumvjq2}) as
\begin{eqnarray}
-\operatorname{Res}_{\pa_2}\Big(\pa_2^l\cdot H_2^{-1}\cdot(\pa_2+q_2)\cdot\pa_2^{-1}\cdot Q_2^{\#}\Big)=
\operatorname{Res}_{\pa_2}\Big(\pa_2^{l-1}\cdot Q_2\cdot H_2^{-1}\cdot(\pa_2+q_2)\Big),\quad l\geq 1,
\end{eqnarray}
which is equivalent to
\begin{eqnarray}
-\Big(\pa_2\cdot H_2^{-1}\cdot(\pa_2+q_2)\cdot\pa_2^{-1}\cdot
  Q_2^{\#}\Big)_{2,<0}=
  \Big(Q_2\cdot H_2^{-1}\cdot(\pa_2+q_2)\Big)_{2,<0}.
\end{eqnarray}
Using the expansion \eqref{h2-expansion} we transform the above relation into
\begin{eqnarray}
  &&
     -\pa_2\cdot\sum_{l=1}^\infty\pa_2^{-l}\cdot v_l\cdot
     \pa_2^{-1}\cdot Q_2^{\#}=
     Q_2\cdot\sum_{l=1}^\infty\pa_2^{-l}\cdot v_l \nonumber\\
  &\Leftrightarrow&
                    \sum_{l=1}^\infty \pa_2^{-l+1}\cdot
                    v_l\cdot(\Lambda-\Lambda^{-1})\cdot\pa_2^{-1}\cdot
                    Q_2=
                    Q_2\cdot\sum_{l=1}^\infty\pa_2^{-l}\cdot v_l\cdot(\Lambda-\Lambda^{-1})\nonumber\\
  &&
     \Big(\text{by $Q_2^{\#}=-(\Lambda-\Lambda^{-1})\cdot
     Q_2\cdot(\Lambda-\Lambda^{-1})^{-1}$}\Big)
  \nonumber\\
  &\Leftrightarrow&
                    v_1\cdot
                    (\Lambda-\Lambda^{-1})\cdot\pa_2^{-1}\cdot
                    Q_2+Q_2\cdot\sum_{l=1}^\infty \pa_2^{-l}\cdot
                    v_l\cdot(\Lambda-\Lambda^{-1})\cdot\pa_2^{-1}\cdot
                    Q_2\nonumber\\ 
&&=Q_2\cdot\sum_{l=1}^\infty\pa_2^{-l}\cdot
   v_l\cdot(\Lambda-\Lambda^{-1}),\quad
   \Big(\text{by formula \eqref{vj-relation}}\Big) \nonumber\\
  &\Leftrightarrow&
                    2\pa_2^{-1}\cdot Q_2+
                    \sum_{l=1}^\infty \pa_2^{-l}\cdot
                    v_l\cdot(\Lambda-\Lambda^{-1})\cdot\pa_2^{-1}\cdot
                    Q_2=
                    \sum_{l=1}^\infty\pa_2^{-l}\cdot v_l\cdot(\Lambda-\Lambda^{-1})\nonumber\\
&&\Big(\text{by $v_1=2Q_2\cdot(\Lambda-\Lambda^{-1})^{-1}$}\Big) \nonumber\\
  &\Leftrightarrow&
                    2\pa_2^{-1}\cdot Q_2(1-\pa_2^{-1}\cdot Q_2)^{-1}=
                    \sum_{l=1}^\infty \pa_2^{-l}\cdot v_l\cdot (\Lambda-\Lambda^{-1})\nonumber\\
  &&
     \Big(\text{by solving for
     $\sum_{l=1}^\infty \pa_2^{-l}\cdot v_l\cdot
     (\Lambda-\Lambda^{-1})$}
     \Big) \nonumber\\
  &\Leftrightarrow&
                    2(1-\pa_2^{-1}\cdot Q_2)^{-1}-1+\Lambda^{-1}=
                    H_2^{-1}\cdot(\pa_2+q_2)\cdot(\Lambda-\Lambda^{-1})\quad\Big(\text{by formula \eqref{h2-expansion}}\Big) \nonumber\\
  &\Leftrightarrow&
                    H_2=(\Lambda-1)\cdot\pa_2-q_2\cdot(\Lambda+1),\quad
                    \Big(\text{by $1-\pa_2^{-1}\cdot Q_2=
                    \pa_2^{-1}\cdot(\Lambda-1)^{-1}\cdot H_2$}\Big), \nonumber
\end{eqnarray}
which is just the definition of $H_2$. \qed

\section{Existence of the derivations}

In this section we are going to prove part a) of Theorem
\ref{thm:Lax}, except for the existence of the derivations of type
$\pa_{0,k}$ ($k\geq 1$). The latter are called {\em extended
  flows}. Their construction, based on dressing operators, requires
some additional work which we postpone until next section.

\subsection{Auxiliarly lemma}
The following simple lemma allows us to construct
derivations of $\R$ commuting with $\pa_x$. 
\begin{lemma}\label{le:der_R}
Suppose that $B\in \E$ is a differential-difference operator, s.t.,
the projections $B^\pm:=\pi_{\pm}(B)$ satisfy the following
conditions: 
\begin{itemize}
\item[(i)]
  $(B^+)^\# = -(\Lambda-\Lambda^{-1})
  B^- \iota_\Lambda (\Lambda-\Lambda^{-1})^{-1}$.
\item[(ii)]
  $\pi_+[B^+,\L]=\pi_-[B^-,\L]\in \E$.
\item[(iii)]
  $\pi_+(H_a B^+)=\pi_-(H_a B^-)\in \E$.
\end{itemize}
Then there exists a unique derivation $\pa$ of $\R$ commuting with
$\pa_x$, such that one of the following two equivalent conditions is satisfied
\begin{align}\label{L-deriv}
  \pa(\L) & =   \pi_+([B^+,\L]), \\
  \label{H_a-deriv}
\pa (H_a) & = -\pi_+(H_a B^+) ,\quad a=2,3,
\end{align}
or
\begin{align}\nonumber
  \pa(\L)-[B,\L] & \in \E H,\\
  \nonumber
  \pa(H_a)-[B,H_a] & \in \E H, \quad a=2,3.
\end{align}
Moreover $\pa\pa_a(q_a)=\pa_a \pa(q_a)$ for $a=2,3$. 
\end{lemma}
\proof
We claim that it is sufficient to prove that there exists a unique
derivation satisfying \eqref{L-deriv} and \eqref{H_a-deriv}. Indeed, if
we knew this, then let us write $B^+= B +\sum_{b=2,3} P_b H_b$, where
$P_b\in \E_{(+)}$. We have
\begin{align}\label{commL-proj}
  [B^+,\L] & =[B,\L] + \sum_{b=2,3} ([P_b,\L] H_b + P_b[\L,H_b]),\\
  \label{commH-proj}
  [B^+,H_a] & =[B,H_a] + \sum_{b=2,3} ([P_b,H_a] H_b + P_b[\L,H_b]).
\end{align}
Since $H_b\L\in \E H$, we get that $\pi_+([B^+,\L])=\pi_+([B,\L]).$
Therefore, $\pa(\L)-[B,\L]\in \E_{(+)}H\cap \E=\E H$. Similarly,
$\pa(H_a)-[B,H_a]\in \E_{(+)}H\cap \E=\E H$. This proves that there
exists a derivation satisfying the scond set of equations. Conversely,
suppose that there exists a derivation $\pa$ satisfying the scond set
of equations. Then the relations \eqref{commL-proj} and
\eqref{commH-proj} prove that the 1st set of equations \eqref{L-deriv}
and \eqref{H_a-deriv} are also satisfied
and hence $\pa$ must be the unique derivation satisfying the firs set
of equations. Our claim is proved.

It remains to prove that there exists a unique derivation satisfying
\eqref{L-deriv} and \eqref{H_a-deriv}.  Note that $\pa$ is uniquely
defined by the equations 
\ben
\left(\pa_{1,k}(\L) -  \pi_+([B_{1,k}^+,\L])\right)_{1,\geq 0} & = & 0, \\
\left(\pa_{i,k} (H_a) +\pi_+(H_a B_{1,k}^+) \right)_{1,[0]}  & = & 0 ,\quad a=2,3
\een
and the requirement that it commutes with $\pa_x$. Indeed, the first
equation determines $\pa(a_i), \pa(c_2),$ and
$\pa(c_3)$, while the remaining two equations determine
$\pa(q_a)$ $(a=2,3)$. We would like to prove that
\ben
\pa(\L) & = &  \pi_+([B^+,\L]), \\
\pa (H_a) & = & -\pi_+(H_a B^+) ,\quad a=2,3.
\een
First we will prove the second equation. Let us
consider only the case $a=2$. The argument for $a=3$ is similar.
Put
\ben
\A:=(\Lambda^{-1}+1) (\pa(H_2) + \pi_+(H_2 B^+)) =
(\Lambda^{-1}+1) \pa(H_2) + (\Lambda-\Lambda^{-1}) \Big(
  \pa_2(B^+) + [B^+,Q_2] \Big).  
  \een
Using $H_2^\#= (\Lambda^{-1}+1) H_2 (\Lambda+1)^{-1}$ and conditions
(i) and (iii) we get that $  \A^\#=\A$. Recalling again condition
(iii) we write
\ben
\pa(H_2) + \pi_+(H_2 B^+) = a_1 \Lambda + a_0 + \sum_{i=1}^m a_{-i} \Lambda^{-i}.
\een
Then
\ben
\A = a_1 \Lambda + (a_0+a_1[-1]) + \sum_{i=1}^{m+1} (a_{-i}+a_{-i+1}[-1])
\Lambda^{-i},
\een
where $a_{-m-1}:=0$. Since $\A^\#=\A$ we get the following $m+1$
equations: $a_{-i}+a_{-i+1}[-1]=0$ for $2\leq
i\leq m+1$ and $a_1=a_{-1}[1]+a_0$. The first $m$ equations imply
that $a_{-i}=0$ for all $1\leq i\leq m$. Furthermore, $a_0=0$
by the definition of $\pa(H_2)$. Finally, we get $a_1=0$ from the last
equation. 

Let us prove the first equation. We will use the following identity
\beqa\label{Lax-invol}
\left(\pa(\L)-\pi_+[B^+,\L]\right)^\# -
(\Lambda-\Lambda^{-1})\,
\left(\pa(\L)-\pi_+[B^+,\L]\right)\,
(\Lambda-\Lambda^{-1})^{-1}  & =  \\
\notag
\sum_{a=2,3} (\Lambda-\Lambda^{-1})\,
\left(
  \pa\pa_a(Q_a) - \pa_a \pa(Q_a)
\right)\,
(\Lambda-\Lambda^{-1})^{-1}, & 
\eeqa
where on both sides we could take either the expansion $\iota_\Lambda$ or
the expansion $\iota_{\Lambda^{-1}}$ and the identity holds in both
cases. Let us prove \eqref{Lax-invol}. The main difficulty is to compute
\beq\label{comm-invol}
\Big(\pi_+([B^+,\L])\Big)^\# =
[B^+,M]^\# -
\sum_{a=2,3}\Big(
\pa_a(B^+)Q_a +\tfrac{1}{2} \pa_a^2(B^+)
\Big)^\#,
\eeq
where $M:=\L-\tfrac{1}{2}(\pa_2^2+\pa_3^2)$.
It is straightforward to check that
\beq\label{L-invol}
\L^\# = (\Lambda-\Lambda^{-1}) \, \left( \L +
  \pa_2(Q_2)+\pa_3(Q_3)\right)\, (\Lambda-\Lambda^{-1})^{-1} ,
\eeq
where on the RHS we take either the expansion in the powers of
$\Lambda^{-1}$ or the expansion in $\Lambda$. 
Recalling conditions (i) and (ii), after a short
computation, we get that \eqref{comm-invol} transforms into
\ben
(\Lambda-\Lambda^{-1})\, \left(
\pi_+[B^+, \L]
+\sum_{a=2,3} \Big(
\pa_a^2(B^-)+[\pa_a(B^-),Q_a]+[B^-,\pa_a(Q_a)]\Big)
\right)\, (\Lambda-\Lambda^{-1})^{-1}.
\een
On the other hand, 
\ben
\pa_a\pa (Q_a) = -\pa_a \pa (\pa_a-Q_a) = \pa_a
\pi_-((\pa_a-Q_a) B^-) = \pa_a^2(B^-)+\pa_a[B^-,Q_a].
\een
Therefore,
\ben
\Big(\pi_+([B^+,\L])\Big)^\# =
(\Lambda-\Lambda^{-1})\, \left( \pi_+[B^+,\L]
  +\pa_a\pa(Q_a) \right) \, (\Lambda-\Lambda^{-1})^{-1}.
\een
The above identity and \eqref{L-invol} yields \eqref{Lax-invol}.

Put $A:=
(\Lambda-\Lambda^{-1})(\pa(\L)-\pi_+[B^+,\L])$. Then
\eqref{Lax-invol} yields
\beq\label{Lax-Lax_invol}
-A^\# - A  = (\Delta_2+\Delta_3) \Lambda + (
\Delta_2+\Delta_3-\Delta_2[-1]-\Delta_3[-1]) + \Lambda^{-1} (\Delta_2+\Delta_3),
\eeq
where $\Delta_a:=\pa\pa_a(q_a)-\pa_a\pa(q_a)$. The
definition of $\pa$ and condition (ii) imply that $\pa(\L)-\pi_+[B^+,\L]=\sum_{i\geq 1} a_i
\Lambda^{-i}$ has only negative powers of $\Lambda$ and that only finitely many
$a_i$'s are not zero. Suppose that there exists at least 
one $i\geq 1$ such that $a_i\neq 0$. Let us denote by $m$ the largest
such $i$. Comparing the coefficients in front
$\Lambda^{-m-1}$ in \eqref{Lax-Lax_invol} we get $a_m=0$ --
contradiction. Therefore $A=0$, i.e.,
$\pa(\L)-\pi_+[B^+,\L]=0$ which is exactly what we wanted
to prove. Note that as a biproduct of our argument we get also that 
$\pa \pa_a(q_a)=\pa_a\pa(q_a)$. \qed

\subsection{The derivation $\pa_{1,k}$}\label{sec:pa_1k}
Let us first prove the existence of the Lax operator $L_1$.
\begin{lemma}\label{le:L1}
  There exists a unique operator $L_1=b_{1,0}\Lambda
  +\sum_{i=1}^\infty b_{1,i}\Lambda^{1-i}$ with coefficients $b_i\in \R$, such that,
\ben
b_{1,0}= \exp\Big( \frac{(\Lambda-1)(n-2)}{\Lambda^{n-2}-1} (\alpha) \Big),
\quad
\pi_+(\L)=\frac{1}{n-2} \, L_1^{n-2}.
\een
\end{lemma}
\proof
In order to avoid cumbersome notation let us write $b_i$ for
$b_{1,i}$. Let us first check that the coefficients in front of
$\Lambda^{n-2}$ in $\pi_+(\L)=\tfrac{1}{n-2} \, L_1^{n-2} $ match,
i.e., the following identity must hold: 
\ben
a_{n-3} = \frac{1}{n-2}b_0b_0[1]\cdots b_0[n-3].
\een
On the other hand, $a_{n-3}=\tfrac{1}{n-2}e^{(n-2)\alpha}$ and 
\ben
\frac{(\Lambda-1)(n-2)}{\Lambda^{n-2}-1} (1+\Lambda+\cdots
+\Lambda^{n-3})(\alpha) = (n-2)\alpha.
\een
Our claim is proved. 

Comparing the coefficients in front of the powers of
$\Lambda^{n-2-i}$ for $i\geq 1$ we get a relation of the following form
\ben
\Big(\sum_{s=0}^{n-3} (b_0\cdots b_0[s-1])(b_0[s+1-i]\cdots b_0[n-3-i])
\Lambda^s \Big) (b_i) = \cdots,
\een
where on the RHS we have an expression involving $b_0,\dots,b_{i-1}$
and the coefficients of $\pi_+(\L)$. The operator on the LHS that is
applied to $b_i$ belongs to $\R[\pa_x][\![\epsilon]\!]$, i.e., it is a
formal power series in $\epsilon$ whose 
coefficients are differential operators in $\pa_x$ with coefficients
in $\R$. The leading order term is the constant $n-2$ because $b_0=
1+O(\epsilon)$ and $\Lambda=1+O(\epsilon)$, so the operator is invertible as an element in
$\R[\pa_x][\![\epsilon]\!]$.  Therefore we can solve for $b_i$, get
a recursion that uniquely determines $L_1$, and using a simple
induction on $i$ prove that $b_i\in \R$.\qed

Let us recall the difference operators $B_{1,k}$ ($k\geq 1$) defined
in Section \ref{sec:Lax_eqns}.  
We will check that $B^+:=B_{1,k}$ and
$B^-:=B_{1,k}$ satisfy conditions (i)--(iii) in Lemma
\ref{le:der_R}. Note that condition (i) is immediate from the
definition of $B_{1,k}^+$. We have to verify only (ii) and (iii). 

To begin with, note that the decomposition \eqref{Lax-decomp} and the symmetry
\eqref{Lax-symm} implies that 
\beq\label{proj_pm-L}
\L=\tfrac{1}{n-2}\, L_1^{n-2} + \iota_{\Lambda^{-1}} \Q =
\tfrac{1}{n-2}\, (L^-_1)^{n-2} + \iota_{\Lambda} \Q, 
\eeq
where
\ben
\Q= \tfrac{1}{2}\Big( (\pa_2+Q_2)(\pa_2-Q_2) +
(\pa_3+Q_3)(\pa_3-Q_3)\Big)
\een
and
\beq\label{L1-}
L_1^- := \iota_{\Lambda} (\Lambda-\Lambda^{-1})^{-1}\,
L_1^\# \, (\Lambda-\Lambda^{-1}). 
\eeq
In particular,
$\pi_-(\L)=\tfrac{1}{n-2}(L_1^-)^{n-2}$.

We have by definition that 
\ben
\pi_+(H_2 B_{1,k})  & = &  (\Lambda-1) \left(
  \pa_2(B_{1,k}) +[B_{1,k},\iota_{\Lambda^{-1}} Q_2]
  \right) \\
  \pi_-(H_2 B_{1,k})  & = &  (\Lambda-1) \left(
    \pa_2(B_{1,k}) +[B_{1,k},\iota_{\Lambda} Q_2]
    \right). 
\een
Subtracting the second equation from the first one we get
\ben
\pi_+(H_2 B_{1,k})- \pi_-(H_2 B_{1,k}) =
(\Lambda-1)\, B_{1,k} \, \left(\sum_{m\in \ZZ} \Lambda^m \right)\,  H_2=0
\een
where we used that $B_{1,k}$
is divisible from the right by $(\Lambda-\Lambda^{-1})$. We claim that
$[\pa_2-Q_2,L_1]=0$. Indeed, using \eqref{proj_pm-L} we get that
$[\pa_2-Q_2,L_1^{n-2}] = \pi_+[\pa_2-Q_2,L_1^{n-2}]=0$. Suppose that
that $[\pa_2-Q_2,L_1]=a_m \Lambda^m + O(\Lambda^{m-1})$ with $a_m\neq
0$, then the vanishing of the coefficients of the highest power  of
$\Lambda$ in 
$[\pa_2-Q_2,L_1^{n-2}] $
yields
\ben
\sum_{i=1}^{n-2} b_{1,0} b_{1,0}[1]\cdots b_{1,0}[i-1] a_m[i] b_{1,0}[m+i]\cdots b_{1,0}[m+n-3]=0. 
\een
Since $b_{1,0}=1+O(\epsilon)$ the above equation implies that $a_m=0$
-- contradiction. In particular, our claim implies that $\pi_+(H_2
L_1^k)=0$. Therefore,  $\pi_+(H_2B_{1,k}) =  \pi_+(H_2 (B_{1,k}-L_1^k))$
is a Laurent series in $\Lambda^{-1}$ that involves only powers of
$\Lambda$  that are $\leq 1$. This completes the proof of condition (iii). 

Let us prove that condition (ii) holds. We have
\ben
\pi_+[B_{1,k},\L] = [B_{1,k},M] + \tfrac{1}{2} \pi_+[B_{1,k},\pa_2^2+\pa_3^2],
\een
where $M:=\L-\tfrac{1}{2}(\pa_2^2+\pa_3^2)$. Similarly,
\ben
\pi_-[B_{1,k},\L] = [B_{1,k},M] + \tfrac{1}{2} \pi_-[B_{1,k},\pa_2^2+\pa_3^2].
\een
It remains only to notice that 
\ben
\tfrac{1}{2} 
\pi_+[B_{1,k},\pa_a^2] = -\pa_a(B_{1,k})\iota_{\Lambda^{-1}}Q_a -\tfrac{1}{2}
\pa_a^2(B_{1,k}) =
\tfrac{1}{2} 
\pi_-[B_{1,k},\pa_a^2],
\een
where we used that $B_{1,k}$ is divisible from the right by
$\Lambda-\Lambda^{-1}$, so $\pa_a(B_{1,k})Q_a\in
\R[\Lambda,\Lambda^{-1}]$ is a difference operator. Finally, note that 
$\pi_+[B_{1,k},\L] = \pi_+[B_{1,k} - L_1^k,\L]\in
\Lambda^{n-2}\R[\Lambda^{-1}]$ because $B_{1,k}-L_1^k$ does not
involve terms with positive powers of $\Lambda$.

\subsection{The derivation $\pa_{a,2l+1}$}
Let us prove the existence of the operators $L_a$.
\begin{lemma}\label{le:La}
  There exists a unique pseudo-differential operator
  $L_a=\pa_a+\sum_{i=1}^\infty b_{a,i} \pa_a^{1-i}$ ($a=2,3$) with coefficients
  in $\R$, such that,
  \ben
  b_{a,1}=0,\quad \pi_a(\L) = \frac{1}{2} L_a^2.
  \een
\end{lemma}
\proof
The proof is straightforward computation by comparing the coefficients
in front of $\pa_a^2$, $\pa_a$, and $\pa_a^{1-i}$ ($i\geq 1$). The
first identity is trivially satisfied. The second one is equivalent to
$b_{a,1}=0$, and the remaining ones, give a recursion for solving
$b_{1,i+1}$ in terms of $b_{1,1},b_{1,2},\dots,b_{1,i}$.
\qed

Let us recall the differential operators $B_{a,2l+1}$ from Section
\ref{sec:Lax_eqns}. Put 
\ben
B^\pm_{a,2l+1} := \frac{1}{2}\iota_{\Lambda^{\mp 1}} 
\Big( L_a^{2l+1} H_a^{-1} (\pa_a+q_a)\Big)_{a,[0]}\,
(\Lambda-\Lambda^{-1}),
\quad a=2,3,
\een
where $H_a^{-1}$ is the inverse of $H_a$ in the ring
$\R(\!(\Lambda)\!)(\!(\pa_a^{-1})\!)$.  Recalling Proposition
\ref{prop:proj_pi}, we get $B^\pm_{a,2l+1}=\pi_{\pm}(B_{a,2l+1})$.
\begin{lemma}\label{le:pa_2-identities} Suppose that $a\in
  \{2,3\}$. Then the following properties are satisfied.
  
  a) The operators $\pa_a^{-1} H_1$, $(\pa_a+q_a)^{-1}H_a $,  and
  $L_a$ pairwise commute.

  b) $L_a^\# = -\pa_a \cdot L_a \cdot \pa_a^{-1}$ and $(L_a^{2l+1})_{a,[0]}=0$.

  c) $\pi_+(H_b B^+_{a,2l+1}) = \pi_-(H_b B^-_{a,2l+1})$ ($b=1,2,3$).

  d) $\pi_+([B_{a,2l+1}^+,\L]) = \pi_-([B_{a,2l+1}^-,\L])$.
  
\end{lemma}
\proof
a) Let us argue in the case  when $a=2$. The case $a=3$ is similar.
We have already proved that
$\widetilde{H}_1 = \pa_2^{-1} H_1 
= \pa_3 +\pa_2^{-1} q_1$ and 
\ben
\widetilde{H}_2:=(\pa_2+q_2)^{-1}H_2 = (\pa_2+q_2)^{-1}(\pa_2-q_2)
\Lambda-1
\een
commute (see Proposition \ref{prop:proj-A}). Let us prove that
$[\widetilde{H}_2,L_2]=0$. The argument for the remaining identity
$[\widetilde{H}_1,L_2]=0$ is similar.  
We have
\ben
[\widetilde{H}_2,L_2^2] = \Big(
(\pa_2+q_2)^{-1}(\pa_2-q_2) (L_2[1])^2- L_2^2
(\pa_2+q_2)^{-1}(\pa_2-q_2)
\Big) \Lambda.
\een
By definition $L_2^2=2 \pi_2(\L)$ and we know that $\widetilde{H}_2
\L\in \E_{(2)}H$. Therefore, $\pi_2([\widetilde{H}_2,L_2^2])=0$. On
the other hand, if $\pi_2(A\Lambda)=0$ for some $A\in \E^0_{(2)}$,
then $A=0$. Indeed, by definition
\ben
\pi_2(A\Lambda) = A(\pa_2-q_2)^{-1} (\pa_2+q_2)
\een
and the operator $(\pa_2-q_2)^{-1} (\pa_2+q_2)$ is invertible. This
proves that $[\widetilde{H}_2,L_2^2] =0$. Suppose that
$[\widetilde{H}_2,L_2]=a_m(\Lambda) \pa_2^m +
O(\pa_2^{m-1})$ with $a_m\neq 0$. The coefficients in front of $\pa_2^{2m+1}$
in $[\widetilde{H}_2,L_2^2] =  [\widetilde{H}_2,L_2] L_2 + L_2
[\widetilde{H}_2,L_2] $ is $2 a_m(\Lambda)$ and it must be 0 --
contradiction.

b) It is enough to check that $(L_2^2)^\# = \pa_2\cdot (L_2^2)\cdot
\pa_2^{-1}$. Indeed, if this was known then $X:=-\pa_2^{-1}\cdot L_2^\#
\cdot \pa_2$ is a pseudo-differential operator of the form
$\pa_2+O(\pa_2^{-1})$ solving the equation $\pi_2(\L) =
\tfrac{X^2}{2}$. However, we know that such a pseudo-differential
operator is unique and by definition it is $L_2$, i.e., $-\pa_2^{-1}\cdot L_2^\#
\cdot \pa_2=L_2$, which is the first formula that we have to
prove. For the second formula $(L_2^{2l+1})_{2,[0]}=0$ we argue as follows. We have 
$(L_2^{2l+1})_{2,[0]}=\operatorname{Res}_{\pa_2} \left(
  L_2^{2l+1}\pa_2^{-1}\right)$. Let us recall the following formula $\left(\operatorname{Res}_{\pa_2}(P)\right)^\# =
-\operatorname{Res}_{\pa_2}(P^\#)$ for every pseudodifferential operator $P\in
\R(\!(\pa_2^{-1})\!)$. In the case $P=L_2^{2l+1}\pa_2^{-1}$, we have
$P^\#=-P$, so $\operatorname{Res}_{\pa_2} (P)=0$. 

Let us prove that $(L_2^2)^\# = \pa_2\cdot (L_2^2)\cdot
\pa_2^{-1}$. By definition $L_2^2=2\pi_2(\L)$. Let us write the
operator $\L= A+\tfrac{1}{2}\pa_2^2+\tfrac{1}{2}\pa_3^2$. Recalling
Proposition \ref{prop:proj_pi}, c) we have
\ben
\pi_2(A) = 2\Big(A_{1,>0}\, \iota_{\Lambda^{-1}}\left(
H_2^{-1}  (\Lambda^{-1}+1)^{-1}
\right) \Big)_{1,[0]} \pa_2 - 
2\Big(A_{1,<0}\, \iota_{\Lambda}\left(
H_2^{-1}  (\Lambda^{-1}+1)^{-1}
\right) \Big)_{1,[0]} \pa_2 + A|_{\Lambda=-1},
\een
where for $A=\sum_i \alpha_i\Lambda^i$ we set 
$A|_{\Lambda=-1}:= \sum_i (-1)^i\alpha_i$. It is straightforward to
check that  $A|_{\Lambda=-1}=c_3$. On the other hand, by definition
$\tfrac{1}{n-2} L_1^{n-2} =\pi_+(\L)$ and the $\pi_+$ projection of
$\pa_2^2$ and $\pa_3^2$ involves only non-positive powers of
$\Lambda$, while $\pi_+(A)=A$. Therefore, $A_{1,>0} = \tfrac{1}{n-2}
(L_1^{n-2})_{1,>0}$. Furthermore,  
\ben
\iota_{\Lambda^{-1}}\left(
H_2^{-1}  (\Lambda^{-1}+1)^{-1}\right) = 
\iota_{\Lambda^{-1}}\left(
(\pa_2-Q_2)^{-1} (\Lambda-\Lambda^{-1})^{-1}
\right) 
\een
involves only negative powers of $\Lambda$. Therefore, replacing
$A_{1,>0}$ by $\tfrac{1}{n-2}
L_1^{n-2}$ does not change the contribution to the projection
$\pi_2(\L)$. Similarly, replacing $A_{1,<0}$ by $\tfrac{1}{n-2}
(L_1^-)^{n-2}$ does not change the contribution to the projection. We
got the following formula for $\pi_2(A)$
\ben
2\Big(
\tfrac{1}{n-2} L_1^{n-2}\iota_{\Lambda^{-1}}\left(
(\pa_2-Q_2)^{-1} (\Lambda-\Lambda^{-1})^{-1}
\right) \Big)_{1,[0]} \pa_2 - 
2\Big( 
\tfrac{1}{n-2} (L_1^-)^{n-2}\iota_{\Lambda}\left(
(\pa_2-Q_2)^{-1} (\Lambda-\Lambda^{-1})^{-1}
\right) \Big)_{1,[0]} \pa_2 + c_3.
\een
Using $L_1^\# = (\Lambda-\Lambda^{-1}) L_1^-
(\Lambda-\Lambda^{-1})^{-1}$, $Q_2^\#=-(\Lambda-\Lambda^{-1}) Q_2
(\Lambda-\Lambda^{-1})^{-1}$, and that $L_1$ (resp. $L_1^-$) commutes
with $\pa_2-\iota_{\Lambda^{-1}}(Q_2)$
(resp. $\pa_2-\iota_{\Lambda}(Q_2)$) we get that 
\ben
(\pi_2(A))^\# = \pa_2 \cdot \pi_2(A)\cdot \pa_2^{-1} -\pa_2(c_3)\cdot
\pa_2^{-1}. 
\een
A direct computation (or using Proposition \ref{prop:proj_pi}, b))
yields 
\ben
\tfrac{1}{2}\, \Big(\pi_2(\pa_3^2)\Big)^\# = \tfrac{1}{2} \, \pa_2 \cdot
\pi_2(\pa_3^2)\cdot \pa_2^{-1} + \pa_3(q_1)\cdot \pa_2^{-1}.
\een
Therefore, we get the following relation 
\ben
(\pi_2(\L))^\#-\pa_2\cdot \pi_2(\L)) \cdot \pa_2^{-1} =
(\pa_3(q_1)-\pa_2(c_3)) \cdot \pa_2^{-1}.
\een
We claim that the RHS is 0. This follows from part a). Namely,
the operators $\tfrac{1}{2} L_3^2=\pi_3(\L) = \tfrac{\pa_3^2}{2} + c_3
+ O(\pa_3^{-1})$ and $\pa_3^{-1} H_1=\pa_2 + \pa_3^{-1} q_1$ commute,
i.e., 
\ben
\Big[ \tfrac{\pa_3^2}{2} + c_3+ O(\pa_3^{-1}), 
\pa_2 + \pa_3^{-1} q_1
\Big]=0.
\een
Comapring the coefficients in front of $\pa_3^0$ we get $\pa_2(c_3) =
\pa_3(q_1)$. 

c) We will prove a stronger statement. Namely, we claim that for each
pair $(a,b)$ there are
$\alpha_{a,b}$ and $\beta_{a,b}\in \R$ such that
$H_b B_{a,2l+1}-\alpha_{a,b}\Lambda-\beta_{a,b}\in \E H$.
Suppose that $a=2$, the case   $a=3$ is similar. Let us prove
the statement for $b=2$. Recalling Proposition \ref{prop:proj_pi} we get 
\ben
\pi_+(H_2 B_{2,2l+1}) = \tfrac{1}{2}\iota_{\Lambda^{-1}} \Big(
H_2 (L_2^{2l+1})_{2,\geq 0} \widetilde{H}_2^{-1} \Big) (\Lambda-\Lambda^{-1}).
\een
Using $(H_2 L_2^{2l+1})_{2,\geq 0} = H_2 (L_2^{2l+1})_{2,\geq 0} +
(\Lambda-1)\operatorname{res}_{\pa_2}(L_2^{2l+1})$ and the fact that
$\widetilde{H}_2^{-1}=(\Lambda-1)^{-1} + O(\pa_2^{-1})$ has only
non-positive powers of $\pa_2$ (see formula \eqref{h2-expansion}) we get
\ben
\pi_+(H_2 B_{2,2l+1}) = \tfrac{1}{2}\iota_{\Lambda^{-1}} \Big(
H_2 L_2^{2l+1}) \widetilde{H}_2^{-1} \Big) (\Lambda-\Lambda^{-1}) -
\tfrac{1}{2} (\Lambda-1) \resi_{\pa_2}(L_2^{2l+1}) (\Lambda^{-1}+1).
\een
According to part a) $L_2$ and $\widetilde{H}_2$ commute, $H_2
\widetilde{H}_2^{-1} = \pa_2+q_2$, and by part b) $(L_2^{2l+1})_{2,[0]}=0$. Therefore
\ben
\pi_+(H_2 B_{2,2l+1}) = \tfrac{r}{2}(\Lambda-\Lambda^{-1})
-(\Lambda-1)\tfrac{r}{2}(\Lambda^{-1}+1) =\tfrac{1}{2}(r-r[1]) (\Lambda+1),
\een
where $r=\operatorname{Res}_{\pa_2} (L_2^{2l+1})$. From here, recalling
Proposition \ref{prop:ker_pi_a}, we get
$H_2 B_{2,2l+1} - \tfrac{1}{2}(r-r[1])(\Lambda+1) \in \E_{(+)}H\cap \E
= \E H.$
In particular, we get the following formulas
\ben
\alpha_{2,2} = \beta_{2,2} = \tfrac{1}{2}(r-r[1]).
\een
Let us move to the case $b=1$. Recalling Proposition
\ref{prop:proj_pi} we get
\ben
\pi_3(H_1 B_{2,2l+1}) = \Big( H_1 B_{2,2l+1}
\widetilde{H}_1^{-1}\Big)_{2,[0]} \, \pa_3.
\een
On the other hand, using that $(H_1 L_2^{2l+1})_{2,\geq 0} = H_1
(L_2^{2l+1})_{2,\geq 0}+ \pa_3\cdot
\operatorname{res}_{\pa_2}(L_2^{2l+1})$ and that
$\widetilde{H}_1^{-1}=\pa_3^{-1} + O(\pa_2^{-1})$ has only
non-positive powers of $\pa_2$ we get
\ben
\pi_3(H_1 B_{2,2l+1}) = \Big( H_1 L_2^{2l+1} \widetilde{H}_1^{-1}
\Big)_{2,[0]} \pa_3 -\pa_3 \cdot r,
\een
where $r=\operatorname{Res}_{\pa_2} (L_2^{2l+1})$. According to part a), the
operators $L_2$ and  $\widetilde{H}_1$ commute. Since $H_1
\widetilde{H}_1^{-1} = \pa_2$ we get
\ben
\pi_3(H_1 B_{2,2l+1}) = r\cdot \pa_3 -\pa_3\cdot r = -\pa_3(r).
\een
Recalling Proposition \ref{prop:ker_pi_a} we get
$H_1 B_{2,2l+1} +\pa_3(r) \in \E_{(3)}H\cap \E = \E H.$ Let us point
out that
\ben
\alpha_{2,1} =0,\quad \beta_{2,1} = -\pa_3(r). 
\een
Finally, for the case $b=3$, let us recall the identity
\ben
(\Lambda-1) H_1 = q_2 H_3 +\pa_3 \cdot H_2.
\een
Multiplying both sides by $B_{2,2l+1}$, taking the projection $\pi_+$,
and using that
\ben
\pi_+(H_1 B_{2,2l+1}) =-\pa_3(r)
\pi_+(H_2 B_{2,2l+1}) = \tfrac{1}{2}(r-r[1]) (\Lambda+1),\quad
\een
we get the following identity
\ben
-(\Lambda-1) \pa_3(r) = q_2 \pi_+(H_3 B_{2,2l+1}) +
\tfrac{1}{2}\pa_3(r-r[1]) (\Lambda+1) + \tfrac{1}{2} (r-r[1]) q_2 (\Lambda-1),
\een
which yields
\ben
q_2 \pi_+(H_3 B_{2,2l+1}) = -\tfrac{1}{2}\Big(
\pa_3(r+r[1]) +q_2(r-r[1])\Big)(\Lambda -1).
\een
Comparing the coefficients in front of $\Lambda$ we get that
$\pi_+(H_3 B_{2,2l+1}) =\alpha_{2,3} \Lambda + \beta_{2,3}$ for some
$\alpha_{2,3},\beta_{2,3}\in \R$ satisfying
\ben
q_2\beta_{2,3}=-q_2 \alpha_{2,3}  =  \tfrac{1}{2}\Big(
\pa_3(r+r[1]) +q_2(r-r[1])\Big).
\een
It remains only to use that $H_3 B_{2,2l+1}-\alpha_{2,3}(\Lambda-1)\in
\E_{(+)}H\cap \E=\E H.$

d) Let us give the argument for $a=2$. The case $a=3$ is
similar. Since $B^\pm_{2,2l+1}=\pi_\pm((L_2^{2l+1})_{2,\geq 0})$ and
$H_i\mathcal{L}\in \E H$, we have
$\pi_\pm[B^\pm_{2,2l+1},\L]=\pi_\pm[(L_2^{2l+1})_{2,\geq 0},\L]$. It is
sufficient to prove that the $\pi_+$-projection of $
[(L_2^{2l+1})_{2,\geq 0},\L]$ is a difference operator, i.e., it
belongs to the ring $\R[\Lambda,\Lambda^{-1}]$. Indeed, if we knew
this fact, then 
\ben
[(L_2^{2l+1})_{2,\geq 0},\L]-\pi_+[(L_2^{2l+1})_{2,\geq 0}),\L]\in
\E\cap \E_{(+)} H = \E H,
\een
where we used Proposition \ref{prop:ker_pi_a}.
The above statement implies that $\pi_+[(L_2^{2l+1})_{2,\geq
  0}),\L]=\pi_-[(L_2^{2l+1})_{2,\geq 0}),\L]$, which would complete
the proof. 

Note that 
\ben
[(L_2^{2l+1})_{2,\geq 0},\L] = 
[L_2^{2l+1},\L]_{2,\geq 0} + \pa_2(\resi_{\pa_2} (L_2^{2l+1})).
\een
Therefore, it is enough to prove that $\pi_+[L_2^{2l+1},\L]_{2,\geq
  0}$ is a difference operator. Let us decompose the operator 
\beq\label{pi_2-decomp}
\L=\tfrac{1}{2}L_2^{2} + A \widetilde{H}_2 + \pa_3 \cdot
\widetilde{H}_1- \tfrac{1}{2}\widetilde{H}_1^2,
\eeq
where $\widetilde{H}_2 =(\pa_2+q_2)^{-1} H_2$ and
$\widetilde{H}_1=\pa_2^{-1} H_1$. The last two terms in
\eqref{pi_2-decomp} coincide with
$\tfrac{1}{2}(\pa_3^2-\pi_2(\pa_3^2))$. The terms of $\L$ that involve
only $\Lambda$ decompose by the following recursive rules:
\ben
\Lambda^i= (\pa_2-q_2[i-1])^{-1} \Lambda^{i-1} H_2 + 
(\pa_2-q_2[i-1])^{-1} (\pa_2+q_2[i-1]) \Lambda^{i-1}
\een 
and
\ben
\Lambda^{-i}=- (\pa_2+q_2[-i])^{-1} \Lambda^{-i} H_2 + 
(\pa_2+q_2[-i])^{-1} (\pa_2-q_2[-i]) \Lambda^{-i+1},
\een
where $i>0$ is an integer. Therefore, the coefficient $A\in
\R[\Lambda,\Lambda^{-1}](\!(\pa_2^{-1})\!)$. According to part a) the
operators $L_2, \widetilde{H_2}$, and $\widetilde{H}_1$ pairwise
commute. Therefore, 
\ben
[L_2^{2l+1},\L]_{2,\geq 0} = 
\left([L_2^{2l+1}, A] \widetilde{H}_2\right)_{2,\geq 0} - 
\left(\pa_3(L_2^{2l+1}\pa_2^{-1}) H_1\right)_{2,\geq 0}.
\een
Note that
\ben
\left(\pa_3(L_2^{2l+1}\pa_2^{-1}) H_1\right)_{2,\geq 0} -
\left(\pa_3(L_2^{2l+1}\pa_2^{-1})\right)_{2,\geq 0} H_1 =
\operatorname{Res}_{\pa_2}(\pa_3(L_2^{2l+1}\pa_2^{-1}) )\cdot \pa_3=0,
\een
where for the last equality we used that $(L_2^{2l+1})_{2,[0]}=0$ (see
part b) ). We get
\ben
\pi_+\left( [L_2^{2l+1},\L]_{2,\geq 0}\right)  =
\left([L_2^{2l+1}, A] \widetilde{H}_2\right)_{2,[0]} +
\pi_+\left( \left([L_2^{2l+1}, A] \widetilde{H}_2\right)_{2,> 0} \right).
\een
The first term on the RHS is a difference operator, while for the
second one we recall Proposition \ref{prop:proj_pi} part a) and we get
\ben
\pi_+\left( \left([L_2^{2l+1}, A] \widetilde{H}_2\right)_{2,> 0} \right) =
\tfrac{1}{2}\left(
  [L_2^{2l+1}, A] -([L_2^{2l+1}, A] \widetilde{H}_2 )_{2,[0]} \widetilde{H}^{-1}_2
  \right)_{2,[0]} (\Lambda-\Lambda^{-1}).
\een
This is clearly a difference operator.
\qed

We claim that the operators $B^\pm:=B^\pm_{a,2l+1}$ satisfy the conditions
(i)--(iii) of Lemma \ref{le:der_R}. Conditions (ii) and (iii) follow
from respectively part d) and c) of Lemma \ref{le:pa_2-identities}. We
need only to verify that condition (i) is satisfied, that is, 
\beq\label{B_2odd-invol}
(B_{a,2l+1}^+)^\# = -(\Lambda-\Lambda^{-1})\,
B_{a,2l+1}^-\,
\iota_{\Lambda}(\Lambda-\Lambda^{-1})^{-1}.
\eeq
Again let us give the argument only for $a=2$. The case $a=3$ is
similar. 
Recall that the residue of a pseudo-differential operator $P=\sum_{j}
p_j \pa_2^{j}$ is defined by $\resi_{\pa_2} (P):= p_{-1}$. We have
\ben
B_{2,2l+1}^{\pm} = \tfrac{1}{2} \iota_{\Lambda^{\mp 1}}
\operatorname{Res}_{\pa_2} \, \Big(
L_2^{2l+1} \widetilde{H}_2^{-1} \pa_2^{-1} \Big)
(\Lambda-\Lambda^{-1}). 
\een
Using the following simple formula
\beq\label{res-invol}
\Big( \operatorname{Res}_{\pa_2} \, P(\Lambda,\pa_2)\Big)^\# =
-\operatorname{Res}_{\pa_2} \, P^\#(\Lambda,\pa_2)
\eeq
we get
\ben
(B_{2,2l+1}^+)^\# =\tfrac{1}{2} (\Lambda-\Lambda^{-1})
\operatorname{Res}_{\pa_2} \Big(
(-\pa_2)^{-1} (\widetilde{H}_2^\#)^{-1} (L_2^\#)^{2l+1}
\Big).
\een
Recalling Lemma \ref{le:pa_2-identities}, a) and b) and
\ben
\widetilde{H}_2^\# = -\pa_2 \cdot
\widetilde{H}_2\cdot (\widetilde{H}_2+1)^{-1} \cdot \pa_2^{-1}
\een
we get
\ben
(B_{2,2l+1}^+)^\# =-\tfrac{1}{2} (\Lambda-\Lambda^{-1})
\operatorname{Res}_{\pa_2} \Big(
L_2^{2l+1} \widetilde{H}_2^{-1} \pa_2^{-1} + L_2^{2l+1} \pa_2^{-1}\Big).
\een
In order to prove \eqref{B_2odd-invol} we need only to recall that
according to Lemma \ref{le:pa_2-identities}, b) 
$\operatorname{Res}_{\pa_2} \Big( L_2^{2l+1} \pa_2^{-1}\Big) = 0$.

\section{Existence of the extended flows}\label{existofextflows}

\subsection{Dressing operators}
Let us recall the notation from Section \ref{sec:dress_op}. 
\begin{proposition}\label{commuteonR1}
  a)
The derivations $\pa_2$ and $\pa_3$ can be extended uniquely to
derivations of $\R_1$ such that $\epsilon\pa_x$, $\pa_2$, and $\pa_3$ pairwise
commute and $\pa_a(S_1) = Q_a S_1$, $a=2,3$.

b) The coefficients of the operator series
$\ell_1:=\epsilon\pa_x(S_1)\cdot S_1^{-1}$ belong to $\R$. 
\end{proposition}
\proof
Let us extend the derivations $\pa_a$ $(a=2,3)$ to 
derivations of $\R_1$ via the given formula $\pa_a(S_1) = Q_a \cdot
S_1$. We need only to prove that $\epsilon\pa_x$, $\pa_2$, and $\pa_3$
pairwise  commute as derivations of $\R_1$. The commutativity of
$\pa_2$ and $\pa_3$ is a direct consequence of the 0-curvature
condition.  Let us prove that $\pa_2$ and $\epsilon\pa_x$ commute. The
argument for $\pa_3$ and $\epsilon\pa_x$ is similar.  
We have to prove that $\epsilon\pa_x \pa_2(S_1) = \pa_2 \epsilon \pa_x (S_1)$.

The first step in our argument is to prove that the derivation $\pa_2$
commutes with the translation operator $\Lambda^{n-2}$.  
Using that
\ben
\frac{1}{n-2} L_1^{n-2} = \pi_+(\L) = 
\L-\frac{1}{2} \Big( 
(\pa_2+Q_2)(\pa_2-Q_2) + (\pa_3+Q_3)(\pa_3-Q_3)\Big) 
\een
and 
\ben
\pa_2(\L) = [Q_2,\L-\tfrac{1}{2}(\pa_2^2+\pa_3^2)] -\sum_{a=2,3}
\Big(\frac{1}{2} \pa_a^2(Q_2) + \pa_a(Q_2) Q_a\Big)
\een
we get $\pa_2 (L_1^{n-2}) = [Q_2,L_1^{n-2}]$.  Let us substitute $L_1^{n-2} =
\sum_{i=0}^\infty b_{n-2,i} \Lambda^{n-2-i}$ and
$S_1=\sum_{j=0}^\infty \psi_{1,j}\Lambda^{-j}$ in $L_1^{n-2} S_1 = S_1
\Lambda^{n-2}$ and differentiate the resulting identity with respect
to $\pa_2$. We get 
\ben
\pa_2(L_1^{n-2}) S_1 +
\sum_{k=0}^\infty \sum_{i=0}^k b_{m,i}
\pa_2(\psi_{1,k-i}[n-2-i])\Lambda^{n-2-k} = \pa_2(S_1) \Lambda^{n-2}. 
\een
Since $\pa_2(L_1^{n-2})=[Q_2,L_1^{n-2}]$, $\pa_2(S_1)=Q_2 S_1$, and
$L_1^{n-2}S_1=S_1 \Lambda^{n-2}$, the above identity yields
\ben
\sum_{i=0}^k b_{m,i}
\pa_2(\psi_{1,k-i}[n-2-i])=
\sum_{i=0}^k b_{m,i}
\pa_2(\psi_{1,k-i}) [n-2-i].
\een
A simple induction on $k$ yields $\pa_2\, \Lambda^{n-2}(\psi_{1,k}) =
\Lambda^{n-2} \, \pa_2(\psi_{1,k})$ for all $k\geq 0$.

Now we are in position to prove that
$(\epsilon\pa_x\circ \pa_2)(\psi_{1,i}) =
(\pa_2\circ\epsilon\pa_x)(\psi_{1,i})$ for all $i\geq 0$. For $i=0$ the equality
is straightforward to check. In particular, it is sufficient to prove
that $(\epsilon\pa_x\circ \pa_2)(\psi_{1,i}/\psi_{1,0}) =
(\pa_2\circ\epsilon\pa_x)(\psi_{1,i}/\psi_{1,0})$. We argue by
induction on $i$. Formula \eqref{x-der-S1} can be written
as $(1-\Lambda)(\psi_{1,i}/\psi_{1,0}) = f_k$, where $f_k\in \R_1$
is a formal power series in $\epsilon$ whose coeffcients are
differential polynomials in $\phi, \psi_{1,1},\dots,\psi_{1,i-1}$ with
coefficients in $\R$. Using the iductive assumption we get
$\pa_2(\epsilon\pa_x)^l(f_k)=(\epsilon\pa_x)^l\pa_2(f_k)$ for all
$l\geq 0$. Therefore
\ben
\pa_2\, (1-\Lambda^{n-2})\Big(\frac{\psi_{1,i}}{\psi_{1,0}} \Big) =
\pa_2 \, \frac{1-\Lambda^{n-2}}{1-\Lambda } (f_k) =
\frac{1-\Lambda^{n-2}}{1-\Lambda }\, \pa_2  (f_k).
\een
Let us take the equality between the first and the last term in the
above formula.  After exchanging the order of $\pa_2$ and
$1-\Lambda^{n-2}$ on the LHS and applying to both sides the
operator $\tfrac{\epsilon\pa_x}{1-\Lambda^{n-2}}$ we get
\ben
\epsilon\pa_x \,\pa_2 \Big(\frac{\psi_{1,i}}{\psi_{1,0}} \Big) =
\frac{ \epsilon\pa_x }{1-\Lambda }\, \pa_2  (f_k) = 
\pa_2\, \frac{ \epsilon\pa_x }{1-\Lambda }(f_k)=
\pa_2\, \epsilon\pa_x \Big(\frac{\psi_{1,i}}{\psi_{1,0}} \Big),
\een 
where for the second equality we used the inductive assumption and for
the thrid equality we recall formula \eqref{x-der-S1}.

b)
The idea of the proof is the same as the proof of Theorem 2.1 in
\cite{CDZ}. Let us first prove the following general statement about
pseudo-difference operators: let 
\ben
\widetilde{\mathcal{R}}=\CC[\pa_x^m\widetilde{\psi}_j\,|\, m\geq 0,\ j\geq
1][\![\epsilon]\!]
\een
be the ring of formal power series in $\epsilon$ whose coefficients are
differential polynomials on the infinite number of formal variables
$\widetilde{\psi}_1, \widetilde{\psi}_2,\dots$. The translation
operator $\Lambda=e^{\epsilon \pa_x}$ 
acts naturally on $\widetilde{\R}$, so we can define the ring of
pseudo-difference operators $\widetilde{\R}(\!(\Lambda^{-1})\!)$.  Let
us define the following pseudo-difference operators:
\ben
\widetilde{S}:=1+\sum_{j\geq 1} \widetilde{\psi}_j \Lambda^{-j},\quad
\widetilde{L} := \widetilde{S}  \Lambda \widetilde{S} ^{-1},\quad
\widetilde{\ell}:=\epsilon \pa_x(\widetilde{S} ) \widetilde{S} ^{-1}.
\een
We claim that the coefficients of $\widetilde{\ell}$ are differential
polynomials in the coefficients of $\widetilde{L}$, that is, the
coefficients of $\widetilde{\ell}$ belong to the differential subring
of $\widetilde{\R}$ generated by the coefficients of $\widetilde{L}$.

Let us prove the above claim. The idea is to compare the coefficient in front of $\Lambda^0$
in the following equations: 
\ben
\epsilon\pa_x (\widetilde{L}^m)=[\widetilde{\ell},
\widetilde{L}^m],\quad m\geq 1.
\een 
Let us write $\widetilde{\ell}=\sum_{i=1}^\infty a_i \Lambda^{-i}$
and $\widetilde{L}^m = \Lambda^m + \sum_{j=1}^\infty
u_{m,j}\Lambda^{m-j}$, then we get 
\beq\label{a_m-recursion}
\epsilon \pa_x(u_{m,m}) = 
(1-e^{m\epsilon\pa_x})(a_m) + 
\sum_{i=1}^{m-1} (1-e^{i\epsilon \pa_x}) (a_i u_{m,m-i} [-i]).
\eeq
In order to complete the proof it is sufficient to prove that
\beq\label{a_m}
a_m= \frac{\epsilon \pa_x}{1-e^{m\epsilon \pa_x}} (u_{m,m}) -
\sum_{i=1}^{m-1}
\frac{1-e^{i\epsilon \pa_x}}{1-e^{m\epsilon \pa_x}}(a_i u_{m,m-i} [-i]).
\eeq
The recursion \eqref{a_m-recursion} implies that the derivatives of
the LHS and the RHS of \eqref{a_m} with respect to $x$ are the same. We have to prove
only that the integration constant is 0. Let us turn
$\widetilde{\mathcal{R}}$ into a graded ring by assigning degree $i$
to $\pa_x^m\widetilde{\psi}_i$ for $i\geq 1$, $m\geq 0$.  We leave it as an
exercise to check that the ring of constants of $\widetilde{R}$, that
is, the set of elements $A\in \widetilde{R}$, such that, $\pa_x(A)=0$,
is $\CC[\![\epsilon]\!]$. In other words, the ring of constants
coincides with the subring of homogeneous elements of degree 0. 
On the other hand,  note that the coefficients $a_i$ and $u_{m,i}$ are
homogeneous of degree $i$. The difference of the LHS and the RHS in
\eqref{a_m} is homogeneous of
degree $m\geq 1$ and its derivative is $0$, so it must be identically
0.

Now the proof of part b) can be completed as follows: 
Put $\widetilde{S}=\psi_{1,0}^{-1} S_1$. Note that the coefficients of the
pseudo-difference operator $\widetilde{L}:= \widetilde{S} \Lambda
\widetilde{S}^{-1} = \psi_{1,0}^{-1} L_1 \psi_{1,0} $ belong to
$\R$, because $e^{\phi[i]-\phi}\in \R$ for all $i\in \ZZ$.  Recalling the
claim from above, we get that the coefficients of 
$\widetilde{\ell}=\epsilon\pa_x(\widetilde{S})\,
\widetilde{S}^{-1}$ also belong to $\R$. Finally, the coefficients of
the operator
\ben
\ell_1 =\epsilon\pa_x(S_1)\cdot S_1^{-1} =
\frac{\epsilon\pa_x(\psi_{1,0})}{\psi_{1,0}} +
\psi_{1,0} \cdot \widetilde{\ell} \cdot \psi_{1,0}^{-1} =
\frac{(n-2)\epsilon\pa_x}{1-\Lambda^{n-2}}(\alpha) +
\psi_{1,0} \cdot \widetilde{\ell} \cdot \psi_{1,0}^{-1} 
\een
belong to $\R$, which is what we have to prove.
\qed

\begin{proposition}\label{commuteonR23}
  Suppose that $a=2,b=3$ or $a=3,b=2$.

  a) The translation operator $\Lambda=e^{\epsilon\pa_x}$ and $\pa_b$
  can be extended uniquely to 
respectively an automorphism and a derivation of the ring $\R_a$ in
such a way that $\Lambda$, $\pa_2$, and $\pa_3$ pairwise commute, and
\ben
T_a[1] = (\pa_a-q_a)^{-1} (\pa_a+q_a) T_a,\quad
\pa_b(T_a)=-\pa_a^{-1}\cdot q_1 \cdot T_a.
\een

b) The derivation $\epsilon \pa_x$ can be extended uniquely to a
derivation on $\R_a$, such that $\epsilon \pa_x, \pa_2$, and $\pa_3$
pairwise commute and $\Lambda=e^{\epsilon\pa_x}$. The coefficients of
the operator $\ell_a:=\epsilon\pa_x(S_a)\cdot S_a^{-1}$ belong to
$\R$. 
\end{proposition}
\proof
a) 
Let us give the argument for the case $a=2$ and $b=3$. The other case,
$a=3$ and $b=2$ is similar. By definition,
\ben
\pa_2(T_2) = (\pa_2-L_2)\cdot T_2. 
\een
Let us extend the translation operator $\Lambda$ and the derivation
$\pa_3$ to $\R_2$ by the given formulas, that is,
\begin{align}\nonumber
  T_2[1] & =   (\pa_2-q_2)^{-1} (\pa_2+q_2) T_2,\\
  \nonumber
\pa_3(T_2) & =  -\pa_2^{-1}\cdot q_1 \cdot T_2.
\end{align}
The commutativity of the extended translation operator and the
extended derivation $\pa_3$ is equivalent to
\ben
\pa_3 \left(  (\pa_2-q_2)^{-1} (\pa_2+q_2) T_2\right) =
- \pa_2^{-1}\cdot q_1[1] \cdot T_2[1]
\een
Substituting the above formulas, after a short computation, we get
that the above identity is equivalent to the fact that
$\widetilde{H}_1:= \pa_3 + \pa_2^{-1} q_1$ and $\widetilde{H}_2:=
(\pa_2+q_2)^{-1} (\pa_2-q_2) \Lambda$ commute. This however, follows
from part a) of Lemma \ref{le:pa_2-identities}. Similarly, the
commutativity of $\Lambda$ and $\pa_2$ is equivalent to
\ben
\pa_2 \left(  (\pa_2-q_2)^{-1} (\pa_2+q_2) T_2\right) =
  (\pa_2-L_2[1])\cdot T_2[1].  
\een
The above identity is equivalent to the fact that $L_2$ and
$\widetilde{H}_2$ commute, which again follows from part a) of Lemma
\ref{le:pa_2-identities}. Finally, the commutativity of $\pa_2$ and
$\pa_3$ is equivalent to the commutativity of $L_2$ and $\widetilde{H}_1$.

b) Let us define pseudo-differential operators 
$\tilde{\ell}_2^{(m)}\in \R[\![\pa_2^{-1}]\!]\pa_2^{-1}$ ($m\geq 1$)  
by the following two conditions:
\begin{enumerate}
	\item[(i)] $\tilde{\ell}_2^{(m+1)}=\tilde{\ell}_2\cdot \tilde{\ell}_2^{(m)}+\epsilon\pa_x\tilde{\ell}_2^{(m)}$ with $\tilde{\ell}_2^{(0)}=1$ and $\tilde{\ell}_2^{(1)}=\tilde{\ell}_2$,
	\item[(ii)] $\sum_{m=0}^{+\infty}\tilde{\ell}_2^{(m)}/m!=(\pa_2+q_2)^{-1}(\pa_2-q_2)$.
        \end{enumerate}
Let us check that conditions (i) and (ii) uniquely determine the
sequence $\widetilde{\ell}_2^{(m)}$. Condition (i) defines $
\widetilde{\ell}_2^{(m)}$ in terms of $\widetilde{\ell}_2$, so we need only to
check that $\widetilde{\ell}_2=    \sum_{i\geq 1} a_i \pa_2^{-i}$ is
uniquely determined by condition (ii). Note that
$
\widetilde{\ell}_2^{(m)} = (\epsilon\pa_x)^{m-1} (\ell_2) + \cdots,
$
where the dots stand for a differential polynomial in
$\widetilde{\ell}_2$ that involves at least quadratic terms in
$\widetilde{\ell}_2$. In particular, the ceofficient in front of
$\pa_2^{-i}$ in $\sum_{m\geq 1} \widetilde{\ell}_2^{(m)}/m!$  has the form
\ben
\frac{e^{\epsilon \pa_x}-1}{\epsilon \pa_x}(a_i) +
\text{terms involving $a_1,\dots,a_{i-1}$} 
\een
Since the operator $ \tfrac{e^{\epsilon \pa_x}-1}{\epsilon \pa_x}$ is
invertible, the identity in condition (ii) is equivalent to a system
of recursions, which uniquely define the coefficients $a_i$ ($i\geq
1$).

Let us extend the derivation $\epsilon\pa_x$ to a derivation on $\R_2$
by the following formula: $\epsilon\pa_x(T_2^{-1})=T_2^{-1}\cdot
\tilde{\ell}_2$. Using the fact that $\pa_2$ commutes with
$\epsilon\pa_x$ on $\mathcal{R}$ and condition (i), we get 
$(\epsilon\pa_x)^m(T_2^{-1})=T_2^{-1}\tilde{\ell}_2^{(m)}$ for all
$m\geq 0$. Recalling condition (ii), we get that the extension of the
translation operator to $\R_2$ coincides with $e^{\epsilon\pa_x}$. It
remains only to prove that $[\pa_2,\epsilon\pa_x]=0$ on $\R_2$.

Put $c_m=T_2\cdot [\pa_2,(\epsilon\pa_x)^m](T_2)$. According to part
a), we have $[\pa_2,\Lambda]=0$ on $\R_2$ $\Rightarrow$
$\sum_{m=1}^{+\infty}c_m/m!=0$. By using $[\pa_2,\epsilon\pa_x]=0$ on
$\mathcal{R}$, we get the following recursion relations for the
commutators $c_m$:
$c_{m+1}=\tilde{\ell}_2c_m+c_1\tilde{\ell}_2^{(m)}+\epsilon\pa_xc_{m}$. Using
this recursion relation, by induction on $m$, we get the following
formula: 
\begin{align}
c_m=\sum_{k=0}^{m-1}\sum_{i=0}^k\binom{k}{i}\tilde{\ell}_2^{(i)}(\epsilon\pa_x)^{k-i}(c_1\cdot\tilde{\ell}_2^{(m-1-k)})=(\epsilon\pa_x)^{m-1}(c_1)+\text{at least quadratic terms}.
\end{align}
Then assume $c_1=\sum_{i=1}^{\infty} c_{1,i}\pa_2^{-i}$. By comparing the
coefficient of $\pa_2^{-i}$ in $\sum_{m=1}^{+\infty}c_m/m!=0$, we get 
\ben
\frac{e^{\epsilon \pa_x}-1}{\epsilon \pa_x}(c_{1,i}) +
\text{terms involving $c_{1,1},\dots,c_{1,i-1}$} =0.
\een
Since the operator $\frac{e^{\epsilon \pa_x}-1}{\epsilon \pa_x}$ is
invertible, by induction on $i$, we get $a_i=0$ for all $i\geq 1$,
that is, $c_1=0$.

Finally, note that $\ell_2=-\widetilde{\ell}_2$, so the coefficients
of $\ell_2$ belong to $\R$ as claimed.
\qed

\begin{corollary}\label{cor:s2-symm}
  The operator
  \ben
  T_a^\# \,\pa_a \, T_a\, \pa_a^{-1}\quad \in \quad \R_a
  [\![\pa_a^{-1}]\!] ,\quad a=2,3,
  \een
  commutes with $\epsilon \pa_x,\pa_2$, and $\pa_3$. 
\end{corollary}
\proof
The commutativity $[P,\epsilon\pa_x]=0$ is equivalent to
$P[1]=P$. Indeed, let us assume that $P[1]=P$ (the other implication
is obvious). Expanding $P=\sum_{k=0}^\infty P_k \epsilon^k$ and
comparing the coefficients in front of $\epsilon^1$ in $P[1]=P$ we get
$\pa_x(P_0)=0$ and hence $P_0[1]=P_0$. Clearly, we can continue
inductively and prove that $\pa_x(P_k)=0$ and $P_k[1]=P_k$ for all
$k$. 

Let us check that
\beq\label{transl-sym}
T_a^\#[1] \,\pa_a \, T_a[1]\, \pa_a^{-1}=T_a^\# \,\pa_a \, T_a\, \pa_a^{-1}.
\eeq
Using that
\ben
T_a[1] = (\pa_a-q_a)^{-1} (\pa_a+q_a) T_a
\een
and
\ben
T_a^\# [1]=\Lambda T_a^\# \Lambda^{-1} = (T_a[1])^\# = T_a^\#
(\pa_a-q_a) (\pa_a+q_a) ^{-1} ,
\een
we get that the LHS of \eqref{transl-sym} is equal to
\ben
T_a^\#
(\pa_a-q_a) (\pa_a+q_a) ^{-1} \pa_a (\pa_a-q_a)^{-1} (\pa_a+q_a) T_a
\pa_a^{-1}.
\een
This expression coincides with the RHS of \eqref{transl-sym}, because
\ben
\pa_a (\pa_a-q_a)^{-1} (\pa_a+q_a)  = (1-q_a \pa_a^{-1})^{-1} (1+q_a
\pa_a^{-1}) \pa_a
\een
and $(\pa_a-q_a) (\pa_a+q_a) ^{-1} = (1-q_a
\pa_a^{-1})(1+q_a\pa_a^{-1})^{-1}$.

The commutativity with $\pa_b$ is proved in a similar way,
using the explicit formula $\pa_b(T_a)=-\pa_a^{-1} q_1 T_a$. Finally,
to prove the commutativity with $\pa_a$, let us conjugate the identity
$L_a T_a=T_a\pa_a$, we get $T_a^\# L_a^\# = -\pa_a T_a^\#$. On the
other hand, according to part b) of Lemma \ref{le:pa_2-identities}, we
have $L_a^\# = -\pa_a L_a \pa_a^{-1} = -\pa_a T_a \pa_a T_a^{-1}
\pa_a^{-1}$. Therefore,
\ben
T_a^\# \pa_a T_a \pa_a T_a^{-1} \pa_a^{-1} = \pa_a T_a^\#
\quad\Rightarrow\quad
[T_a^\# \pa_a T_a \pa_a^{-1} , \pa_a ]=0.\qed
\een

Using Corollary \ref{cor:s2-symm}, we can give a proof of Proposition
\ref{prop:S_a} as follows:
Put $A:=T_a^\# \pa_a T_a \pa_a^{-1}\in 1+ \R_a[\![\pa_a^{-1}]\!]$. Note
that $A^\# = \pa_a^{-1} A \pa_a = A$, where for the second identity we
used that $A$ commutes with $\pa_2$. Therefore, $A=1+
\sum_{i=1}^\infty a_i \pa_2^{-2i}$. Recalling Corollary
\ref{cor:s2-symm} we also have that the derivatives with respect to
$\pa_x$, $\pa_2$ and $\pa_3$ of $a_i$ are 0. There
exists a unique operator $B=1+\sum_{j=1}^\infty b_j \pa_2^{-2j}$ such that
$B^2= A$. Indeed, comparing the coefficients in front of
$\pa_2^{-2i}$ we get a system of recursion relations for $b_i$ of the
form $2b_i + \cdots =a_i$, where the dots stand for terms involving
$b_1,\dots, b_{i-1}$. Clearly, the operator $B$ commutes with the 3
derivations and $B^\#=B$. Therefore the operator $S_a:= T_a B^{-1}$
has the required properties. \qed

\subsection{Extension of the coefficient ring $\R$}
It will be convenient to extended the coefficient ring $\R$ to
$\hat{\R}:=\R[\pa_x]$, that is, let us define
\ben
\hat{\E} = \hat{R}[\Lambda^{\pm 1},\pa_2,\pa_3],
\een
\ben
\hat{\mathcal{E}}_{(\pm)}=\hat{\R}[\pa_2,\pa_3]((\Lambda^{\mp 1})),\quad
\hat{\mathcal{E}}_{(2)}=\hat{\R}[\Lambda,\Lambda^{-1},\pa_3]((\pa_2^{-1})),\quad
\hat{\mathcal{E}}_{(3)}=\hat{\R}[\Lambda,\Lambda^{-1},\pa_2]((\pa_3^{-1}))
\een
and
\ben
\hat{\mathcal{E}}_{(\pm)}^{0}=\hat{\R}((\Lambda^{\mp 1})),\quad
\hat{\mathcal{E}}_{(2)}^{0}=\hat{\R}((\pa_2^{-1})),\quad
\hat{\mathcal{E}}_{(3)}^{0}=\hat{\R}((\pa_3^{-1})).
\een
We still have a direct sum decomposition $\hat{\A}=\hat{\A}^0\oplus \hat{\A} H$ for all
$\A\in \{\mathcal{E}_{(\pm)}, \mathcal{E}_{(2)},\mathcal{E}_{(3)}\}$,
where $\hat{\A}H$ is the left ideal in $\hat{\A}$ generated by
$H_1,H_2,$ and $H_3$. Therefore, we can define the projections
$\pi_\alpha: \hat{\A}_\alpha\to \hat{\A}^0_\alpha.$ Moreover, the
natural generalization of both
Proposition \ref{prop:ker_pi_a} and Proposition \ref{prop:proj_pi}
still hold.   
\subsection{The extended flows}

We have the following symmetry
\beq\label{A_ak-symm}
(A_{a,k})^\# = -\pa_a A_{a,k} \pa_a^{-1},\quad a=2,3.
\eeq
Indeed, since $S_a^\# = \pa_a S_a^{-1}\pa_a^{-1}$ and by definition,
\ben
A_{a,k} = S_a \tfrac{\pa_a^{2k}}{2^k k!} (\epsilon\pa_x) S_a^{-1},
\een
the identity \eqref{A_ak-symm} is straightforward to prove.
Let us define
\ben
B_{0,l}^+= (B_{0,l,1}^++B_{0,l,2}^++B_{0,l,3}^+)(\Lambda-\Lambda^{-1}),
\een
where
\ben
B_{0,l,1}^+:= \sum_{m=0}^\infty  \left(\left(
  \frac{L_1^{(n-2)l}}{(n-2)^l l!} (\log L_1 -h_l)
  \Lambda^{-2m-1} \right)_{1,\geq 0} +\left(
  \Lambda^{2m+1} \Big( \frac{L_1^{(n-2)l}}{(n-2)^l l!} (\log L_1
  -h_l)\Big)^\#
  \right)_{1,<0}\right),
  \een
  \ben
  B_{0,l,2}^+:=\frac{1}{2}\,\,\iota_{\Lambda^{-1}} \left(
    \frac{L_2^{2l}}{2^l l!}\, H_2^{-1} (\pa_2+q_2)\,
    \log \Big( (\pa_2+q_2)^{-1} H_2 +1\Big)
    \right)_{2,[0]},
  \een
  and
  \ben
  B_{0,l,3}^+:=\frac{1}{2}\,\,\iota_{\Lambda^{-1}} \left(
    \frac{L_3^{2l}}{2^l l!}\, H_3^{-1} (\pa_3+q_3)\,
    \log \Big( (\pa_3+q_3)^{-1} H_3 -1\Big)
    \right)_{3,[0]}.
    \een
Note that $B_{0,k,1}=B_{0,k,1}^+(\Lambda-\Lambda^{-1})$. Recalling
parts a) and b) of Proposition \ref{prop:proj_pi} we get that
$B^+_{0,k}=\pi_+(B_{0,k})$. 

We claim that the operators $B^+_{0,k}$ and $B^-_{0,k}=\pi_-(B_{0,k})$
satisfy conditions (i)--(iii) in Lemma \ref{le:der_R} and hence the
extended flows $\pa_{0,k}$ exist. We divide the argument into two
parts. First, we show that conditions (i)--(iii) of Lemma
\ref{le:der_R} hold provided that the coefficient ring $\R$ is
replaced with $\hat{\R}$. Second, we will prove that the coefficients
of $B^\pm_{0,k}$ belong to $\R$, that is, the terms that involve the
differential operator $\epsilon\pa_x$ cancel out.

Let us check condition (i). The identity
\ben
B_{0,k,1}^\# = -(\Lambda-\Lambda^{-1})B_{0,k,1}
(\Lambda-\Lambda^{-1})^{-1} 
\een
is obvious.  We have
\ben
\pi_+(B_{0,k,2}) = \frac{1}{2} \iota_\Lambda \operatorname{Res}_{\pa_2} \Big(
A_{2,k} \widetilde{H}_2^{-1} \pa_2^{-1} \Big) (\Lambda-\Lambda^{-1}). 
\een
Note that $\widetilde{H}_2=S_2(\Lambda-1)S_2^{-1}$ and $A_{2,k}$
commute, and that the conjugate of $A_{2,k}$ can be computed by
formula \eqref{A_ak-symm}. Therefore,  the same argument used in the
proof of formula \eqref{B_2odd-invol} yields
\ben
(\pi_+(B_{0,k,2}))^\# =- (\Lambda-\Lambda^{-1}) \pi_-(B_{0,k,2})
\iota_\Lambda (\Lambda-\Lambda^{-1})^{-1}. 
\een
Similarly,
\ben
(\pi_+(B_{0,k,3}))^\# =- (\Lambda-\Lambda^{-1}) \pi_-(B_{0,k,3})
\iota_\Lambda (\Lambda-\Lambda^{-1})^{-1}. 
\een
Let us check condition (iii). It is sufficient to prove that each $B_{0,k,i}$
($i=1,2,3$) satisfies condition (iii) with $B^+=B^-=B_{0,k,i}$ and
$\R$ replaced by $\hat{\R}$. For $i=1$ the statement is
obvious, because $B_{0,k,1}$ is a difference operator divisible by
$\Lambda-\Lambda^{-1}$. For $i=2$ we will use the same argument as in the
proof of part c) of Lemma \ref{le:pa_2-identities} to  prove that
$\pi_+(H_a B_{0,k,2})\in \hat{\R}[\Lambda^{\pm 1}]$ for $a=1,2,3$. This
would imply that $H_a B_{0,k,2}-  \pi_+(H_a B_{0,k,2}) \in
\hat{\E}^0_{(+)}H\cap \hat{\E} = \hat{\E} H$, that is,
\beq\label{ha-bok2}
H_a B_{0,k,2}\in \hat{\R}[\Lambda^{\pm 1}] \oplus \hat{\E} H,
\eeq
so the projections $\pi_+$ and $\pi_-$ of $H_a B_{0,k,2}$ coincide
with the projection on the first factor of the above direct sum
\eqref{ha-bok2}. Let us see how the argument from the proof of part c)
of Lemma \ref{le:pa_2-identities} is modified in order to prove that
$\pi_+(H_a B_{0,k,2})\in \hat{\R}[\Lambda^{\pm 1}]$ for $a=2$. The
modification for the other two cases $a=1$ and $a=3$ is analogous. 
To begin with, we have
\beq\label{pi-0-term}
\frac{1}{2}\pi_+(H_2 (A_{2,k})_{2,[0]} (\Lambda^{-1}+1) ) =
\frac{1}{2}(\Lambda-1)\Big(\pa_2(A_{2,k})_{2,[0]} (\Lambda^{-1}+1)  +
[(A_{2,k})_{2,[0]} (\Lambda^{-1}+1),Q_2]
\Big).
\eeq
Note that
\ben
H_2 (A_{2,k})_{2,>0} = (H_2 A_{2,k})_{2,\geq 0} - H_2 \cdot (A_{2,k})_{2,[0]}
\een
Recalling part a) of Proposition \ref{prop:proj_pi} and the expansion
$\widetilde{H}_2^{-1}= (\Lambda-1)^{-1} + 2 Q_2
(\Lambda-\Lambda^{-1})^{-1} \pa_2^{-1} + \cdots$,  we get that
\begin{align}\nonumber
\pi_+\left((A_{2,k})_{2,>0} \right) & = \frac{1}{2}\iota_{\Lambda^{-1}} \Big(
H_2 A_{2,k} \widetilde{H}_2^{-1} \Big)_{2,0} (\Lambda-\Lambda^{-1}) +
  \\
  \nonumber
  &
+ \frac{1}{2} \Big(q_2 (\Lambda+1) (A_{2,k})_{2,[0]} - \pa_2
(A_{2,k})_{2,[0]} \Big)(\Lambda^{-1}+1)  -
(\Lambda-1)  (A_{2,k})_{2,[0]} Q_2.
\end{align}
Adding up the above formula and \eqref{pi-0-term}, we get that
\ben
\pi_+(H_2B_{0,k,2})=
\frac{1}{2}\iota_{\Lambda^{-1}} \Big(
H_2 A_{2,k} \widetilde{H}_2^{-1} \Big)_{2,0} (\Lambda-\Lambda^{-1}) +
\frac{1}{2} (\Lambda-1)  (A_{2,k})_{2,[0]} (\Lambda^{-1}-1) Q_2.
\een
The second term on the RHS is a difference operator, because
$(\Lambda^{-1}-1)Q_2=-q_2[-1] (\Lambda^{-1}+1)$. The first term is
also a difference operator, because
$A_{2,k}$ and $\widetilde{H}_2$ commute  and  $H_2
\widetilde{H}_2^{-1} = \pa_2 + q_2$.

The proof that $B_{0,k,3}$ satisfies condition (iii) is completely
analogous. Let us move to the last step, i.e., we will prove that
$B_{0,k,i}$ ($i=1,2,3$) satisfy condition (ii) with
$B^+=B^-=B_{0,k,i}$ and $\R$ replaced by $\hat{\R}$. Again, this fact is
obvious for $i=1$ and the arguments for $i=2$ and $i=3$ are identical,
so let us give the details only for $i=2$. We follow the same idea is
in the proof of part d of Lemma \ref{le:pa_2-identities}. Since
$[B_{0,k,2},\L]\in \hat{\E}$, it is sufficient to prove that the
projection $\pi_+([B_{0,k,2},\L])\in \hat{\R}[\Lambda,\Lambda^{-1}]$.
We have
\beq\label{comm-b0k2}
[B_{0,k,2},\L] =( [A_{2,k},\L])_{2,\geq 0} +
\frac{1}{2}[(A_{2,k})_{2,[0]} (\Lambda^{-1}-1),\L]+
\pa_2 (\resi_{\pa_2} (A_{k,2})).
\eeq
The $\pi_+$-projections of the 3rd term is difference
operators, so we have to check that the projection of the 1st and
the 2nd terms add up to a difference operator. Let us decompose $\L$ as in
\eqref{pi_2-decomp}. Since $A_{2,k}$-commutes  with $\widetilde{H}_2$
and $\widetilde{H}_1$ (see the conjugation formulas for $S_2$), we get
\ben
[A_{2,k},\L]_{2,\geq 0} = ([A_{2,k},A] \widetilde{H}_2)_{2,\geq 0} -
\Big(\pa_3(A_{2,k} \pa_2^{-1}) H_1\Big)_{2,\geq 0}.
\een
Note that
\ben
\pi_+(([A_{2,k},A] \widetilde{H}_2)_{2,>
  0})=\frac{1}{2}\iota_{\Lambda^{-1}}
\Big(([A_{2,k},A] \widetilde{H}_2)_{2,>0}
\widetilde{H}_2^{-1}\Big)_{2,[0]}\, (\Lambda-\Lambda^{-1})
\een
is a difference operator, because in the above formula we can replace
$([A_{2,k},A] \widetilde{H}_2)_{2,>0}$ by $([A_{2,k},A]
\widetilde{H}_2 - ([A_{2,k},A] \widetilde{H}_2)_{2,[0]}$, where only
the second term could have a contribution which is not a difference
operator. However, this contribution  involves only the 0-th order
term of $\widetilde{H}_2^{-1}(\Lambda-\Lambda^{-1})$ which is
$\Lambda^{-1}+1$, that is, a difference operator. Clearly,
$\pi_+(([A_{2,k},A] \widetilde{H}_2)_{2,[0]})$ is a difference
operator. Note that
\ben
\pi_+\Big(\pa_3(A_{2,k} \pa_2^{-1}) H_1\Big)_{2,\geq 0} =
\pa_3(A_{2,k} )_{2,[0]} Q_3.
\een
Recalling formula \eqref{comm-b0k2}, we get that up to terms that are
difference operators, the projection $\pi_+([B_{0,k,2},\L])$ coincides
with
\ben
-\pa_3(A_{2,k} )_{2,[0]} Q_3 +\frac{1}{4} \pi_+(
[A_{2,k} )_{2,[0]}, \pa_3^2] (\Lambda^{-1}-1)) =
-\frac{1}{2} \pa_3(A_{2,k} )_{2,[0]} (\Lambda^{-1}+1)Q_3. 
\een
However, $(\Lambda^{-1}+1)Q_3 = q_3[-1] (1-\Lambda^{-1})$ is a
difference operator, so $\pi_+([B_{0,k,2},\L])$ is also a difference
operator.

Let us prove that the coefficients of $B^\pm_{0,k}$ belong to
$\R$. We will consider only the case $B^+_{0,k}$. The other
case is similar. Since $\pi_+ ([A_{1,k},\L])=\pi_+(H_a
A^{+}_{1,k})=0$, it will be sufficient to prove that the coefficients
of the operator $B_{0,k}-A^+_{1,k}\in \hat{\E}_{(+)}$ belong to
$\R$. The coefficients of the operator $B_{0,k}-A^+_{1,k}$ are at most
1st order differential operators in $\epsilon\pa_x$. The vanishing of
the coefficient in front of $\epsilon\pa_x$ is equivalent to the
following identity
\ben
\sum_{m=0}^\infty \left(
  \tfrac{L_1^{(n-2)k}}{(n-2)^k}  \Lambda^{-2m-1} +
  \Lambda^{2m+1} \tfrac{(L_1^\#)^{(n-2)k}}{(n-2)^k}
\right)_{1,<0} (\Lambda-\Lambda^{-1}) = 
\frac{1}{2} \iota_{\Lambda^{-1}} \left(\Big( \tfrac{L_2^{2k}}{2^k} 
  \widetilde{H}_2^{-1} \Big)_{2,[0]}+
\Big( \tfrac{L_3^{2k}}{2^k} 
  \widetilde{H}_3^{-1} \Big)_{3,[0]}\right)
  (\Lambda-\Lambda^{-1}). 
  \een
  \begin{proposition}\label{prop:ext_rel}
    If $k\geq 0$ is an integer, then the following identity holds
    \ben
    \sum_{m=0}^\infty \left(
  \tfrac{L_1^{(n-2)k}}{(n-2)^k}  \Lambda^{-2m-1} +
  \Lambda^{2m+1} \tfrac{(L_1^\#)^{(n-2)k}}{(n-2)^k}
\right) =
\frac{1}{2} \sum_{m\in \ZZ} \left(\Big( \tfrac{L_2^{2k}}{2^k} 
  S_2 \Lambda^m S_2^{-1} \Big)_{2,[0]}-
\Big( \tfrac{L_3^{2k}}{2^k} 
  S_3 (-\Lambda)^m S_3^{-1}\Big)_{3,[0]}\right).
    \een
  \end{proposition}
  The above proposition yields the identity that we want to
  prove. Indeed, since
\ben
\sum_{m<0} S_2 \Lambda^m S_2^{-1} = 
\iota_{\Lambda^{-1}} S_2 (\Lambda-1)^{-1} S_2^{-1} = 
\iota_{\Lambda^{-1}} \widetilde{H}_2^{-1}
\een
and
\ben
-\sum_{m<0} S_3 (-\Lambda)^m S_3^{-1} = 
\iota_{\Lambda^{-1}} S_3 (\Lambda+1)^{-1} S_3^{-1} = 
\iota_{\Lambda^{-1}} \widetilde{H}_3^{-1},
\een
we just have to remove the terms that involve non-negative powers of
$\Lambda$ and multiply both sides by $(\Lambda-\Lambda^{-1})$. 

\subsection{Proof of Proposition \ref{prop:ext_rel}}\label{sec:ext_rel}
We are going to prove a more general statement. Namely, put
\begin{align}
  \nonumber
M_{+,k}^{r,s} & := \sum_{m=0}^\infty
                S_1 \cdot\mu_{1,k} \cdot (\pa_2^r \pa_2^s(S_1^{-1})) \cdot\Lambda^{-2m-1},\\
  \nonumber
M_{-,k}^{r,s} & := \sum_{m=0}^\infty \Lambda^{2m+1}\cdot
\Big((\pa_2^r\pa_3^s (S_1)) \cdot\mu_{1,k} \cdot S_1^{-1}\Big)^\#,\\
  \nonumber
M_{2,k}^{r,s} & := \Big(
                S_2 \cdot\mu_{2,k} \cdot\pa_3^s(S_2^{-1}) \cdot (-\pa_2)^r \Big)_{2,[0]} ,\\
   \nonumber
M_{3,k}^{r,s} & := \Big(
                S_3 \cdot\mu_{3,k} \cdot\pa_2^r(S_3^{-1}) \cdot (-\pa_3)^s \Big)_{3,[0]},
\end{align}
where
\ben
\mu_{1,k}:= \tfrac{\Lambda^{(n-2)k}}{(n-2)^k},\quad
\mu_{2,k}:=\sum_{m\in \ZZ} \tfrac{\pa_2^{2k}}{2^k}
\Lambda^m,\quad
\mu_{3,k}:=-\sum_{m\in \ZZ} \tfrac{\pa_3^{2k}}{2^k}
(-\Lambda)^m.
\een
We are going to prove the following identity
\beq\label{ext-id}
M_{+,k}^{r,s} + M_{-,k}^{r,s} = \frac{1}{2}\Big( M_{2,k}^{r,s} + M_{3,k}^{r,s}\Big).
\eeq
The identity stated in the proposition is the case when $r=s=0$. Let
us introduce the following notation: If $P(\Lambda,\pa_2,\pa_3)\in \E$
is a differential-difference operator, then we denote by $\vec{P}$ the
operator, acting on $\E_{(\pm)}$, such that $\Lambda$ acts by
multiplication, while $\pa_2$ and $\pa_3$ act by derivations, e.g.,
\ben
\vec{H}_2 (M) = (\Lambda-1) \cdot \pa_2(M) -q_2(\Lambda+1)\cdot M =
H_2\cdot M - (\Lambda-1) \cdot M \cdot \pa_2. 
\een
\begin{lemma}\label{le:hl-action}
  Suppose that $i\in \{+,-,2,3\}$ and $k,r,s\geq 0$ are arbitrary
  integers.
  
  a) The following formulas hold:
  \ben
  \vec{H}_2(M_{i,k}^{r,s}) = (\Lambda-1) M_{i,k}^{r+1,s},\quad
  \vec{H}_3(M_{i,k}^{r,s}) = (\Lambda+1) M_{i,k}^{r,s+1}.
  \een
  
  b) The following formulas hold:
  \ben
  \vec{\L}(M_{i,k}^{r,s})=M_{i,k+1}^{r,s} +\pa_2(M_{i,k}^{r+1,s} )+
  \pa_3(M_{i,k}^{r,s+1}) -\frac{1}{2} M_{i,k}^{r+2,s} -\frac{1}{2} M_{i,k}^{r,s+2}.
  \een
\end{lemma}
\proof
a) Let us prove the formulas for $\vec{H}_2(M_{+,k}^{r,s})$ and
$\vec{H}_2(M_{3,k}^{r,s})$. The argument in all other cases is
similar. We have
\ben
\vec{H}_2(M_{+,k}^{r,s}) = H_2 \cdot M_{+,k}^{r,s} -
(\Lambda-1)M_{+,k}^{r,s}\cdot \pa_2.
\een
On the other hand, $H_2 S_1 = (\Lambda-1) S_1 \pa_2$. Therefore, the
above identity is transformed into
\ben
\vec{H}_2(M_{+,k}^{r,s}) = \sum_{m=0}^\infty
(\Lambda-1) S_1 \mu_{1,k}
(\pa_2 \cdot (\pa_2^r\pa_3^s(S_1^{-1}))-
(\pa_2^r\pa_3^s(S_1^{-1}))\cdot \pa_2) =
(\Lambda-1) M_{+,k}^{r+1,s}.
\een
We have
\ben
\vec{H}_2(M_{3,k}^{r,s}) =
\Big(
H_2 S_3 \mu_{3,k}\pa_2^r(S_2^{-1})  (-\pa_3)^s  -
(\Lambda-1) S_3 \mu_{3,k}\pa_2^r(S_2^{-1})  (-\pa_3)^s  \pa_2\Big)_{3,[0]}.
\een
Recall that $H_2= (\Lambda-1)\pa_3^{-1} H_1 -\pa_3^{-1} q_2
H_3$. Recalling also the conjugation formulas for $S_3$, we get
$H_3 S_3= (\pa_3+q_3) S_3 (\Lambda+1)$ and $\pa_3^{-1}
H_1S_3=S_3\pa_2$. Finally, note that $(\Lambda+1)
\mu_{3,k}=0$. Putting these facts together we get 
\ben
\vec{H}_2(M_{3,k}^{r,s}) = (\Lambda-1)\Big(
S_3 \cdot \mu_{3,k} \cdot \pa_2\cdot \pa_2^r(S_2^{-1}) \cdot (-\pa_3)^s -
S_3 \cdot \mu_{3,k}\cdot \pa_2^r(S_2^{-1})  \cdot (-\pa_3)^s \cdot \pa_2\Big)_{3,[0]}.
\een
The RHS of the above formula coincides with $M_{2,k}^{r+1,s}$.

b) Let us prove the formula for $\vec{\L}(M_{2,k}^{r,s})$. The
remaining ones are proved with the same technique. The action of
$\vec{\L}$ can be computed as follows:
\beq\label{vecL-action}
\vec{\L}(M) = \Big(\L-\tfrac{\pa_2^2}{2}-\tfrac{\pa_3^2}{2}\Big) \cdot
M + \tfrac{\pa_2^2(M)}{2} + \tfrac{\pa_3^2(M)}{2} =
\L\cdot M -\sum_{a=2,3} \Big(
\pa_a\cdot M\cdot \pa_a -
M \cdot \tfrac{\pa_a^2}{2}\Big)
\eeq
Let us apply formula \eqref{vecL-action} to
$M:=S_2 \cdot \mu_{2,k} \cdot \pa_3^s(S_2^{-1}) (-\pa_2)^r$. The
operator $\L$ has the following decomposition (see formula
\eqref{pi_2-decomp}):
\ben
\L=\tfrac{L_2^2}{2} + A \widetilde{H}_2 +
\pa_3 \cdot \widetilde{H}_1-
\tfrac{\widetilde{H}_1^2}{2}.
\een
Recalling the conjugation formulas for $S_2$, we get
\ben
\tfrac{L_2^2}{2} \cdot M = S_2 \cdot \mu_{2,k+1} \cdot \pa_3^s(S_2^{-1}) (-\pa_2)^r,
\een
$\widetilde{H}_2 \cdot M=0$, because $\widetilde{H}_2 S_2=(\pa_2+q_2 )
S_2(\Lambda-1)$ and $(\Lambda-1)\mu_{2,k}=0$, and finally,
\ben
\Big(\pa_3 \cdot \widetilde{H}_1-
\tfrac{\widetilde{H}_1^2}{2}\Big)\cdot M =
\Big(\pa_3 \cdot S_2\cdot \pa_3- S_2 \cdot \tfrac{\pa_3^2}{2} \Big)
\cdot
\mu_{2,k} \cdot \pa_3^s(S_2^{-1}) (-\pa_2)^r
\een
Substituting these formulas in \eqref{vecL-action}, we get
\begin{align}
  \nonumber
  \vec{\L}(M) & =S_2 \cdot \mu_{2,k+1} \cdot \pa_3^s(S_2^{-1})
                (-\pa_2)^r +
      \Big(\pa_3 \cdot S_2\cdot \pa_3- S_2 \cdot \tfrac{\pa_3^2}{2} \Big)
\cdot
\mu_{2,k} \cdot \pa_3^s(S_2^{-1}) (-\pa_2)^r          
-\sum_{a=2,3} \Big(
\pa_a\cdot M\cdot \pa_a -
M \cdot \tfrac{\pa_a^2}{2}\Big) .               
\end{align}
Substituting also the formula for $M$, after a short computation we
get the following formula:
\begin{align}
  \nonumber
  \vec{\L}(M) & =S_2 \cdot \mu_{2,k+1} \cdot \pa_3^s(S_2^{-1})
                (-\pa_2)^r \\
  \nonumber
  & +\pa_3\Big(S_2\cdot \mu_{2,k} \cdot \pa_3^{s+1}(S_2^{-1}) (-\pa_2)^r \Big)
  +\pa_2\Big(S_2\cdot \mu_{2,k} \cdot \pa_3^{s}(S_2^{-1}) (-\pa_2)^{r+1} \Big)+\\
  \nonumber
  & -\frac{1}{2}\Big(S_2\cdot \mu_{2,k} \cdot \pa_3^{s+2}(S_2^{-1})
    (-\pa_2)^r \Big)
    -\frac{1}{2}\Big(S_2\cdot \mu_{2,k} \cdot \pa_3^{s}(S_2^{-1}) (-\pa_2)^{r+2} \Big).
\end{align}
Since $\vec{\L}(M_{2,k}^{r,s}) = (\vec{\L}(M))_{2,[0]}$, we see that
the 1st line of the RHS of the above formula contributes
$M_{2,k+1}^{r,s}$, the 2nd line contributes
$\pa_3(M_{2,k}^{r,s+1})+\pa_2(M_{2,k}^{r+1,s})$ and the 3rd line
contributes $-\tfrac{1}{2} M_{2,k}^{r,s+2}-\tfrac{1}{2}
M_{2,k}^{r+2,s}$. This is exactly the formula that we wanted to prove.
\qed

We prove \eqref{ext-id} by induction on the lexicographical order of
the tripple $(k,r,s)$. If $k=r=s=0$, then the identity is obvious. Let
\ben
P_{k,r,s}:=M_{+,k}^{r,s}+M_{-,k}^{r,s}-\frac{1}{2} \Big( M_{2,k}^{r,s}
+ M_{3,k}^{r,s}.
\een
Note that the operator
\ben
P_k:= \sum_{r,s=0}^\infty P_{k,r,s} (-1)^r \pa_2^{-r-1} (-1)^s \pa_3^{-s-1}
\een
coincides with
\begin{align}
  \nonumber
& \sum_{m=0}^\infty \Big(
S_1 \mu_{1,k} \pa_2^{-1} \pa_3^{-1} S_1^{-1} \Lambda^{-2m-1}+
\Lambda^{2m+1} \left(
  S_1 \mu_{1,k} \pa_2^{-1} \pa_3^{-1} S_1^{-1} \right)^\# \Big)+ \\
  \nonumber
  &
  -\frac{1}{2} \Big(S_2\mu_{2,k} \pa_3^{-1} S_2^{-1} \pa_2^{-1} \Big)_{2,<0}
-\frac{1}{2} \Big(S_3\mu_{3,k} \pa_2^{-1} S_3^{-1} \pa_3^{-1} \Big)_{3,<0}.
\end{align}
This formula implies that $P_k^\#=P_k$. Let us prove that $P_0=0$. If
$k=0$, then let us check that \eqref{ext-id} holds if $r=0$ or
$s=0$. In terms of the operator $P_0$, this statement is equivalent to
$\resi_{\pa_2}(P_0)=\resi_{\pa_3}(P_0)=0$. We have 
\ben
\resi_{\pa_2} (P_0)_{1,<0} = S_1 \pa_3^{-1} S_1^{-1}
(\Lambda-\Lambda^{-1})^{-1} -\frac{1}{2}\Big(
\pa_3^{-1} (\Lambda-1)^{-1} + S_3 (\Lambda+1)^{-1} S_3^{-1} \pa_3^{-1}
\Big),
\een
where the rational expressions in $\Lambda$ should be expanded in the
powers of $\Lambda^{-1}$. Recalling the conjugation formulas for $S_1$
and $S_3$, we get
\ben
(\pa_3-Q_3)^{-1} (\Lambda-\Lambda^{-1})^{-1} - \frac{1}{2}\Big(
\pa_3^{-1} (\Lambda-1)^{-1}  + H_3^{-1} (\pa_3+q_3)\pa_3^{-1}
\Big)=0.
\een
Since $P_0^\#=P_0$, the residue satisfies $(\resi_{\pa_2}(P_0))^\# =
-\resi_{\pa_2}(P_0)$. Therefore, $\resi_{\pa_2}(P_0)=0$ as
claimed. The vanishing of the other residue is proved in a similar
way. Suppose now that $(r_0,s_0)$ is a lexicographically minimal pair
such that $P_{0,r_0,s_0}\neq 0$. Then $r_0>0$ and $s_0>0$ by what we
have just proved. Recalling part a) of Lemma \ref{le:hl-action}, we get
\ben
(\Lambda-1)P_{0,r_0,s_0} = \vec{H}_2(P_{0,r_0-1, s_0}) =0
\een
and
\ben
(\Lambda+1)P_{0,r_0,s_0} = \vec{H}_3(P_{0,r_0, s_0-1}) =0.
\een
The first identity implies that $P_{0,r_0,s_0} = (\sum_{m\in \ZZ}
\Lambda^m ) a$ for some $a\in \R$, while the second one implies that
$P_{0,r_0,s_0}= (\sum_{m\in \ZZ} (-\Lambda)^m ) b$ for some $b\in
\R$. Comparing the coefficients in front of $\Lambda^0$ and
$\Lambda^1$, we get respectively that $a=b$ and $a=-b$. This is
possible only if $a=b=0$ -- contradiction. Therefore, $P_0=0$, that is
\eqref{ext-id} holds for $k=0$ and for all $r,s\in \ZZ_{\geq
  0}$. Recalling part b) of Lemma \ref{le:hl-action}, we get that if
\eqref{ext-id} holds for some $k$ and for all $r,s\in \ZZ_{\geq 0}$,
then it must hold for $k+1$ and for all $r,s\in \ZZ_{\geq 0}.$ This
completes the inductive step, so formula \eqref{ext-id} is proved.

\section{Commutativity of flows}
The main goal in this section is to prove part b) of Theorem
\ref{thm:Lax}.
\subsection{Evolution of the dressing operators}\label{sec:dress_evol}
\begin{proposition}\label{corevolution}
	The projections of $\mathcal{L}$ satisfy the following
        differential equations: 
	\begin{eqnarray*}
		\pa_{i,k}\pi_b(\mathcal{L})=\pi_{b}([B_{i,k},\pi_b(\mathcal{L})])=[\pi_{b}(B_{i,k}),\pi_b(\mathcal{L})],
	\end{eqnarray*}
	where $i=0,1,2,3$, $b=\pm, 2,3$ and $k$ is odd when $i=2,3$.
\end{proposition}
\proof
First, let us compute the derivatives of the operator $H_1$.  Using
the relation $H_1=\frac{1}{2}(\pa_2-q_2)\cdot 
H_3-\frac{1}{2}(\pa_3-q_3)\cdot H_2$, we get 
$\pa_{i,k}H_1\in\frac{1}{2}(\pa_2-q_2)\cdot \pa_{i,k}
H_3-\frac{1}{2}(\pa_3-q_3)\cdot \pa_{i,k} H_2+\mathcal{E}H=-H_1
B_{i,k}+\mathcal{E}H$. Therefore 
\ben
\pa_{i,k}H_1+H_1B_{i,k}\in\mathcal{E}H.
\een
The identity
$\pi_{b}([B_{i,k},\pi_{b}(\mathcal{L})])=[\pi_{b}(B_{i,k}),\pi_b(\mathcal{L})]$
follows from the definition of the projection $\pi_b$ and part b) of
Proposition \ref{prop:proj-A}.

By definition, the projection
$\pi_b(\mathcal{L})=\mathcal{L}+\sum_{i=1}^3A_iH_i$, where 
$A_i\in\mathcal{E}_{(b)}$. Recalling Theorem \ref{thm:Lax}, we get
\begin{eqnarray*}
	\pa_{i,k}\pi_b(\mathcal{L})&\in&[B_{i,k},\mathcal{L}]-\sum_{i=1}^3 A_iH_iB_{i,k}+\mathcal{E}_{(b)}H\nonumber\\
	&=& [B_{i,k},\pi_b(\mathcal{L})]-[B_{i,k},\sum_{i=1}^3A_iH_i]-\sum_{i=1}^3 A_iH_iB_{i,k}+\mathcal{E}_{(b)}H\nonumber\\
	&=& [\pi_b(B_{i,k}),\pi_b(\mathcal{L})]+\mathcal{E}_{(b)}H.
\end{eqnarray*}
Therefore, $\pa_{i,k}\pi_b(\mathcal{L})-[\pi_b(B_{i,k}),\pi_b(\mathcal{L})]\in \mathcal{E}_{(b)}^0\cap\mathcal{E}_{(b)}H=\{0\}$. 
\qed

Furthermore, recalling the definitions \eqref{projs-L}, we get the
following corollary:
\begin{corollary}\label{corevolpmla}
	The operators $L_1^{\pm}$ and $L_a$ $(a=2,3$) satisfy the
        following differential equations:
       	\begin{eqnarray*}
		&&\pa_{i,k}L_1^{\pm}=\pi_{\pm}([B_{i,k},L_1^{\pm}])=[\pi_{\pm}(B_{i,k}),L_1^{\pm}],\\
		&&\pa_{i,k}L_a=\pi_{a}([B_{i,k},L_a])=[\pi_{a}(B_{i,k}),L_a],
	\end{eqnarray*}
        where $L_1^+:=L_1$ and $L_1^-$ is defined by formula
        \eqref{L1-}.
        \qed
\end{corollary}
We would like to extend the derivations $\pa_{i,k}$ of $\R$ to
derivations of $\R_i$ ($i=1,2,3$), that is, the rings that contain the
coefficients of the dressing operators. Recall that the Lax operators
can be expressed in terms of the dressing operators as follows:
$S_1\cdot \Lambda \cdot S_1^{-1} =L_1^+$ and
$S_a\cdot\pa_a\cdot S_a^{-1}=L_a$ ($a=2,3$). Comparing these formulas
with the formulas for the derivatives of the Lax operators (see
Corollary \ref{corevolpmla}), we get that a natural choice for the
extension is given by the following formulas:

  a) The derivation $\pa_{1,k}$ ($k\geq 1$) 
  \begin{align}\nonumber
    \pa_{1,k}S_1 & =(-(L_1^+)^k+B_{1,k})\cdot S_1,\\
    \nonumber
    \pa_{1,k}S_2& =\pi_{2}(B_{1,k})\cdot S_2,\\
    \nonumber
    \pa_{1,k}S_3 & =\pi_{3}(B_{1,k})\cdot S_3.
  \end{align}
  b) The derivation $\pa_{a,2l+1}$ ($a=2,3$, $l\geq 0$) 
  \begin{align}\nonumber
    \pa_{a,2l+1}S_1 & =\pi_+(B_{a,2l+1})\cdot S_1,\\
    \nonumber
    \pa_{a,2l+1}S_a& =(-L_a^{2l+1}+B_{a,2l+1})\cdot S_a,\\
    \nonumber
    \pa_{a,2l+1} S_b& =\pi_b(B_{a,2l+1})\cdot S_b,
    \end{align}
    where $b=\{2,3\}\setminus{\{a\}}$.

    \noindent
    c) The derivation $\pa_{0,k}$ ($k\geq 1$) 
  \begin{align}\nonumber
	\pa_{0,k}S_1& =(-A_{1,k}^++\pi_{+}(B_{0,k}))\cdot
                           S_1,\\
    \nonumber
    \pa_{0,k}S_2 & =(-A_{2,k}+\pi_{2}(B_{0,k}))\cdot S_2,\\
    \nonumber
    \pa_{0,k}S_3 & =(-A_{3,k}+\pi_{3}(B_{0,k}))\cdot S_3.
  \end{align}
  Recall that the rings $\R_i$ ($i=1,2,3$) are equipped with the
  derivations $\epsilon\pa_x,\pa_2,\pa_3$ and the action of the translation
  operator $\Lambda=e^{\epsilon\pa_x}$.

  \begin{proposition}\label{corevodres}
    The following commutators in $\R_i$ ($i=1,2,3$) vanish:
    $[\pa_{j,k},\epsilon\pa_x]= [\pa_{j,k},\pa_2]=[\pa_{j,k},\pa_3]$, that is, 
    the extended derivation $\pa_{j,k}$ commutes with $\epsilon\pa_x,\pa_2,$
    and $\pa_3$. 
  \end{proposition}
  \proof
The vanishing of the commutators
  $[\pa_{j,k},\pa_x]= [\pa_{j,k},\pa_2]=[\pa_{j,k},\pa_3]=0$ in
  $\mathcal{R}_1$ can be proved in the same way as in Proposition
  \ref{commuteonR1}, a). We only prove the vanishing of the
  commutators in $\mathcal{R}_2$. The proof of the vanishing in $\mathcal{R}_3$ is
  similar.
  
Let us prove that $[\pa_{j,k},\pa_2]=0$ in $\mathcal{R}_2$. Put
$L_2=\sum_{i=0}^{+\infty} b_{2,i}\pa_2^{1-i}$ and
$S_2=\sum_{j=0}^{+\infty}\psi_{2,j}\pa_2^{-j}$, where $b_{2,0}=1$,
$b_{2,1}=0$ and $\psi_{2,0}=1$. Applying $\pa_{j,k}$ to the relation
$L_2S_2=S_2\pa_2$, we get
\begin{align*}
\pa_{j,k}(L_2)S_2+\sum_{l=0}^{+\infty}\sum_{i=0}^l\sum_{p=0}^{l-i}b_{2,i}\binom{1-i}{l-i-p}
\pa_{j,k}\left(\psi_{2,p}^{(l-i-p)}\right)\pa_2^{1-l}=\pa_{j,k}(S_2)\cdot \pa_2.
\end{align*}
Here $\psi_{2,p}^{(s)}=\pa_2^s(\psi_{2,p})$.
Since $\pa_{j,k}L_2=[\pi_2(\tilde{B}_{j,k}),L_2)]$ and
$\pa_{j,k}(S_2)=\pi_2(\tilde{B}_{j,k})S_2$, where
\begin{align*} 
\tilde{B}_{j,k}=\left\{
\begin{array}{ll}
	B_{j,k}, & \hbox{$j\neq 2$ and $j\neq 0$;} \\
	-L_2^k+B_{j,k}, & \hbox{$j=2$ and $j\neq 0$;}\\
	-A_{2,k}+B_{0,k}, & \hbox{$j=0$,}
\end{array}
\right.
\end{align*}
we get
\begin{align}
\sum_{i=0}^l\sum_{p=0}^{l-i}b_{2,i}\binom{1-i}{l-i-p}
\cdot\pa_{j,k}\left(\psi_{2,p}^{(l-i-p)}\right)=\sum_{i=0}^l\sum_{p=0}^{l-i}b_{2,i}\binom{1-i}{l-i-p}
\left(\pa_{j,k}\psi_{2,p}\right)^{(l-i-p)}.\label{pajkpa2psi2}
\end{align}
Formula (\ref{pajkpa2psi2}) for $l=2$ yields
$\pa_{j,k}(\psi_{2,1}^{(1)})=(\pa_{j,k}\psi_{2,1})^{(1)}$. Arguing by
induction on $l$ we get
$\pa_{j,k}\pa_2(\psi_{2,l})=\pa_2\pa_{j,k}(\psi_{2,l})$. Thus
$[\pa_{j,k},\pa_2]=0$ in $\mathcal{R}_2$. 

Let us prove that $[\pa_{j,k},\Lambda]=0$ and $[\pa_{j,k},\pa_3]=0$ in
$\mathcal{R}_2$. For brevity, put $A_2=(\pa_2-q_2)^{-1}\cdot
(\pa_2+q_2)$ and $A_2=\Lambda+A\cdot H_2$, where
$A=-(\pa_2-q_2)^{-1}\in\mathcal{E}_{(2)}$. Recalling Proposition
\ref{commuteonR23} and using that
$[\pa_{j,k},\pa_2]=[\Lambda,\pa_2]=0$, we get
\begin{align*}
  [\pa_{j,k},\Lambda](S_2)=&
                             \Big(\pa_{j,k}(A_2)+A_2\cdot\pi_2(\tilde{B}_{j,k})-\pi_2(\tilde{B}_{j,k})[1]\cdot A_2\Big)\cdot S_2
\end{align*}
Therefore, we only need to prove that
$\pa_{j,k}(A_2)+
A_2\cdot\pi_2(\tilde{B}_{j,k})-
\pi_2(\tilde{B}_{j,k})[1]\cdot A_2=0$.
In fact, since $\pa_{j,k}H_2=-H_2B_{j,k}+\mathcal{E}H$, we have 
\begin{align*}
\pa_{j,k}(A_2)+&A_2\cdot\pi_2(\tilde{B}_{j,k})-\pi_2(\tilde{B}_{j,k})[1]\cdot A_2\\
=& A\cdot\pa_{j,k}H_2+\Lambda\cdot\pi_2(\tilde{B}_{j,k})+A\cdot H_2\cdot\pi_2(\tilde{B}_{j,k})\\
&-\pi_2(\tilde{B}_{j,k})[1]\cdot\Lambda+\mathcal{E}_{(2)}H
\in \mathcal{E}_{(2)}H\cap \mathcal{E}_{(2)}^0=\{0\}.
\end{align*}
The vanishing of $[\pa_{j,k},\pa_3]$ in $\mathcal{R}_2$ is proved
similarly. Finally, let us prove that $[\pa_{j,k},\epsilon\pa_x]=0$ in
$\mathcal{R}_2$. Put $c_m=S_2\cdot
[\pa_{j,k},(\epsilon\pa_x)^m](S_2)$. Since $[\pa_{j,k},\Lambda]=0$ in
$\R_2$, we have $\sum_{m=1}^{+\infty}c_m/m!=0$. Using that
$[\pa_{j,k},\epsilon\pa_x]=0$ in $\mathcal{R}$,  we get the following
recursion relations:
$c_{m+1}=\tilde{\ell}_2c_m+c_1\tilde{\ell}_2^{(m)}+\epsilon\pa_xc_{m}$. 
The rest of the proof is the same as the proof of
$[\pa_2,\epsilon\pa_x]=0$ in $\mathcal{R}_2$  
(see Proposition \ref{commuteonR23}). \qed

\begin{remark} 
    Proposition \ref{corevodres} is very important, because it implies
    that the Leibniz rule holds for pseudo-differential-difference operators
    with  coefficients in $\R_i$. In particular, the following formula holds:
    \ben
    \pa_{j,k} \Big(S_i \cdot \Lambda^{m_1} \pa_2^{m_2}\pa_3^{m_3} \cdot S_i^{-1} \Big)=
    \pa_{j,k} (S_i) \cdot \Lambda^{m_1} \pa_2^{m_2}\pa_3^{m_3} \cdot  S_i^{-1} +
    S_i \cdot \Lambda^{m_1} \pa_2^{m_2}\pa_3^{m_3} \cdot \pa_{j,k} (S_i^{-1}). \qed
    \een
\end{remark}   
We will use quite frequently the following dressing formula:
$L_1^-=(S_1^-)^\# \Lambda^{-1} ((S_1^-)^{\#})^{-1}$, where
$S_1^-:= (\Lambda-\Lambda^{-1}) S_1^{-1} \iota_{\Lambda^{-1}}
(\Lambda-\Lambda^{-1})^{-1} .$ Using the formulas for the derivatives
of $S_1$, we get that the following formulas for the derivatives of the dressing
operator $S_1^-$:
\begin{align}
  \nonumber
  \pa_{1,k}(S_1^{-})^\# & =((L_1^-)^k+B_{1,k})\cdot (S_1^{-})^\#,\\
  \nonumber
  \pa_{a,2l+1}(S_1^{-})^\# & =\pi_-(B_{a,2l+1})\cdot
                             (S_1^{-})^\#,\\
  \nonumber
  \pa_{0,k}(S_1^{-})^\# & =(A_{1,k}^-+\pi_{-}(B_{0,k}))\cdot
                            (S_1^{-})^\#, 
\end{align}
where $k\geq 1$ and $a=2,3$. The operator series 
$A_{1,l}^{\pm}$ and $A_{a,l}$ (see Section \ref{sec:Lax_eqns}) can be
expressed in terms of the dressing operators, which gives immediately
the following corollary: 
\begin{corollary}\label{corAl}
	The operator series $A_{1,l}^{\pm}$ and $A_{a,l}$ satisfy the
        following differential equations: 
	\begin{eqnarray*}
		&&\pa_{i,k}A_{1,l}^{\pm}=\pi_{\pm}([B_{i,k},A_{1,l}^{\pm}])=[\pi_{\pm}(B_{i,k}),A_{1,l}^{\pm}],\\
		&&\pa_{i,k}A_{a,l}=\pi_{a}([B_{i,k},A_{a,l}])=[\pi_{a}(B_{i,k}),A_{a,l}].\qed
	\end{eqnarray*}
\end{corollary}
In the next two Lemmas we derive formulas for the operator series
$B_{1,k}$ and $B_{0,k}$, which will be needed for the proof of Theorem
\ref{thm:Lax}, b).
\begin{lemma}\label{b1klemma}
  The following formula holds:
	\begin{eqnarray*}
          B_{1,k}=(L_1^+)^k_{1,\geq 1}-(L_1^-)^k_{1,\leq -1}-
          \left.\Big(
          (L_1^+)^k_{1,\geq 1}-(L_1^-)^k_{1,\leq -1}\Big)\right|_{\Lambda=1}
	\end{eqnarray*}
\end{lemma}
\proof
Recall the definition of $B_{1,k}$ from Section
\ref{sec:Lax_eqns}. Note that 
\begin{eqnarray*}
(B_{1,k})_{1,\geq 1}=\left(\left((L_1^+)^k_{1,\geq 1}\sum_{m=0}^{+\infty}\Lambda^{-2m-1}\right)_{1,\geq 0}\cdot(\Lambda-\Lambda^{-1})\right)_{1,\geq 1}
\end{eqnarray*}
Let us rewrite the term $\left((L_1^+)^k_{1,\geq
    1}\sum_{m=0}^{+\infty}\Lambda^{-2m-1}\right)_{1,\geq 0}$ in the
form $(\quad)-(\quad)_{1,<0}$. Since
$\left(A_{1,<0}\cdot (\Lambda-\Lambda^{-1})\right)_{1,\geq 1}=0$ for
any operator $A$, we get $(B_{1,k})_{1,\geq 1}=(L_1^+)^k_{1,\geq
  1}$. Similarly, using the conjugation relation
\ben
\Big( (L_1^+)^k\sum_{m=0}^{+\infty}\Lambda^{-2m-1}\Big) ^\#=
\Big( (L_1^-)^k\sum_{m=0}^{+\infty}\Lambda^{2m+1}\Big),
\een 
we get $(B_{1,k})_{1,\leq -1}=-(L_1^-)^k_{1,\leq -1}$.
Finally, since $\left.B_{1,k}\right|_{\Lambda=1}=0$, we get
$\left.\left((L_1^+)^k_{1,\geq 1}-(L_1^-)^k_{1,\leq
    -1}\right)\right|_{\Lambda=1}+\left(B_{1,k}\right)_{1,[0]}=0$. Solving
for $\left(B_{1,k}\right)_{1,[0]} $, we get that $B_{1,k}=
(B_{1,k})_{1,\geq 1} +\left(B_{1,k}\right)_{1,[0]} + (B_{1,k})_{1,\leq
  -1}$ coincides with the  RHS of the formula that we wanted to prove.
\qed

\begin{lemma} \label{b0klemma} a) The following formula holds:
	\begin{eqnarray}
	B_{0,k}=(A_{1,k}^+)_{1,\geq 1}-(A_{1,k}^{-})_{1,\leq -1}+(A_{2,k})_{2,\geq 1}+(A_{3,k})_{3,\geq 1}+(B_{0,k})_{[0]},
	\end{eqnarray}
	where
        \begin{align}
          \nonumber
          (B_{0,k})_{[0]}= & -\left.\Big(
                             (A_{1,k}^+)_{1,\geq1}-(A_{1,k}^{-})_{1,\leq
                             -1}\Big)\right|_{\Lambda=1}+(A_{2,k})_{2,[0]} \\
          \nonumber
          = &  -\left.\Big(
              (A_{1,k}^+)_{1,\geq 1}-(A_{1,k}^{-})_{1,\leq -1}
              \Big)\right|_{\Lambda=-1}+(A_{3,k})_{3,[0]}.
        \end{align}

        b)	The following formula holds:
        \ben
        B_{0,k}=B_{k,\Lambda}+(A_{2,k})_{2,\geq
          0}+(A_{3,k})_{3,\geq 1}=B_{k,-\Lambda}+(A_{2,k})_{2,\geq
          1}+(A_{3,k})_{3,\geq 0},
        \een
        where
        \begin{align}
          \nonumber
          B_{k,\Lambda}= & (A_{1,k}^+)_{1,\geq 1}-
                           (A_{1,k}^-)_{1,\leq -1}-
                           \left.\Big((A_{1,k}^+)_{1,\geq
                           1}-(A_{1,k}^-)_{1,\leq -1}\Big)\right|_{\Lambda=1} \\
          \nonumber
          B_{k,-\Lambda}= & (A_{1,k}^+)_{1,\geq 1}-
                            (A_{1,k}^-)_{1,\leq -1}-
                            \left.\Big((A_{1,k}^+)_{1,\geq 1}-(A_{1,k}^-)_{1,\leq
                            -1}\Big)\right|_{\Lambda=-1}.
        \end{align}

c) The projections of $B_{0,k}$ are given by the following formulas:
\ben
\pi_2(B_{0,k})_{2,\geq 0}=(A_{2,k})_{2,\geq 0},\quad 
\pi_3\left(\left. B_{0,k}\right|_{\Lambda\rightarrow-\Lambda}\right)_{3,\geq
  0}=(A_{3,k})_{3,\geq 0}.
\een
\end{lemma}
\proof
a) 
Recalling the definition of $B_{0,k}$, we get $(B_{0,k})_{a,\geq
  1}=(A_{a,k})_{a,\geq 1}$ for $a=2,3$.  Using the identity from
Proposition \ref{prop:ext_rel} we get $(B_{0,k})_{1,\geq 
  1}=(A^+_{1,k})_{1,\geq 1}$ and $(B_{0,k})_{1,\leq
  -1}=-(A^-_{1,k})_{1,\leq -1}$, where both identities are proved in
the same ways as the corresponding identities for $B_{1,k}$ were
proved in Lemma \ref{b1klemma}. Note that
\begin{eqnarray*}
B_{0,k}|_{\Lambda=1}=(A_{2,k})_{2,[0]},\quad B_{0,k}|_{\Lambda=-1}=(A_{3,k})_{2,[0]}.
\end{eqnarray*}
The formula for $(B_{0,k})_{[0]}$ follows easily.

b) This is just a reformulation of a).

c) Follows again from a) and the following relations:
\begin{align}
  \nonumber
  \pi_2(\Lambda^j)= & 1+O(\pa_2^{-1})+\mathcal{E}_{(2)}H_2,
  \\
  \nonumber
  \pi_3(\Lambda^j)= & (-1)^j+O(\pa_3^{-1})+\mathcal{E}_{(3)}H_3.
     \qed
\end{align}

Next proposition is the main result of this section
\begin{proposition}\label{zsequation}
	The following Zakhalov-Shabat equations hold
	\begin{eqnarray}
	\pa_{a,k}B_{b,l}-\pa_{b,l}B_{a,k}+[B_{b,l},B_{a,k}]\in \mathcal{A}H,\label{paakbblineh}
	\end{eqnarray}
	where $a,b\in\{0,1,2,3\}$, $k$ (resp. $l$) is odd when $a=2,3$ (resp.
        $b=2,3$), $\mathcal{A}=\mathcal{E}$ if $ab\neq 0$, and
        $\mathcal{A}=\mathcal{\hat E}$ if 
        $ab=0$.
\end{proposition}
\subsection{Proof of the case $a=1$ and $b=1$}
In this case, we have to prove that 
\begin{eqnarray}
\pa_{1,k}B_{1,l}-\pa_{1,l}B_{1,k}+[B_{1,l},B_{1,k}]=0.\label{zs11}
\end{eqnarray}
Let us substitute the formulas for $B_{1,k}$ and $B_{1,l}$ from Lemma
\ref{b1klemma} in the projection $\left(\text{LHS of
    (\ref{zs11})}\right)_{1,\geq
  1}=\left[B_{1,k},(L_1^+)^l\right]_{1,\geq
  1}-\left[B_{1,l},(L_1^+)^k\right]_{1,\geq
  1}+[B_{1,l},B_{1,k}]_{1,\geq 1}$. The result can be written as the
sum of the following three parts:
\begin{eqnarray*}
	&\bullet& a_{11}= \left[\left((L_1^+)^k\right)_{1,\geq 1},(L_1^+)^l\right]_{1,\geq 1}-\left[\left((L_1^+)^l\right)_{1,\geq 1},(L_1^+)^k\right]_{1,\geq 1}+\left[\left((L_1^+)^l\right)_{1,\geq 1},\left((L_1^+)^k\right)_{1,\geq 1}\right],\\
	&\bullet& b_{11}= -\left[\left((L_1^+)^k\right)_{1,\leq -1},(L_1^+)^l\right]_{1,\geq 1}+\left[\left((L_1^+)^l\right)_{1,\leq -1},(L_1^+)^k\right]_{1,\geq 1}\\
	&&\quad\quad-\left[\left((L_1^+)^l\right)_{1,\leq -1},\left((L_1^+)^k\right)_{1,\geq 1}\right]_{1,\geq 1}-\left[\left((L_1^+)^l\right)_{1,\geq 1},\left((L_1^+)^k\right)_{1,\leq -1}\right]_{1,\geq 1},\\
	&\bullet&
                  c_{11}=\left[\left(B_{1,k}\right)_{1,[0]},(L_1^+)^l\right]_{1,\geq 1}-\left[\left(B_{1,l}\right)_{1,[0]},(L_1^+)^k\right]_{1,\geq 1}+ \\
	&&\quad\quad+\left[\left(B_{1,l}\right)_{1,[0]},\Big((L_1^+)^k\Big)_{1,\geq 1}\right]-\left[\left((L_1^+)^l\right)_{1,\geq 1},\left(B_{1,k}\right)_{1,[0]}\right].
\end{eqnarray*}
Note that $b_{11}=c_{11}=0$ and 
\begin{eqnarray*}
	a_{11}&=&\left[\left((L_1^+)^k\right)_{1,\geq 1},\left((L_1^+)^l\right)_{1,\leq 0}\right]_{1,\geq 1}-\left[\left((L_1^+)^l\right)_{1,\geq 1},(L_1^+)^k\right]_{1,\geq 1}\\
	&=&\left[(L_1^+)^k,\left((L_1^+)^l\right)_{1,\leq 0}\right]_{1,\geq 1}-\left[\left((L_1^+)^l\right)_{1,\geq 1},(L_1^+)^k\right]_{1,\geq 1}=\left[(L_1^+)^k,(L_1^+)^l\right]_{1,\geq 1}=0,
\end{eqnarray*}
where we have used that $\left[A_{1,\geq 1},B_{1,\geq 1}\right]=\left[A_{1,\geq 1},B_{1,\geq 1}\right]_{1,\geq 1}$
and $\left[A_{1,\leq 0},B_{1,\leq 0}\right]_{1,\geq 1}=0$ for any two operators $A$ and $B$.
Therefore $\Big(\text{LHS of (\ref{zs11})}\Big)_{1,\geq
  1}=0$. Similarly, one can show $\Big(\text{LHS of
  (\ref{zs11})}\Big)_{1,\leq -1}=0$. Finally, for the zero-order term
of (\ref{zs11}), using Lemma \ref{b1klemma}, we get
$$\left(\text{LHS of (\ref{zs11})}\right)_{1,[0]}=-\left(\left(\text{LHS of (\ref{zs11})}\right)_{1,\geq 1}+\left(\text{LHS of (\ref{zs11})}\right)_{1,\leq-1}\right)(1)=0,$$
where we have used the relation
$[B_{1,l},B_{1,k}]_{1,[0]}=-\left([B_{1,l},B_{1,k}]_{1,\geq
    1}+[B_{1,l},B_{1,k}]_{1,\leq -1}\right)(1)$, which is obtained
from $[B_{1,l},B_{1,k}](1)=0$. Summarize the the above results, we get
(\ref{zs11}).  
\subsection{Proof of the case $a=1$ and $b=2,3$}
Suppose that $a=1$ and $b=2$. The case when $b=3$ is similar. We have
to prove the following relation:
\begin{eqnarray}
\pa_{1,k}B_{2,l}-\pa_{2,l}B_{1,k}+[B_{2,l},B_{1,k}]\in\mathcal{E}H.
\end{eqnarray}
According to Corollary \ref{corevolpmla} and Lemma \ref{b1klemma}, we
have 
\begin{eqnarray*}
&&\pa_{2,l}\Big(B_{1,k}\Big)_{1,\geq 1}= \pi_+\left([B_{2,l},B_{1,k}]\right)_{1,\geq 1},\nonumber\\
&&\pa_{2,l}\Big(B_{1,k}\Big)_{1,\leq -1}= \pi_-\left([B_{2,l},B_{1,k}]\right)_{1,\leq -1},\nonumber\\
&&\pa_{1,k}B_{2,l}=\pi_2 \left([B_{1,k},B_{2,l}]\right)_{2,\geq 1},\label{pa2l1kb2b1}
\end{eqnarray*}
where we used that $\pi_+(\pa_2^{k}\Lambda^{-l})\in (\mathcal{E}_{(+)}^0)_{1,\leq 0}$,
$\pi_-(\pa_2^{k}\Lambda^{l})\in (\mathcal{E}_{(-)}^0)_{1,\geq 0}$,
and $\pi_2(\pa_2^{-k}\Lambda^{p})\in (\mathcal{E}_{(2)}^0)_{1,\leq 0}$
for $k,l\geq 0, p\in \mathbb{Z}$.  Moreover, using the relation
$\Lambda \pa_2=\pa_2+q_2\Lambda+q_2+H_2$, we can remove the terms
involving $\Lambda\pa_2$ in $[B_{2,l},B_{1,k}]$  and obtain
$$[B_{2,l},B_{1,k}]=h+
\pi_+\left([B_{2,l},B_{1,k}]\right)_{1,\geq 1}+
\pi_-\left([B_{2,l},B_{1,k}]\right)_{1,\leq -1}-
\pi_2 \left([B_{2,l},B_{1,k}]\right)_{2,\geq 1}+\mathcal{E}H_2,
$$
where $h$ is some function. Therefore, we only need to prove that $\pa_{2,l}\left(B_{1,k}\right)_{1,[0]}=h$.
In fact, note that
$[B_{2,l},B_{1,k}]|_{\Lambda\rightarrow-\Lambda}(1)=H_2|_{\Lambda\rightarrow-\Lambda}(1)=0$.
Therefore, $h=-\left(\pi_+\left([B_{2,l},B_{1,k}]\right)_{1,\geq
    1}+\pi_-\left([B_{2,l},B_{1,k}]\right)_{1,\leq
    -1}\right)(1)$. Recalling Lemma \ref{b1klemma}, we get
$\pa_{2,l}\left(B_{1,k}\right)_{1,[0]}=h$. 

\subsection{Proof of the case $a, b=2 $ or $3$}

Firstly, in the case $a=b=2$ or $3$, one only needs to check
$\pa_{a,k}B_{a,l}-\pa_{a,l}B_{a,k}+[B_{a,l},B_{a,k}]=0$. The proof is
very standard just like the KP case \cite{Di}.  

Then in the case $a=2$ and $b=3$, the goal is to show
\begin{eqnarray}
\pa_{2,k}B_{3,l}-\pa_{3,l}B_{2,k}+[B_{3,l},B_{2,k}]\in \mathcal{E}H\label{b2kb3lineh}
\end{eqnarray}
According to Corollary \ref{corevolpmla} and the relations
$\pa_a^{k+1}\pa_b^{-l}\in
(\mathcal{E}_{(b)}^0)_{b,<0}+\mathcal{E}_{(b)}H$, $k,l\geq 0$, $a\neq
b$, we get $\pa_{2,k}B_{3,l}=(\pa_{2,k} L_3^l)_{3,\geq
  1}=[\pi_3(B_{2,k}),L_3^l]_{3,\geq 1}=\pi_3([B_{2,k},L_3^l])_{3,\geq
  1}=\pi_3([B_{2,k},B_{3,l}])_{3,\geq 1}$. Similarly,
$\pa_{3,l}B_{2,k}=\pi_2([B_{3,l},B_{2,k}])_{2,\geq 1}$. Furthermore,
let us rewrite $[B_{3,l},B_{2,k}]=\sum_{i,j\geq
  1}a_{ij}\pa_2^i\pa_3^j$ into $h+\pi_2([B_{3,l},B_{2,k}])_{2,\geq
  1}+\pi_3([B_{3,l},B_{2,k}])_{3,\geq 1}+\mathcal{E}H$ by using the
relation $\pa_2\pa_3=H_1-q_1$. We only need to prove that $h=0$. On
the other hand, recalling part b) of Proposition \ref{prop:proj-A} and
part b) of Lemma \ref{le:pa_2-identities}, we get
$h=
\left(\pi_2([B_{2,k},B_{3,l}])\right)_{2,[0]}=
\left(\pi_2([L_2^k,B_{3,l}])\right)_{2,[0]}=
[L_2^k,\pi_2(B_{3,l})]_{2,[0]}=
-(\pa_{3,l}L_2^k)_{2,[0]}=0$. 
\subsection{Proof of the case $a=0$ and $b=2,3$}
Suppose that $a=0$ and $b=2$. The case when $b=3$ is similar. We have
to prove that 
\begin{eqnarray}
\pa_{0,k}B_{2,l}-\pa_{2,l}B_{0,k}+[B_{2,l},B_{0,k}]\in\mathcal{\hat E}H.\label{zs02}
\end{eqnarray}
According to Corollary \ref{corAl} and Lemma \ref{b0klemma}, $\pi_+\left(\text{LHS of (\ref{zs02})}\right)_{1,\geq 1}=-\left(\pa_{2,l}A_{1,k}^+\right)_{1,\geq 1}+\pi_+\left([B_{2,l},B_{0,k}]\right)_{1,\geq 1}=-\pi_+\left([B_{2,l},A_{1,k}^+]\right)_{1,\geq 1}+\pi_+\left([B_{2,l},(A_{1,k}^+)_{1,\geq 1}]\right)_{1,\geq 1}=0$. Here $\pi_{+}\left([B_{2,\geq 1},A_{1,\leq 0}]\right)_{1,\geq 1}=\pi_{+}\left([B_{2,\geq 1},A_{b,\geq 0}]\right)_{1,\geq 1}=0$ for any operators $A$ and $B$. Similarly, $\pi_-\left(\text{LHS of (\ref{zs02})}\right)_{1,\leq-1}=0$.

Using Lemma \ref{b0klemma}, we get that $\pi_2\left(\text{LHS of
    (\ref{zs02})}\right)_{2,\geq 0}$ can be written as a sum of the
following three parts: 
\begin{eqnarray*}
	a_{02}&=&[(A_{2,k})_{2,\geq 0},L_2^l]_{2,\geq 0}-[B_{2,l},A_{2,k}]_{2,\geq 0}+[B_{2,l},(A_{2,k})_{2,\geq 0}]_{2,\geq 0}=[A_{2,k},L_2^l]_{2,\geq 0}=0,\\
	b_{02}&=&\pi_2([B_{k,\Lambda},L_2^l])_{2,\geq 0}+\pi_2([B_{2,l},B_{k,\Lambda}])_{2,\geq 0}=0,\\
	c_{02}&=&\pi_2([(A_{3,k})_{3,>0},L_2^l])_{2,\geq 0}+\pi_2([B_{2,l},(A_{3,k})_{3,>0}])_{2,\geq 0}=0.
\end{eqnarray*}
Therefore $\pi_2\Big(\text{LHS of (\ref{zs02})}\Big)_{2,\geq 0}=0$. As for the projection of $\pi_3$, one can obtain $\pi_3\left(\text{LHS of (\ref{zs02})}\right)_{3,\geq 1}=-\pa_{2,l}(A_{3,k})_{3,\geq 1}+\pi_3([B_{2,l},B_{0,k}])_{3,\geq 1}=-\pi_3([B_{2,l},A_{3,k}])_{3,\geq 1}+\pi_3([B_{2,l},(A_{3,k})_{3,\geq 1}])_{3,\geq 1}=0$.

Note that by the relation $\Lambda\pa_2=\pa_2+q_2\cdot\Lambda+q_2+H_2$ and $\pa_2\pa_3=-q_1+H_1$, one can obtain $[B_{2,l},B_{0,k}]=h+\pi_+([B_{2,l},B_{0,k}])_{1,\geq 1}+\pi_-([B_{2,l},B_{0,k}])_{1,\leq -1}++\pi_2([B_{2,l},B_{0,k}])_{2,\geq 1}+\pi_3([B_{2,l},B_{0,k}])_{3,\geq 1}+\mathcal{\hat E}H$.
Thus by the same method as the one in Lemma \ref{b0klemma}, $\pi_2([B_{2,l},B_{0,k}])_{2,[0]}=h+\pi_+([B_{2,l},B_{0,k}])_{1,\geq 1}(1)+\pi_-([B_{2,l},B_{0,k}])_{1,\leq -1}(1).$
Therefore
\begin{eqnarray*}
[B_{2,l},B_{0,k}]&=&\pi_+([B_{2,l},B_{0,k}])_{1,\geq 1}+\pi_-([B_{2,l},B_{0,k}])_{1,\leq -1}\nonumber\\
&&-\Big(\pi_+([B_{2,l},B_{0,k}])_{1,\geq 1}+\pi_-([B_{2,l},B_{0,k}])_{1,\leq -1}\Big)(1)\nonumber\\
&&+\pi_2([B_{2,l},B_{0,k}])_{2,\geq 0}+\pi_3([B_{2,l},B_{0,k}])_{3,\geq 1}+\mathcal{\hat E}H. 
\end{eqnarray*}
Finally, it remains only to check that
$-\pa_{2,l}(B_{k,\Lambda})_{1,[0]}-\left(\pi_+([B_{2,l},B_{0,k}])_{1,\geq
    1}+\pi_-([B_{2,l},B_{0,k}])_{1,\leq -1}\right)(1)=0$, which is
correct thanks to Lemma \ref{b0klemma}.

\subsection{Proof of the case $a=0$ and $b=1$}
The task is to prove
\begin{eqnarray}
\pa_{0,k}B_{1,l}-\pa_{1,l}B_{0,k}+[B_{1,l},B_{0,k}]\in\mathcal{\hat E}H.\label{zs01}
\end{eqnarray}
By Corollary \ref{corAl}, Lemmas \ref{b1klemma} and \ref{b0klemma},
one can rewrite $\pi_+\left(\text{LHS of (\ref{zs01})}\right)_{1,\geq
  1}$ as the sum of the following terms:
\begin{eqnarray*}
	a_{01}&=&[(A_{1,k}^+)_{1,\geq 1},(L_1^+)^l]_{1,\geq 1}-[((L_1^+)^l)_{1,\geq 1}, A_{1,k}^+]_{1,\geq 1}+[((L_1^+)^l)_{1,\geq 1}, (A_{1,k}^+)_{1,\geq 1}]\\
	&=&[A_{1,k}^+,(L_1^+)^l]_{1,\geq 1}=0,\\
	b_{01}&=&-[(A_{1,k}^-)_{1,\leq -1},(L_1^+)^l]_{1,\geq 1}+[((L_1^-)^l)_{1,\leq -1},A_{1,k}^+]_{1,\geq 1}\\
	&&-[((L_1^-)^l)_{1,\leq -1},(A_{1,k}^+)_{1,\geq 1}]_{1,\geq 1}-[((L_1^+)^l)_{1,\geq1},(A_{1,k}^-)_{1,\leq -1}]=0,\\
	c_{01}&=& [(B_{k,\Lambda})_{1,[0]},(L_1^+)^l]_{1,\geq 1}-[(B_{1,l})_{1,[0]}, A_{1,k}^+]_{1,\geq 1}\\
	&&+[(B_{1,l})_{1,[0]},(A_{1,k}^+)_{1,\geq 1}]_{1,\geq 1}+[((L_1^+)^l)_{1,\geq1},(B_{k,\Lambda})_{1,[0]}]=0,\\
	d_{01}&=& \pi_{+}([(A_{2,k})_{2,\geq 0},(L_1^+)^l])_{1,\geq 1}+\pi_{+}([((L_1^+)^l)_{1,\geq 1},(A_{2,k})_{2,\geq 0}])_{1,\geq 1}=0,\\
	e_{01}&=& \pi_{+}([(A_{3,k})_{3,\geq 1},(L_1^+)^l])_{1,\geq 1}+\pi_{+}([((L_1^+)^l)_{1,\geq 1},(A_{3,k})_{3,\geq 1}])_{1,\geq 1}=0.
\end{eqnarray*}
So $\pi_+\Big(\text{LHS of (\ref{zs01})}\Big)_{1,\geq 1}=0$. And similarly $\pi_-\Big(\text{LHS of (\ref{zs01})}\Big)_{1,\leq -1}=0$. In order to discuss $\pi_2$ and $\pi_3$, one has to know the structure of $[B_{1,l},B_{0,k}]$. In fact, 
similar to the case of $a=0$ and $b=2$, 
\begin{eqnarray*}
[B_{1,l},B_{0,k}]&=&\pi_+([B_{1,l},B_{0,k}])_{1,\geq 1}+\pi_-([B_{1,l},B_{0,k}])_{1,\leq -1}\nonumber\\
&&-\Big(\pi_+([B_{1,l},B_{0,k}])_{1,\geq 1}+\pi_-([B_{1,l},B_{0,k}])_{1,\leq -1}\Big)(1)\nonumber\\
&&+\pi_2([B_{1,l},B_{0,k}])_{2,\geq 0}+\pi_3([B_{1,l},B_{0,k}])_{3,\geq 1}+\mathcal{\hat E}H. 
\end{eqnarray*}
Thus $\pi_2\left(\text{LHS of (\ref{zs01})}\right)_{2,\geq
  0}=-\pi_2([B_{1,l},A_{2,k}])_{2,\geq
  0}+\pi_2([B_{1,l},(A_{2,k})_{2,\geq 0}])_{2,\geq 0}=0$. Similarly,
$\pi_3\left(\text{LHS of (\ref{zs01})}\right)_{3,\geq 1}=0.$ 

Finally,
$\pa_{0,k}(B_{1,l})_{1,[0]}-\pa_{1,l}(B_{k,\Lambda})_{1,[0]}$-$\left(\pi_+([B_{1,l},B_{0,k}])_{1,\geq
    1}+\pi_-([B_{1,l},B_{0,k}])_{1,\leq -1}\right)(1)=0$ can be proved by Lemmas \ref{b1klemma} and \ref{b0klemma}.
Summarize the above results completes the proof of (\ref{zs01}).

\subsection{Proof of the case $a=0$ and $b=0$}
We have to prove that 
\begin{eqnarray}
\pa_{0,k}B_{0,l}-\pa_{0,l}B_{0,k}+[B_{0,l},B_{0,k}]\in\mathcal{\hat E}H.\label{zs00}
\end{eqnarray}
Using Corollary \ref{corAl} and Lemma \ref{b0klemma}, we can write
$\pi_+\left(\text{LHS of (\ref{zs00})}\right)_{1,\geq 1}$ as the sum 
of the following terms:
\begin{eqnarray*}
	a_{001}&=&[(A_{1,k}^+)_{1,\geq 1},A_{1,l}^+]_{1,\geq 1}-[(A_{1,k}^+)_{1,\geq 1}, A_{1,l}^+]_{1,\geq 1}+[(A_{1,l}^+)_{1,\geq 1}, (A_{1,k}^+)_{1,\geq 1}]\\
	&=&[A_{1,k}^+,A_{1,l}^+]_{1,\geq 1}=0,\\
	b_{001}&=&-[(A_{1,k}^-)_{1,\leq -1},A_{1,l}^+]_{1,\geq 1}+[(A_{1,l}^-)_{1,\leq -1},A_{1,k}^+]_{1,\geq 1}\\
	&&-[(A_{1,l}^-)_{1,\leq -1},(A_{1,k}^+)_{1,\geq 1}]_{1,\geq 1}-[(A_{1,l}^+)_{1,\geq1},(A_{1,k}^-)_{1,\leq -1}]=0,\\
	c_{001}&=& [(B_{k,\Lambda})_{1,[0]},A_{1,l}^+]_{1,\geq 1}-[(B_{l,\Lambda})_{1,[0]}, A_{1,k}^+]_{1,\geq 1}\\
	&&+[(B_{l,\Lambda})_{1,[0]},(A_{1,k}^+)_{1,\geq 1}]_{1,\geq 1}+[(A_{1,l}^+)_{1,\geq1},(B_{k,\Lambda})_{1,[0]}]=0,\\
	d_{001}&=& \pi_{+}([(A_{2,k})_{2,\geq 0},A_{1,l}^+])_{1,\geq 1}-\pi_{+}([(A_{2,l})_{2,\geq 0},A_{1,k}^+])_{1,\geq 1}\\
	&&+\pi_{+}([(A_{1,l}^+)_{1,\geq 1},(A_{2,k})_{2,\geq 0}])_{1,\geq 1}+\pi_{+}([(A_{2,l})_{2,\geq 0},(A_{1,k}^+)_{1,\geq 1}])_{1,\geq 1}=0,\\
	e_{001}&=& \pi_{+}([(A_{3,k})_{3,\geq 1},A_{1,l}^+])_{1,\geq 1}-\pi_{+}([(A_{3,l})_{3,\geq 1},A_{1,k}^+])_{1,\geq 1}\\
	&&+\pi_{+}([(A_{1,l}^+)_{1,\geq 1},(A_{3,k})_{3,\geq 1}])_{1,\geq 1}+\pi_{+}([(A_{3,l})_{3,\geq 1},(A_{1,k}^+)_{1,\geq 1}])_{1,\geq 1}=0.
\end{eqnarray*}
Similarly, $\pi_-\left(\text{LHS of (\ref{zs00})}\right)_{1,\leq
  -1}=0$. Next, we have $\pi_2\left(\text{LHS of
    (\ref{zs00})}\right)_{2,\geq 0}=a_{002}+b_{002}+c_{002}$, where
\begin{eqnarray*}
	a_{002}&=&[(A_{2,k})_{2,\geq 0},A_{2,l}]_{2,\geq 0}-[(A_{2,l})_{2,\geq 0},A_{2,k}]_{2,\geq 0}+[(A_{2,l})_{2,\geq 0},(A_{2,k})_{2,\geq 0}]_{2,\geq 0}\\
	&=&[A_{2,k},A_{2,l}]_{2,\geq 0}=0,\\
	b_{002}&=&\pi_2([B_{k,\Lambda},A_{2,l}])_{2,\geq 0}-\pi_2([B_{l,\Lambda},A_{2,k}])_{2,\geq 0}\\
	&&+\pi_2([(A_{2,l})_{2,\geq 0},B_{k,\Lambda}])_{2,\geq 0}+\pi_2([B_{l,\Lambda},(A_{2,k})_{2,\geq 0}])_{2,\geq 0}=0,\\
	c_{002}&=&\pi_2([(A_{3,k})_{3,>0},A_{2,l}])_{2,\geq 0}-\pi_2([(A_{3,l})_{3,>0},A_{2,k}])_{2,\geq 0}\\
	&&+\pi_2([(A_{2,l})_{2,\geq 0},(A_{3,k})_{3,>0}])_{2,\geq 0}+\pi_2([(A_{3,l})_{3,>0},(A_{2,k})_{2,\geq 0}])_{2,\geq 0}=0.
\end{eqnarray*}
Similarly, we have $\pi_3\left(\text{LHS of
    (\ref{zs00})}\right)_{3,\geq 1}=0$. Just as in the case of $a=0$
and $b=2$, we have
\begin{eqnarray*}
[B_{0,l},B_{0,k}]&=&\pi_+([B_{0,l},B_{0,k}])_{1,\geq 1}+\pi_-([B_{0,l},B_{0,k}])_{1,\leq -1}\nonumber\\
&&-\left(\pi_+([B_{0,l},B_{0,k}])_{1,\geq 1}+\pi_-([B_{0,l},B_{0,k}])_{1,\leq -1}\right)(1)\nonumber\\
&&+\pi_2([B_{0,l},B_{0,k}])_{2,\geq 0}+\pi_3([B_{0,l},B_{0,k}])_{3,\geq 1}+\mathcal{\hat E}H. 
\end{eqnarray*}
Finally,
$\pa_{0,k}(B_{l,\Lambda})_{1,[0]}-\pa_{0,l}(B_{k,\Lambda})_{1,[0]}-\left(\pi_+([B_{0,l},B_{0,k}])_{1,\geq
    1}+\pi_-([B_{0,l},B_{0,k}])_{1,\leq -1}\right)(1)=0$ can be
obtained by Lemma \ref{b0klemma} and the above results, which
completes the proof of \eqref{zs00}. This is the last possible case,
so the proof of Proposition \ref{zsequation} is completed. 
\subsection{Commutativity of $\pa_{a,k}$ and $\pa_{b,l}$}
Proposition \ref{zsequation}, part a) of Theorem \ref{thm:Lax}, and
Proposition \ref{prop:ker_pi_a}, yield the following corollary: 
\begin{corollary}\label{projzsequation}
Under projections $\pi_c$,
\begin{eqnarray}
\pa_{a,k}\pi_c(B_{b,l})-\pa_{b,l}\pi_c(B_{a,k})+[\pi_c(B_{b,l}),\pi_c(B_{a,k})]=0.
\end{eqnarray}
Here $a,b=0,1,2,3$, $c=\pm, 2,3$ and $k$ (resp. $l$) is odd when
$a=2,3$ (resp. $b=2,3$).
\end{corollary}
\proof
Put $\pi_c(B_{b,l})=B_{b,l}+\sum_{j=1}^3A_{cbj}H_j$, $\pi_c(B_{a,k})=B_{a,k}+\sum_{j=1}^3A_{caj}H_j$ and $\pa_{a,k}B_{b,l}-\pa_{b,l}B_{a,k}+[B_{b,l},B_{a,k}]=\sum_{j=1}^3A_{cj}H_j$ with $A_{cbj}, A_{caj}, A_{cj}\in\mathcal{A}_{(c)}$. Then by using part a) of Theorem \ref{thm:Lax},
\begin{align*}
&\pa_{a,k}\pi_c(B_{b,l})-\pa_{b,l}\pi_c(B_{a,k})+[\pi_c(B_{b,l}),\pi_c(B_{a,k})]\nonumber\\
\in&\pa_{a,k}B_{b,l}+\sum_{j=1}^3A_{cbj}\cdot \pa_{a,k}H_j-\pa_{b,l}B_{a,k}-\sum_{j=1}^3A_{caj}\cdot\pa_{b,l}H_j+[B_{b,l}+\sum_{j=1}^3A_{cbj}H_j,B_{a,k}+\sum_{j=1}^3A_{caj}H_j]+\mathcal{A}_{(c)}H\nonumber\\
=&\pa_{a,k}B_{b,l}-\pa_{b,l}B_{a,k}+[B_{b,l},B_{a,k}]-\sum_{j=1}^3A_{cbj}\cdot H_j\cdot B_{a,k}+\sum_{j=1}^3A_{caj}\cdot H_j \cdot B_{b,l}\nonumber\\
&+\sum_{j=1}^3A_{cbj}\cdot H_j\cdot B_{a,k}-\sum_{j=1}^3A_{caj}\cdot H_j\cdot B_{b,l}+\mathcal{A}_{(c)}H=\mathcal{A}_{(c)}H.
\end{align*}
Therefore $\pa_{a,k}\pi_c(B_{b,l})-\pa_{b,l}\pi_c(B_{a,k})+[\pi_c(B_{b,l}),\pi_c(B_{a,k})]\in \mathcal{A}_{(c)}^0\cap \mathcal{A}_{(c)}H=\{0\}$.
\qed

Now we can prove part b) of Theorem \ref{thm:Lax}, that is,
$[\pa_{a,k},\pa_{b,l}]=0$ in $\R$. We have to prove that
$[\pa_{a,k},\pa_{b,l}]\mathcal{L}=0$ and
$[\pa_{a,k},\pa_{b,l}]H_c=0$. Note that 
\begin{align*}
&\pa_{a,k}\mathcal{L}=[B_{a,k},\mathcal{L}]+\sum_{j=1}^3A_{akj}H_j,\quad \pa_{b,l}\mathcal{L}=[B_{b,l},\mathcal{L}]+\sum_{j=1}^3A_{blj}H_j\\
&\pa_{a,k}H_c=-H_c B_{a,k}+\sum_{j=1}^3A_{akj}'H_j,\quad 
\pa_{b,l}H_c=-H_c B_{b,l}+\sum_{j=1}^3A_{blj}'H_j,
\end{align*}
for some $A_{akj},A_{blj},A_{akj}',A_{blj}'\in \mathcal{A}$. The
second derivatives take the form
\begin{align*}
&\pa_{a,k}(\pa_{b,l}\mathcal{L})\in [\pa_{a,k}B_{b,l},\mathcal{L}]+[B_{b,l},[B_{a,k},\mathcal{L}]]-\sum_{j=1}^3\Big(A_{akj}H_{j}B_{b,l}+A_{blj}H_{j}B_{a,k}\Big)+\mathcal{A}H,\\
&\pa_{a,k}(\pa_{b,l}H_c)\in H_c(B_{a,k}B_{b,l}-\pa_{a,k}B_{b,l})-\sum_{j=1}^3(A_{akj}'H_j B_{b,l}+A_{blj}'H_j B_{a,k})+\mathcal{A}H.
\end{align*}
Therefore, according to  Proposition \ref{zsequation} we have
\begin{align*}
&[\pa_{a,k},\pa_{b,l}]\mathcal{L}\in \left[\pa_{a,k}B_{b,l}-\pa_{b,l}B_{a,k}+[B_{b,l},B_{a,k}],\mathcal{L}\right]+\mathcal{A}H=\mathcal{A}H,\\
&[\pa_{a,k},\pa_{b,l}]H_c\in- H_c\cdot\Big(\pa_{a,k}B_{b,l}-\pa_{b,l}B_{a,k}+[B_{b,l},B_{a,k}]\Big)+\mathcal{A}H=\mathcal{A}H.
\end{align*}
Furthermore, we have
$[\pa_{a,k},\pa_{b,l}]\mathcal{L}\in\mathcal{A}_{(+)}^0\cap\mathcal{A}H=\{0\}$
and
$[\pa_{a,k},\pa_{b,l}]H_c\in\mathcal{A}_{(+)}^0\cap\mathcal{A}H=\{0\}$,
because $[\pa_{a,k},\pa_{b,l}]\mathcal{L}\in\mathcal{A}_{(+)}^0$ and
$[\pa_{a,k},\pa_{b,l}]H_c\in\mathcal{A}_{(+)}^0$. 

We claim that $[\pa_{a,k},\pa_{b,l}]L_c=0$. Let us prove only the case
$[\pa_{a,k},\pa_{b,l}]L_1=0$. The remaining two cases are essentially
the same. In fact,
$[\pa_{a,k},\pa_{b,l}]L_1=[\pa_{a,k}\pi_+(B_{b,l})-\pa_{b,l}\pi_+(B_{a,k})+[\pi_+(B_{b,l}),\pi_+(B_{a,k})],L_1]$,
which is zero according to Corollary \ref{projzsequation}.

 \section{Hirota bilinear equations and Lax operators}\label{sec:HBE-lax-1}
 
Let us recall the settings from Section \ref{sec:HBEs}. Our goal is to
construct Lax operators $\L$ and $H_a$ ($1\leq a\leq 3$) with
coefficients in $\O_\epsilon[\![\mathbf{t}]\!]$ and establish some of
their properties, which will be needed for the proof of Theorem
\ref{thm:HBE-sol}.

\subsection{Wave operators}
The first step is to transform the Hirota Bilinear Equations into an
equivalent system of quadratic equations involving operator series. 
Let us introduce the wave operators
\ben
&&
W_1^+(x,\mathbf{t},\Lambda)=S_1^+(x,\mathbf{t},\Lambda)e^{\xi_1(\mathbf{t},\Lambda)},
\quad S_1^+(x,\mathbf{t},\Lambda)=\sum_{j=0}^\infty
\psi_{1,j}^+(x,\mathbf{t})\Lambda^{-j}\\
&&
W_1^-(x,\mathbf{t},\Lambda)=e^{-\xi_1(\mathbf{t},\Lambda)}S_1^-(x,\mathbf{t},\Lambda),\quad 
S_1^-(x,\mathbf{t},\Lambda)=\sum_{j=0}^\infty\Lambda^{-j} \psi_{1,j}^-(x,\mathbf{t})\\
&& 
\xi_1(\mathbf{t},\Lambda)=
\sum_{k=1}^\infty\left(
t_{1,k}\, \Lambda^k+
t_{0,k}\, (\epsilon\pa_x-h_k)\frac{\Lambda^{(n-2)k}}{(n-2)^k k!} \right)
\een
and
\ben
&&
W_a(x,\mathbf{t},\partial_a)=
S_a(x,\mathbf{t},\pa_a)e^{\overline{\xi}_a(\mathbf{t},\pa_a)},\quad 
S_a(x,\mathbf{t},\pa_a)= \sum_{j=0}^\infty
\psi_{a,j}^+(x,\mathbf{t})\pa_a^{-j} \\
&& 
\overline{\xi}_a(\mathbf{t},\pa_a)=
\sum_{k=1}^\infty\left(
t_{a,2k+1}\pa_a^{2k+1}+t_{0,k} \, \epsilon\pa_x\, 
\frac{\pa_a^{2k}}{2^k k!} \right).
\een
Note that in the definition of $\overline{\xi}_a$ compared to the
definition of $\xi_a$ ($a=2,3$) the terms involving $y_a=t_{a,1}$ are
missing. The operator series $W^{\pm }_1$ and $W_a$ ($a=2,3$) take values
respectively in $\D_\epsilon(\!(\Lambda^{-1})\!) [\![\mathbf{t}]\!]$ and
$\D_\epsilon(\!(\partial_a^{-1})\!)[\![\mathbf{t}]\!]$. 

Given a pseudo-differential operator
$P(x,\mathbf{t},\partial_a)=\sum_j p_j(x,\mathbf{t}) \partial_a^{-j}
\in  \D_\epsilon(\!(\partial_a^{-1})\!)[\![\mathbf{t}]\!]$ then we
define its action on $(y_a-y_a')^0$ by the following rule
\ben
\partial_a^{-j}(y_a-y_a')^0=
\begin{cases}
\frac{1}{j!} \, (y_a-y_a')^j & ,\mbox{ if } j\geq 0,\\
0  & ,\mbox{ if } j<0,
\end{cases}
\een
In other words,  $P(x,\mathbf{t},\pa_a)\, (y_a-y_a')^0 = \sum_{j\geq 0} p_j(x,\mathbf{t})
\frac{(y_a-y_a')^j}{j!}$.  
\begin{proposition}\label{HQE-W}
The Hirota Bilinear Equations are equivalent to the following system
of quadratic equations for the operator series $W^{\pm}_1$ and $W_a$
($a=2,3$):
\ben
&&
W_1^+(x,\mathbf{t},\Lambda)
\frac{\Lambda^{(n-2)k-1}}{(n-2)^k k!} W_1^-(x,\mathbf{t}',\Lambda)+
\left(W_1^+(x,\mathbf{t}',\Lambda)
\frac{\Lambda^{(n-2)k-1}}{(n-2)^k k!} 
W_1^-(x,\mathbf{t},\Lambda)\right)^\#\\
&=&\sum_{m\in \mathbb{Z}}
\frac{1}{2}\left(
\Big(
W_2(x,\mathbf{t},\pa_2)\frac{\pa_2^{2k-1}}{2^k k!} 
W_2(x+m\epsilon,\mathbf{t}',\pa_2)^\#\pa_2
\left.\Big)\right|_{y'_2=y_2}(y_2-y_2')^0\right.\\ 
&&
\left.-(-1)^{m}\Big(
W_3(x,\mathbf{t},\pa_3)\frac{\pa_3^{2k-1}}{2^k k!} 
W_3(x+m\epsilon,\mathbf{t}',\pa_3)^\#\pa_3
\left.\Big)\right|_{y'_3=y_3}(y_3-y_3')^0\right)\Lambda^m,
\een
where $k\geq 0$ is an arbitrary non-negative integer.
\end{proposition}
The proof of Proposition \ref{HQE-W} is based on the following two
lemmas. 
\begin{lemma}\label{ABlambda}
  Let $A(x,\Lambda), B(x,\Lambda)$ be two operator series in
  $\D_\epsilon(\!(\Lambda^{-1})\!)$. Then 
  \ben
  A(x,\Lambda)\cdot B(x,\Lambda)^\#=
  \sum_{j\in \mathbb{Z}} \operatorname{Res}_{z=0} \frac{dz}{z}\left(
    A(x,\Lambda)(z^{\pm x/\epsilon})\,  \Big(B(x+j,\Lambda)(z^{\mp(x/\epsilon+j)})\Big)^\#\right)\Lambda^j.
  \een
\end{lemma}
\proof
Similar formula can be found in \cite{AM}. Put $A=\sum_k A_k(x)\Lambda^k$
and $B=\sum_l B_l(x)\Lambda^l$, then the RHS of the identity takes the
form
\ben
\sum_{k,l,j} \operatorname{Res}_{z=0} \frac{dz}{z}\Big(
A_k(x) z^{\pm(x/\epsilon+k)} \Big(B_l(x+j\epsilon) z^{\mp(x/\epsilon+l+j)}\Big)^\# \Lambda^j\Big) =
\sum_{l,j} A_{l+j}(x) \Lambda^j B_l(x) =
A(x,\Lambda)\cdot B(x,\Lambda)^\# ,
\een
where we used that $\Lambda^j f(x) = f(x+j\epsilon) \Lambda^j$.\qed

\begin{lemma}\label{ABpa}
	Let $A(y,\partial)=\sum_{i} a_i(y) \pa^i$ and
        $B(y,\partial)=\sum_j b_j(y)\pa^j$, where
        $\partial:=\frac{\pa}{\pa y}$,  be two 
        pseudo-differential operators. Then we have 
	\ben
	\Big(A(y,\partial)(B(y,\partial))^\#\pa\Big)\,
        (y-y')^{0}=\operatorname{Res}_{z=0}dz
        \Big(A(y,\partial)(e^{yz})B(y',\partial')(e^{-y'z})\Big), 
	\een
	where $\pa'=\frac{\pa}{\pa y'}$ and for $P(y,\partial)=\sum_k
        p_k(y)\pa^{-k}$ we define $P(y,\pa) (y-y')^0 = \sum_{k\geq 0}
        \frac{1}{k!}\, p_k(y) \, (y-y')^k$. 
\end{lemma}
\proof
        Let us recall the following formula 
	\ben
        \operatorname{Res}_{z=0}dz\Big(
        A(e^{yz})\, B(e^{-yz}) \Big) =
        \operatorname{Res}_{\pa}AB^{\#},
        \een
        where $\operatorname{Res}_\pa P(y,\pa)$ denotes the
        coefficient in front of $\pa^{-1}$ in $P(y,\partial).$ The
        proof is straightforward (see also \cite{Di}).
	Using the Taylor's expansion formula at $y'=y$ we get
	\ben
	&&\operatorname{Res}_{z=0}
        dz\Big(
        A(y,\partial)(e^{yz})\,
        B(y',\partial'))(e^{-y'z})\Big) \\ 
	&=&\operatorname{Res}_{z=0}
        dz\Big(
        A(y,\partial)(e^{yz})\sum_{n=0}^{\infty}
        \frac{(y'-y)^{n}}{n!}\pa_y^{n}B(y,\partial)(e^{-yz})
        \Big)   \\ 
	&=&
        \sum_{n=0}^{\infty}\frac{(y'-y)^{n}}{n!}
        \operatorname{Res}_{\pa}A(y,\partial)B^{\#}(y,\partial)
        (-1)^{n}\pa_y^{n}\\ 
	&=&
        \sum_{n=0}^{\infty}\frac{(y-y')^{n}}{n!}
        \operatorname{Res}_{\pa}A(y,\partial)B^{\#}(y,\partial)\pa^{n}
        \\
	&=&
        \sum_{n=0}^{\infty}\pa^{-n} \, (y-y')^{0}\,
        \operatorname{Res}_{\pa}A(y,\partial)B^{\#}(y,\partial)\pa^{n}.
	\een
	Put $A(y,\pa)B^{\#}(y,\pa)=\sum_{i\in \mathbb{Z}}c_{i}\pa^{i}$
        and notice that 
	\ben
	&&
        \operatorname{Res}_{\pa}A(y,\pa)B^{\#}(y,\pa)\pa^{n}=
        \operatorname{Res}_{\pa}
        \sum_{i\in \mathbb{Z}}c_{i}\pa^{i+n}=c_{-n-1}.
	\een
        We get
        \ben
        &&
        \operatorname{Res}_{z=0}\frac{dz}{z}\Big(
        A(y,\pa)(e^{yz})B(y',\pa')(e^{-y'z})\Big) \\
	&=&
        \sum_{n=0}^{\infty}c_{-n-1}\pa^{-n}\, (y-y')^{0}=
	\sum_{k\in \mathbb{Z}}c_{k}\pa^{k+1}\, (y-y')^{0}\\
	&=&
        \Big(A(y,\pa)B^{\#}(y,\pa)\pa\Big)\, (y-y')^{0}.
        \qed
	\een
        
{\em Proof of Proposition \ref{HQE-W}.}
Note that
\ben
\Psi^+_1(x,\mathbf{t},z) = W_1^+(x,\mathbf{t},\Lambda)
(z^{x/\epsilon-\frac{1}{2}})\quad \mbox{and}\quad
\Psi^-_1(x,\mathbf{t},z) ^\#= W_1^-(x,\mathbf{t},\Lambda)^\#
(z^{-x/\epsilon-\frac{1}{2}}).
\een
Using Lemma \ref{ABlambda} we get
\ben
\sum_{m\in \ZZ} {\rm Res}_{z=0}
\frac{z^{(n-2)k}}{(n-2)^kk!}\frac{dz}{z}\Big(
(\Psi_1^+(x,\mathbf{t},z)\Psi_1^-(x+m\epsilon,\mathbf{t}',z)
\Big)\Lambda^m=
W^+_1(x,\mathbf{t},\Lambda) \,
\frac{\Lambda^{(n-2)k-1}}{(n-2)^kk!}\,
W^-_1(x,\mathbf{t}',\Lambda)
\een
and 
\ben
\sum_{m\in \ZZ} {\rm Res}_{z=0}
\frac{z^{(n-2)k}}{(n-2)^kk!}\frac{dz}{z}\Big(
\Psi_1^+(x+m\epsilon,\mathbf{t}',z)\Psi_1^-(x,\mathbf{t},z)
\Big)^\# \Lambda^m =
\left(
W^+_1(x,\mathbf{t}',\Lambda) \,
\frac{\Lambda^{(n-2)k-1}}{(n-2)^kk!}\,
W^-_1(x,\mathbf{t},\Lambda)\right)^\#.
\een
Similarly, for $a=2,3$,  using that
\ben
\Psi^+_a(x,\mathbf{t},z) = W_a(x,\mathbf{t},\pa_a) e^{y_a z}
\quad\mbox{and}\quad
\Psi^-_a(x,\mathbf{t},z) = W_a(x,\mathbf{t},\pa_a) e^{-y_a z}
\een
and Lemma \ref{ABpa} we get
\ben
&&
{\rm Res}_{z=0}\frac{z^{2k}}{2^k
  k!}\frac{dz}{2z}\left(
  \Psi_a^+(x,\mathbf{t},z)\Psi_a^-(x+m\epsilon,\mathbf{t}',z)\right)=\\
&&
\Big(
W_a(x,\mathbf{t},\pa_a) \, \frac{\pa_a^{2k-1}}{2^k k!}\,
W_a(x+m\epsilon,\mathbf{t}',\pa_a)^\#\,\pa_a
\left.\Big)\right|_{y_a=y_a'} (y_a-y_a')^0.
\qed
\een

\begin{proposition}\label{prop:dres}
  a) The operators $S^{\pm}_1$ satisfy the following relation
  \ben
  S^+_1(x,\mathbf{t},\Lambda) \Lambda^{-1}
  S^-_1(x,\mathbf{t},\Lambda)
  = \sum_{m=0}^\infty \Lambda^{-2m-1}.
  \een
 b) The operators $S_a(x,\mathbf{t},\pa_a)$ ($a=2,3$) satisfy the
 following relations
 \ben
 S_a(x,\mathbf{t},\pa_a) \pa_a^{-1} S_a(x,\mathbf{t},\pa_a)^\# =
 \pa_a^{-1}.
 \een
\end{proposition}
\proof
a) Let us substitute in the identity of Proposition \ref{HQE-W} $k=0$,
$\mathbf{t}'=\mathbf{t}$ and compare the coefficients in front of the
negative powers of $\Lambda$. The identity that we want to prove
follows.

b) Let us make the same substitution as in a) except for
$t_{a,1}'=y'_a$, i.e., we keep $y_a'\neq y_a$. Comparing the
coefficients in front of $\Lambda^0$ yields the identity that we want
to prove.
\qed

For brevity put  $B=\O_\epsilon[\![\mathbf{t}]\!]$ and let us denote
by $\E=B[\Lambda,\Lambda^{-1},\pa_2,\pa_3]$ the ring of differential
difference operators. Note that the ring $\E$ acts on
$\Psi:=(\Psi^+_1,(\Psi^-_1)^\#,\Psi_2,\Psi_3)$ by acting on each
component.
\begin{proposition}\label{prop:anul}
  The annihilator $\operatorname{Ann}(\Psi)=\{P\in \E\ |\ P(\Psi)=0\}$
  is a left ideal in $\E$ generated by
  \ben
  H_1=\pa_2\cdot\pa_3+q_1,\quad
  H_2=(\Lambda-1)\cdot\pa_2-q_2\cdot (\Lambda+1),\quad
  H_3=(\Lambda+1)\cdot\pa_3-q_3\cdot (\Lambda-1),
  \een
where $q_1=2\pa_2\pa_3(\log \tau)$,
$q_a=\pa_a (\log f)$ ($a=2,3$), where 
$f:=\frac{\tau(x,\mathbf{t})}{\tau(x+\epsilon,\mathbf{t})}.$ 
\end{proposition}
\proof
It is enough to prove that $H_i\in \operatorname{Ann}(\Psi)$. Indeed,
suppose that this is proved. Let $P\in \operatorname{Ann}(\Psi)$ be
arbitrary.  The operator $P$ can be
decomposed as a sum of an operator in the left ideal $\E
(H_1,H_2,H_3)$ generated by $H_i$ 
and
$Q=P_1+P_2+P_3$, where
$P_1\in B[\Lambda^{\pm 1}]$, $P_2\in B[\pa_2]$, and $P_3\in
B[\pa_3]$. If $P_2\neq 0$ then let $a(x,\mathbf{t})\pa_2^m$ be its
highest order term. We get that $Q(\Psi_2) = a(x,\mathbf{t})
z^m(1+O(z^{-1}))e^{\xi_2(\mathbf{t},z)}$ can not be zero --
contradiction. Similar 
argument shows that $P_3=0$, and finally $0=P(\Psi_1^+) = (P_1\cdot
W_1^{+}) (z^{x/\epsilon-1/2})$ implies $P_1\cdot W_1^+=0$ $\Rightarrow$
$P_1=0$.

Let us prove that $H_1(\Psi)=0$. The argument in the remaining two
cases is similar. Let us apply the differential operator
$\partial_{t_{2,1}'}\partial_{t_{3,1}'}$ to the identiy in Proposition
\ref{HQE-W} with $k=0$, substitute $\mathbf{t}'=\mathbf{t}$, and
compare the coefficients in front of the negative powers of $\Lambda$. We get
\ben
S^+_1(x,\mathbf{t},\Lambda) \Lambda^{-1} \pa_2\pa_3
S^-_1(x,\mathbf{t},\Lambda) = \frac{1}{2}\sum_{m<0} \Lambda^m
\Big( \pa_3 \psi^+_{2,1}(x,\mathbf{t}) -(-1)^m \pa_2
\psi^+_{3,1}(x,\mathbf{t})
\Big).
\een
Note that
$
\pa_3\psi^+_{2,1}(x,\mathbf{t}) =
\pa_2\psi^+_{3,1}(x,\mathbf{t}) = -q_1.
$
Recalling Proposition \ref{prop:dres}, a) we get  
$
\pa_2\pa_3 S^-_1 = - S^-_1\, q_1
$
and
$
\pa_2\pa_3 S^+_1 = - q_1 S^+_1.
$
These two relations are equivalent to $H_1(\Psi^+_1)=0$ and
$H_1((\Psi_1^{-} ) ^\#)=0$.  In order to
prove that $H_1(\Psi_2)=0$, we proceed as above except that we leave
$t_{2,1}'$ arbitrary, set the remaining components of $\mathbf{t}'$
and $\mathbf{t}$ to be equal, and compare the coefficients in front of
$\Lambda^0$. We get $(S_2\cdot \pa_2^{-1} \cdot (\pa_3(S_2))^\#\cdot  \pa_2^2
)(y_2-y_2')^0= q_1(y_2')$,
where we suppressed the dependance of $q_1(x,\mathbf{t})$ on
$x$ and on the remaining components of $\mathbf{t}$. 
Taking the Taylor's expansion of $-q_1(y_2')$ at $y_2'=y_2$
and comparing the coefficients in front of the powers of $y_2'-y_2$ we
get 
\ben
S_2\cdot \pa_2^{-1} \cdot (\pa_3(S_2))^\#\cdot  \pa_2^2 =
\sum_{j=0}^\infty (-1)^j (\pa_2^jq_1(x,\mathbf{t})) \pa_2^{-j} =
\pa_2^{-1}\cdot q_1(x,\mathbf{t}) \cdot \pa_2.
\een
Recalling Proposition \ref{prop:dres}, b) we get $(\pa_3S_2)^\#\pa_2 =
S_2^\#\, q_1$ $\Rightarrow$ $\pa_2\cdot(\pa_3 S_2(x,\mathbf{t},\pa_2)) =
-q_1 S_2(x,\mathbf{t},\pa_2)$. This identity is equivalent to
$H_1(\Psi_2)=0$. Finally, the argument for $H_1(\Psi_3)=0$ is
completely analogus.
\qed

To avoid cumbersome notation let us put $S_1:=S^+_1$.
Proposition \ref{prop:anul} is equivalent to a set of relations involving the
operators $H_i$ ($1\leq i\leq 3$) and the wave operators $S_j$ ($1\leq
j\leq 3$). These relations take the following form.
\begin{corollary}\label{cor:anul}
a) The operator $H_1$ and the wave operators satisfy the following
relations:
\beqa
&&
\label{h1s1}
H_1\cdot S_1=\pa_2\cdot S_1\cdot \pa_3+\pa_3\cdot S_1\cdot\pa_2-S_1\cdot\pa_2\cdot\pa_3,\\
&&H_1\cdot S_2=\pa_2\cdot S_2\cdot\pa_3,\nonumber \\
&&
H_1\cdot S_3=\pa_3\cdot S_3\cdot\pa_2.\nonumber 
\eeqa
b) The operator $H_2$ and the wave operators satisfy the following
relations:
\beqa
&&H_2\cdot S_1=(\Lambda-1)\cdot S_1\cdot\pa_2,\notag\\
&&H_2\cdot S_2=(\pa_2+q_2) \cdot S_2\cdot(\Lambda-1),\notag\\
&&
\label{h2s3}
H_2\cdot S_3=(\Lambda\cdot S_3+S_3\cdot\Lambda)\cdot \pa_2-(\pa_2+q_2)\cdot S_3\cdot(\Lambda+1).
\eeqa
c) The operator $H_3$ and the wave operators satisfy the following
relations:
\beqa
&&H_3\cdot S_1=(\Lambda+1)\cdot S_1\cdot\pa_3,\nonumber\\
&&
\label{h3s2}
H_3\cdot S_2=(\Lambda\cdot S_2+S_2\cdot\Lambda)\cdot \pa_3-(\pa_3+q_3)\cdot S_2\cdot(\Lambda-1),\\
&&H_3\cdot S_3=(\pa_3+q_3)\cdot S_3\cdot(\Lambda+1).\nonumber\qed
\eeqa
\end{corollary}
A straightforward computation yields that the complicated looking
relations \eqref{h1s1}, \eqref{h2s3}, and \eqref{h3s2} can be replaced
equvalently with the following simple relation
\beq\label{efg-relation}
q_1+q_1[1]+2q_2q_3 = 0,
\eeq
where the functions $q_1,q_2$, and $q_3$ are defined in Proposition
\ref{prop:anul} and for a pseudo-differential-difference operator $P$
we put $P[m] := \Lambda^m(P)=\Lambda^m\cdot P\cdot \Lambda^{-m}$ for
its translation by $m$, i.e., substituting $x\mapsto x+m\epsilon$.

\subsection{The Lax operator}\label{sec:lop}
Following Shiota's construction in \cite{Shi} we 
introduce the following rings of pseudo-differential-difference operators:  
\ben
\mathcal{E}_{(\pm)}=B[\pa_2,\pa_3]((\Lambda^{\mp 1})),\quad
\mathcal{E}_{(2)}=B[\Lambda,\Lambda^{-1},\pa_3]((\pa_2^{-1})),\quad
\mathcal{E}_{(3)}=B[\Lambda,\Lambda^{-1},\pa_2]((\pa_3^{-1})),
\een
\ben
\mathcal{E}_{(\pm)}^{0}=B((\Lambda^{\mp 1})),\quad
\mathcal{E}_{(2)}^{0}=B((\pa_2^{-1})),\quad
\mathcal{E}_{(3)}^{0}=B((\pa_3^{-1})),
\een
and the quotient rings
$\mathcal{A}''=\mathcal{A'}/\mathcal{A}H,$
for $\mathcal{A}=\mathcal{E}$, $\mathcal{E}_{(\pm)}$,
$\mathcal{E}_{(2)}$ or $\mathcal{E}_{(3)}$, where
$\mathcal{A}H=\mathcal{A}H_1+\mathcal{A}H_2+\mathcal{A}H_3$ is the
left ideal of $\A$ generated by $H_i$ ($1\leq i\leq 3$) and 
\ben
\mathcal{A}'=\{P\in\mathcal{A}\ |\ H_iP\in \A H\ \forall \ i=1,2,3\}.
\een
Let us introduce the following auxiliarly Lax operators
\ben
&&
L_1^+(x,\mathbf{t},\Lambda):=
S_1^+(x,\mathbf{t},\Lambda)\cdot\Lambda\cdot S_1^+(x,\mathbf{t},\Lambda)^{-1}=:
u_{1,0}^+(x,\mathbf{t})\Lambda+\sum_{j=1}^\infty u_{1,j}^+(m,\mathbf{t})\Lambda^{1-j},\\
&&
L_1^-(x,\mathbf{t},\Lambda):=
S_1^-(x,\mathbf{t},\Lambda)^\#\cdot\Lambda^{-1}\cdot
\left(S_1^-(x,\mathbf{t},\Lambda)^\#\right)^{-1}=:
u_{1,0}^-(x,\mathbf{t})\Lambda^{-1}+\sum_{j=1}^\infty u^-_{1,j}(x,\mathbf{t})\Lambda^{j-1},\\
&&
L_a(x,\mathbf{t},\pa_a):=
S_a(x,\mathbf{t},\pa_a)\cdot\pa_a^{-1}\cdot
S_a(x,\mathbf{t},\pa_a)^{-1}=:
\pa_a+\sum_{j=1}^\infty u_{a,j}(x,\mathbf{t})\pa_a^{-j}\quad (a=2,3).
\een
To avoid cumbersome notation we put $L_1(x,\mathbf{t},\Lambda) =
L_1^+(x,\mathbf{t},\Lambda)$ and $u_{1,j} = u^+_{1,j}$. Recalling
Proposition \ref{prop:dres} we get the following relations
\ben
&&
L_1^+(x,\mathbf{t},\Lambda)^\#=
(\Lambda^{-1}-\Lambda)\cdot L_1^-(x,\mathbf{t},\Lambda)\cdot\
\iota_{\Lambda} \left(\Lambda^{-1}-\Lambda\right)^{-1},\\ 
&&
L_a(x,\mathbf{t},\pa_a)^\#=-\pa_a\cdot
L_a(x,\mathbf{t},\pa_a)\cdot\pa_a^{-1}\quad (a=2,3),
\een
where $\iota_\Lambda $ denotes the operation that takes the Laurent
series expansion at $\Lambda=0$. Note that the coefficient of $L_a$
in front of $\pa_a^0$ must be $u_{a,0}=0$.  
\begin{lemma}\label{le:c_a}
  The following formula holds
  \ben
  c_a(x,\mathbf{t}) :=\Big(L_a(x,\mathbf{t},\pa_a)^2/2\Big)_{a,[0]} =
  2\pa_a^2(\log \tau(x,\mathbf{t}) )\quad (a=2,3). 
  \een
\end{lemma}
\proof
By definition $S_a = 1 + \psi^+_{a,1} \pa_a^{-1}+\psi^+_{a,2}\pa_a^{-2}+\cdots$ with
$\psi^+_{a,1} = -2\pa_a(\log \tau)$. Put $L_a^2= \pa_a^2 + 2 c_a
+\cdots$. By definition $L_a S_a = S_a \pa_a^2$. Comparing the
coefficients in front of $\pa_a^0$ we get
\ben
\psi^+_{a,2} + 2 \pa_a(\psi^+_{a,1}) + 2 c_a = \psi_{a,2}^+.
\een
Solving for $c_a$ we get
$c_a = -\pa_a (\psi^+_{a,1}) = 2\pa_a^2 (\log \tau)$.
\qed

Let us define the following differential-difference operator
\ben
\mathcal{L} & = & \frac{\pa_2^2}{2}+\frac{\pa_3^2}{2}+\left(
  \left( \tfrac{L_1^{n-2}}{n-2} \, \sum_{m=0}^\infty
  \Lambda^{-2m-1} \right)_{1,>0} -
  \left(\sum_{m=0}^\infty
  \Lambda^{2m+1}\,  \tfrac{(L_1^\#)^{n-2} }{n-2}\right)_{1,<0} \right)\,
(\Lambda-\Lambda^{-1})+ \\
&&
+\frac{1}{4}(c_2-c_3)(\Lambda+\Lambda^{-1}) +\frac{1}{2}(c_2+c_3),
\een
where $c_a$ $(a=2,3)$ are the same as in Lemma \ref{le:c_a}. Note that
in both sums only finitely many terms contribute and  that the term
involving $(\ )_{1,<0}$ is conjugated via the anti-involution $\#$ to
the term involving $(\ )_{1,>0}$. 
\begin{proposition}\label{prop:lax-rel}
  Let us introduce the following two rational difference operators
  \ben
  Q_2(x,\mathbf{t},\Lambda) :=
  (\Lambda-1)^{-1} \cdot q_2(x,\mathbf{t}) \cdot (\Lambda+1)
  \quad \mbox{and}\quad
  Q_3(x,\mathbf{t},\Lambda) :=
  (\Lambda+1)^{-1} \cdot q_3(x,\mathbf{t})\cdot  (\Lambda-1)
  \een
  Then the following formulas hold
\ben
\mathcal{L} = \tfrac{(L^\pm_1)^{n-2}}{n-2} + \tfrac{1}{2}
\iota_{\Lambda^{\mp 1}} \Big(
(\pa_2+Q_2)(\pa_2-Q_2) +
(\pa_3+Q_3)(\pa_3-Q_3) \Big). 
\een
\end{proposition}
\proof
Let us recall the Hirota quadratic equations from Proposition
\ref{HQE-W} with $k=1$ and $\mathbf{t}'=\mathbf{t}$. We get
\ben
&&
S_1 \Big(\tfrac{\Lambda^{n-2}}{n-2}\Big) S_1^{-1} \sum_{m=0}^\infty
\Lambda^{-2m-1} + \Big(   S_1 \Big(\tfrac{\Lambda^{n-2}}{n-2}\Big) S_1^{-1} \sum_{m=0}^\infty
\Lambda^{-2m-1}\Big)^\# =\\
&&
\frac{1}{2} \sum_{m\in \ZZ}
\left(
  \Big(\tfrac{L_2^2}{2} \, S_2 \Lambda^m S_2^{-1}\Big)_{2,[0]} -
  \Big(\tfrac{L_3^2}{2} \, S_3 (-\Lambda)^m S_3^{-1}\Big)_{3,[0]}
\right)
\een
Using that
\ben
\sum_{m\in \ZZ}\Lambda^m=
\iota_{\Lambda^{-1}} (\Lambda-1)^{-1} -
\iota_{\Lambda} (\Lambda-1)^{-1}
\quad\mbox{and}\quad
\sum_{m\in \ZZ}(-\Lambda)^m=
-\iota_{\Lambda^{-1}} (\Lambda+1)^{-1} +
\iota_{\Lambda} (\Lambda+1)^{-1}
\een
we write the above identity as
\beqa\label{hqe_{k=1}_{t'=t}}
&&
S_1 \Big(\tfrac{\Lambda^{n-2}}{n-2}\Big) S_1^{-1} \sum_{m=0}^\infty
\Lambda^{-2m-1} + \Big(   S_1 \Big(\tfrac{\Lambda^{n-2}}{n-2}\Big) S_1^{-1} \sum_{m=0}^\infty
\Lambda^{-2m-1}\Big)^\# =\\
\nonumber
&&
\frac{1}{2}(\iota_{\Lambda^{-1}} -\iota_{\Lambda})
\left(
  \Big(\tfrac{L_2^2}{2} \, S_2 (\Lambda-1)^{-1} S_2^{-1}\Big)_{2,[0]} +
  \Big(\tfrac{L_3^2}{2} \, S_3 (\Lambda+1)^{-1} S_3^{-1}\Big)_{3,[0]}
\right).
\eeqa
Recalling the second identity from Corollary \ref{cor:anul}, b) we get
\ben
S_2(\Lambda-1) S_2^{-1} = (\pa_2 + q_2)^{-1} H_2.
\een
We will need the coefficients $v_1$ and $v_2$ from the expansion
\ben
S_2(\Lambda-1)^{-1} S_2^{-1} =
(\Lambda-1)^{-1} + v_1 \pa_2^{-1} + v_2 \pa_2^{-2} + \cdots , 
\een
where the coefficients $v_i$ on the  RHS should be interpreted as
rational difference operators. Since $H_2= (\Lambda-1) \pa_2 -q_2
(\Lambda+1)$ the above formula implies that
\ben
((\Lambda-1) \pa_2 -q_2 (\Lambda+1)\cdot
((\Lambda-1)^{-1} + v_1 \pa_2^{-1} + v_2 \pa_2^{-2} + \cdots )=
\pa_2 + q_2.
\een
Comparing the coefficients in front of $\pa_2^{0}$ and $\pa_2^{-1}$ yields 
\beq\label{v12}
v_1 = 2 Q_2 \, (\Lambda-\Lambda^{-1})^{-1},\quad
(v_2+\pa_2(v_1)) = 2 Q_2^2\,  (\Lambda-\Lambda^{-1})^{-1}.
\eeq
Similar argument using the 3rd identity in Corollary \ref{cor:anul},
c) implies that the coefficients $w_1$ and $w_2$ in the expansion
\ben
S_3(\Lambda+1)^{-1} S_3^{-1} =
(\Lambda+1)^{-1} + w_1 \pa_3^{-1} + w_2\pa_3^{-2} +\cdots
\een
are gven by
\beq\label{w12}
w_1 = 2 Q_3 \, (\Lambda-\Lambda^{-1})^{-1},\quad
(w_2+\pa_2(w_1)) = 2 Q_3^2\,  (\Lambda-\Lambda^{-1})^{-1}.
\eeq
Formulas \eqref{v12} and \eqref{w12} imply that
\ben
\Big(\tfrac{L_2^2}{2} \, S_2 (\Lambda-1)^{-1} S_2^{-1}\Big)_{2,[0]} =
c_2 (\Lambda-1)^{-1}+
\Big(\pa_2(Q_2) + Q_2^2\Big)(\Lambda-\Lambda^{-1})^{-1}
\een
and
\ben
\Big(\tfrac{L_3^2}{2} \, S_3 (\Lambda+1)^{-1} S_3^{-1}\Big)_{3,[0]} =
c_3 (\Lambda+1)^{-1}+
\Big(\pa_3(Q_3) + Q_3^2\Big)(\Lambda-\Lambda^{-1})^{-1}.
\een
The identity \eqref{hqe_{k=1}_{t'=t}} takes the form
\beqa\label{lax-rel-1}
&&
\left(
  \left( \tfrac{L_1^{n-2}}{n-2} \, \sum_{m=0}^\infty
    \Lambda^{-2m-1} \right) +
  \left( \tfrac{L_1^{n-2}}{n-2} \, \sum_{m=0}^\infty
    \Lambda^{-2m-1} \right)^\#  
\right)\, (\Lambda-\Lambda^{-1}) =\\
&&
\frac{1}{2}(\iota_{\Lambda^{-1}} -\iota_{\Lambda})
\left(
  c_2 (1+\Lambda^{-1})+ c_3(1-\Lambda^{-1}) +
  \Big(\pa_2(Q_2) + Q_2^2+
  \pa_3(Q_3) + Q_3^2
  \Big)
\right).
\eeqa
Let us proof the formula in the proposition for $L_1^+$. The other case
is similar. By definition we have
\ben
\mathcal{L} -\tfrac{L_1^{n-2}}{n-2} & = &
\frac{\pa_2^2}{2}+\frac{\pa_3^2}{2}-\left(
  \left( \tfrac{L_1^{n-2}}{n-2} \, \sum_{m=0}^\infty
  \Lambda^{-2m-1} \right)_{1,<0} +
  \left(\sum_{m=0}^\infty
  \Lambda^{2m+1}\,  \tfrac{(L_1^\#)^{n-2} }{n-2}\right)_{1,<0} \right)\,
(\Lambda-\Lambda^{-1})+ \\
&&
+c_1(\Lambda-\Lambda^{-1})+
\frac{1}{4}(c_2-c_3)(\Lambda+\Lambda^{-1}) +
\frac{1}{2}(c_2+c_3),
\een
where
\ben
c_1:=- \left( \tfrac{L_1^{n-2}}{n-2} \, \sum_{m=0}^\infty
  \Lambda^{-2m-1} \right)_{1,[0]} = -\frac{1}{4} (c_2-c_3), 
\een
where the second equality in the above formula is proved by comparing
the coefficients in front of $\Lambda^0$ in \eqref{hqe_{k=1}_{t'=t}}. 
The identity \eqref{lax-rel-1} allows us to transform the above
formula into
\ben
\mathcal{L} -\tfrac{L_1^{n-2}}{n-2} & = &
\frac{\pa_2^2}{2}+\frac{\pa_3^2}{2} -
\frac{1}{2}\iota_{\Lambda^{-1}}
\Big(
\pa_2(Q_2) + Q_2^2+ \pa_3(Q_3) + Q_3^2
\Big)
-
\frac{1}{2} \Big(
c_2 (1+\Lambda^{-1})+ c_3(1-\Lambda^{-1}) 
\Big) +\\
&&
+c_1(\Lambda-\Lambda^{-1})+
\frac{1}{4}(c_2-c_3)(\Lambda+\Lambda^{-1}) +
\frac{1}{2}(c_2+c_3) .
\een
It remains only to check that the terms that involve $c_1,c_2$, and
$c_3$ add up to 0 and to use that
\ben
(\pa_a+Q_a)(\pa_a-Q_a) = \pa_a^2 -\pa_a(Q_a) - Q_a^2\quad
(a=2,3).
\qed
\een
\begin{proposition}\label{prop:ef-L}
  The wave functions are eigenfunctions of $\mathcal{L}$
  \ben
  \mathcal{L}(\Psi_1^+) = \tfrac{z^{n-2}}{n-2} \Psi_1^+ ,\quad
  \mathcal{L}((\Psi_1^-)^\#) = \tfrac{z^{n-2}}{n-2} (\Psi_1^-)^\# ,\quad
  \mathcal{L}(\Psi_a) = \tfrac{z^2}{2} \Psi_a\quad
  (a=2,3). 
  \een
\end{proposition}
\proof
The first identity is a direct consequence from Propositon
\ref{prop:lax-rel}. Indeed, using  that $\Psi_1^+(x,\mathbf{t},z) =
S_1(x,\mathbf{t},\Lambda) \left(e^{\xi_1(\mathbf{t},z)}
  z^{x/\epsilon-\frac{1}{2}} \right)$ and that the expansions at $\Lambda=0$
of $(\pa_2 -Q_2) = (\Lambda-1)^{-1} H_2$ and
$\pa_3-Q_3 = (\Lambda+1)^{-1} H_3$ anihilate $\Psi_1^+$
we get 
\ben
\mathcal{L}(\Psi_1^+) =
\frac{1}{n-2} (L_1^{n-2} \cdot S_1)
\left(e^{\xi_1(\mathbf{t},z)} z^{x/\epsilon-\frac{1}{2}} \right)=
\frac{1}{n-2} (S_1 \cdot \Lambda^{n-2} )
\left(e^{\xi_1(\mathbf{t},z)} z^{x/\epsilon-\frac{1}{2}} \right)=
\frac{z^{n-2}}{n-2}  \, \Psi_1^+.
\een
The formula for $\mathcal{L}((\Psi^-_1)^\#)$ is proved in a similar way.

Let us prove the formula for $\mathcal{L}(\Psi_2)$. The computation
for $\mathcal{L}(\Psi_3)$ is similar. Let us apply the operator
$\mathcal{L}(\Lambda M^{-1},\pa_2,\pa_3)$ to the set of Hirota bilinear
equations (HBEs) \eqref{hbe-psi}, where $M$ acts on the set of 
HBEs via the shift $m\mapsto m+1$. Note that
\ben
\mathcal{L}(\Lambda M^{-1},\pa_2,\pa_3)
\left(
  F(x,\mathbf{t},z ) \, F(x+m\epsilon,\mathbf{t}',z)
\right)  =
\left(
\mathcal{L}(\Lambda,\pa_2,\pa_3) (F(x,\mathbf{t},z ))
\right)\,
F(x+m\epsilon,\mathbf{t}',z),
\een
and
\ben
\mathcal{L}(\Lambda M^{-1},\pa_2,\pa_3) \left(
  (-1)^m
  F(x,\mathbf{t},z ) \, F(x+m\epsilon,\mathbf{t}',z)\right)  =
\left(
  \mathcal{L}(-\Lambda,\pa_2,\pa_3) (F(x,\mathbf{t},z) )
\right)
\,
F(x+m\epsilon,\mathbf{t}',z).
\een
In other words the action of $\mathcal{L}(\Lambda M^{-1},\pa_2,\pa_3)$
on the set of HBEs is equivalent to acting by
$\mathcal{L}(\Lambda,\pa_2,\pa_3)$ (or
$\mathcal{L}^*:=\mathcal{L}(-\Lambda,\pa_2,\pa_3)$ whenever the sign
factor $(-1)^m$ is present) only on the first factors of the 
bilinear identities.  The LHS of \eqref{hbe-psi} with $k=0$ is
transformed into the LHS of \eqref{hbe-psi} with $k=1$, because
$\Psi_1^+$ and $(\Psi_1^-)^\# $ are eigenfunctions of $\mathcal{L}$
with eigenvalue $\tfrac{z^{n-2}}{n-2}$. Therefore using also \eqref{hbe-psi}
with $k=1$ we get that the residue 
\ben
{\rm Res}_{z=0}
\frac{dz}{2z}
\left(
  \Big(\mathcal{L}-\tfrac{z^2}{2}\Big)
  \left(
  \Psi_2^+(x,\mathbf{t},z)
\right)
\Psi_2^-(x+m\epsilon,\mathbf{t}',z)-
(-1)^{m}
\Big(\mathcal{L}^*-\tfrac{z^2}{2}\Big)
\left(
  \Psi_3^+(x,\mathbf{t},z)
  \right)
\Psi_3^-(x+m\epsilon,\mathbf{t}',z)
\right)
\een
is $0$.  Let us substitute now $m=0$ and $\mathbf{t}'=\mathbf{t}$,
except for $t_{2,1}'$, i.e., we keep $y_2'=t_{2,1}'$ and $y_2=t_{2,1}$
arbitrary. 
Note that
\ben
\Big(\mathcal{L}-\tfrac{z^2}{2}\Big)
  \left(
  \Psi_2^+
\right) =
\Big(
\mathcal{L} (S_2^+) + \pa_2(S_2^+) \cdot
\pa_2
\Big) (e^{\xi_2(\mathbf{t},z)}) = \sum_{j=1}^\infty
\widetilde{\psi}^+_{2,j} (x,\mathbf{t}) z^{-j} \, e^{\xi_2(\mathbf{t},z)},
\een
where we used that
$
\mathcal{L} (1) + \pa_2 (\psi^+_{2,1}) =
c_2 + \pa_2 (\psi^+_{2,1})  = 0. 
$
Similarly,
\ben
\Big(\mathcal{L}^*-\tfrac{z^2}{2}\Big)
  \left(
  \Psi_3^+
\right) =
\Big(
\mathcal{L}^*(S_3^+) + \pa_3(S_3^+) \cdot
\pa_3
\Big) (e^{\xi_3(\mathbf{t},z)}) = \sum_{j=1}^\infty
\widetilde{\psi}^+_{3,j} (x,\mathbf{t}) z^{-j} \, e^{\xi_3(\mathbf{t},z)},
\een
where we used that $\mathcal{L}^*(1)+\pa_3(\psi^+_{3,1}) = c_3+
\pa_3(\psi^+_{3,1})=0$. 
Therefore the vanishing of the above residue yields the following
identity
\ben
\sum_{j=1}^\infty 
{\rm Res}_{z=0}
\frac{dz}{2z}
\left(
  \widetilde{\psi}^+_{2,j} \, z^{-j} \, e^{(y_2-y_2') z} \,
  \psi^-_2(x,\mathbf{t}',z)-
  \widetilde{\psi}^+_{3,j} \, z^{-j} \,
  \psi^-_3(x,\mathbf{t}',z)
  \right) = 0
  \een
The terms involving $\psi_3$ do not contribute. Expanding the
exponential in the powers of 
$y_2-y_2'$ we get a sequence of vanishing residues. Since
$\psi^-_2=1+O(z^{-1} )$, a simple induction
on $j$ shows that $\widetilde{\psi}^+_{2,j}=0$ for all $j\geq 1$. This
is exactly what we had to prove.
\qed

Let us point out that in the above proof we have established also the
identities 
$
\mathcal{L} (S_2^+)  = -\pa_2(S_2^+)\cdot \pa_2 
$
and
$
\mathcal{L}^* (S_3^+)  = -\pa_3(S_3^+)\cdot \pa_3.  
$
\begin{proposition}
The differential-difference operator $\mathcal{L}$ satisfies the
following relations

a) $\mathcal{L}-\tfrac{(L^{\pm}_1)^{n-2}}{n-2} \in \E_{(\pm)} H$.

b) $\mathcal{L}-\tfrac{L_2^2}{2} \in \E_{(2)} H$.

c) $\mathcal{L}-\tfrac{L_3^2}{2} \in \E_{(3)} H$.
\end{proposition}
\proof
Part a) is a direct consequence of Proposition \ref{prop:lax-rel}. For
part b). First, note that $\mathcal{L} (S_2)  = -\pa_2(S_2)\cdot
\pa_2 $ implies
\ben
\Big(\mathcal{L}-\frac{\pa_2^2}{2}\Big)(S_2)=
\frac{1}{2} (S_2 \,\pa_2^2-\pa_2^2\, S_2) = 
\frac{1}{2}(L_2^2)_{2,\leq 0}\, S_2,
\een
where we used that $(L_2^2)_{2,>0}=\pa_2^2.$ Using the above identity
we get
\beq\label{L-L_2^2/2}
\mathcal{L}-\frac{1}{2} L_2^2 =
\left(
  \Big(\mathcal{L}-\frac{1}{2}\pa_2^2\Big) \cdot S_2-
  \Big(\mathcal{L}-\frac{1}{2}\pa_2^2\Big) (S_2)
\right)\cdot S_2^{-1}. 
\eeq
The operator $\mathcal{L}-\frac{1}{2}\pa_2^2$ is a sum of terms of the
form $ a(x,\mathbf{t}) \Lambda^i$ and $\pa_3^2/2$. We claim that each
of these terms contributes to \eqref{L-L_2^2/2} a term that belongs to
$\E_{(2)}H$. Indeed, we have
\ben
\Big( a \Lambda^i \cdot S_2 - a\, \Lambda^i(S_2)\Big) \cdot S_2^{-1} =
a S_2[i] (\Lambda^i-1)  S_2^{-1}
\een
and
\ben
\left(\pa_3^2 \cdot S_2- \pa_3^2(S_2)\right)S_2^{-1}  =
\left( \pa_3(S_2) +S_2 \, \pa_3\right) \cdot \pa_3 \cdot S_2^{-1}
\een
On the other hand, recalling Corollary \ref{cor:anul} we get
\ben
(\Lambda-1) S_2^{-1} =  S_2^{-1} (\pa_2+q_2)^{-1} H_2\quad \in \quad
\E_{(2)} H
\een
and 
\ben
\pa_3 \cdot S_2^{-1} = S_2^{-1}\cdot \pa_2^{-1}\cdot H_1 \quad \in
\quad \E_{(2)} H. 
\een
This completes the proof of part b). The proof of part c) is similar,
except for the following point. It is convenient to prove that
$\L^*-\tfrac{L_3^2}{2}\in \E_{(3)}H^*$ where ${}^*$ is the involution
$\Lambda\mapsto -\Lambda$ and $\E_{(3)}H^*$ is the left ideal of
$\E_{(3)}$ generated by $H_i^*$ ($1\leq i\leq 3)$. The argument goes
in the same way as in part b). The only point worth mentioning is
the identity
\ben
(\Lambda-1) S_3^{-1} =  -\Big( (\Lambda+1) S_3^{-1}\Big)^* =
-S_3^{-1} (\pa_3+q_3)^{-1} H_3^*\quad \in \quad
\E_{(3)} H^*.
\qed
\een

\begin{corollary}\label{cor:proj-L}
  The operator $\mathcal{L}$ has the following properties.
  
  a) $\mathcal{L}\in \E'_{(\pm)}$ and
  $\mathcal{L}=\tfrac{(L_1^\pm)^{n-2}}{n-2}$ in the quotient ring $\E''_{(\pm)}$.

b) $\mathcal{L}\in \E'_{(2)}$ and
$\mathcal{L}=\tfrac{L_2^{2}}{2}$ in the quotient ring
$\E''_{(2)}$.

c) $\mathcal{L}\in \E'_{(3)}$ and 
$\mathcal{L}=\tfrac{L_3^2}{2}$ in the quotient ring
$\E''_{(3)}$.
\end{corollary}
\proof
We need only to check that  $L_1^\pm \in \E'_{(\pm)}$ and that
$L_a\in \E'_{(a)}$ ($a=2,3$). This however is a straightforward
computation using Corollary \ref{cor:anul}.
\qed

\section{The evolution of $\mathcal{L}$, $H_1$, $H_2$, and $H_3$}
\label{sec:HBE-lax-2}
Note that the operator $\mathcal{L}$ has the form
\ben
\mathcal{L} =   \Big(\sum_{i=1}^{n-3} (a_i\Lambda^i-\Lambda^{-i} a_i) \Big)
(\Lambda-\Lambda^{-1}) +
\frac{1}{2}\,\pa_2^2 + \frac{1}{2}\,\pa_3^2 
+\frac{1}{4}(c_2-c_3)(\Lambda+\Lambda^{-1}) +\frac{1}{2}(c_2+c_3)
\een
and that the coefficient $a_{n-3}=\tfrac{1}{n-2} e^{(n-2)\alpha}$
for some $\alpha\in \O_\epsilon[\![\mathbf{t}]\!]$. 
Let us denote by $\mathcal{M}$
the set of operators $\mathcal{L}$, $H_2$, and $H_3$. They depend on
the $n+1$ functions
\beq\label{parameters}
\Xi:=\{
a_i(x,\mathbf{t})\,  (1\leq i\leq n-4),\quad
\alpha(x,\mathbf{t}),\quad
q_j(x,\mathbf{t}) \, (j=2,3),\quad
c_k(x,\mathbf{t})\, (k=2, 3)\},
\eeq
which will be viewed as coordinates on $\mathcal{M}$. Let us point out
that the operator $H_1=\pa_2\pa_3+q_1$ is uniquely determined from
$H_2$ and $H_3$, because the relation $q_1+q_1[1]+2q_2q_3=0$ can be
solved uniquely for $q_1\in \O_\epsilon[\![\mathbf{t}]\!]$ in terms of
$q_2$ and $q_3$. 

We would like to prove that the operators
$\mathcal{L}$, $H_1$, $H_2$, and $H_3$
satisfy the Lax equations \eqref{Lax-pi_+} of the Extended D-Toda
  hierarchy. In other words, the Extended D-Toda hierarchy defines an
  infinite sequence of commuting flows on the manifold $\M$.
Let us give an outline of our argument. Put 
\ben
\mathcal{R}:=\CC[\pa_x^i \pa_2^j\pa_3^k (\xi)\,
(\xi\in \Xi,\, i,j,k\in \ZZ_{\geq 0}), e^{\pm \alpha} ][\![\epsilon]\!]
\een
for the ring of formal power series in $\epsilon$ whose coefficients are polynomials in the partial
derivatives  $\pa_x^i\pa_2^j\pa_3^k(\xi)$ where $\xi$ is one of the
functions in \eqref{parameters}. We refer to the elements of
$\mathcal{R}$ as differential polynomials. For brevity, put
$\pa_{i,k}=\tfrac{\pa }{\pa t_{i,k}}$. First, we are going to prove that
$P_{i,k,\xi}:=\pa_{i,k}(\xi)\in \mathcal{R}$ for all $\xi\in \Xi$. Then we will argue that
$P_{a,\xi}:=\pa_{a}(\xi)$ for $a=2,3$ can be expressed in terms of the
partial derivatives of \eqref{parameters} with respect to $x$ only. In
other words the ring $\R$ is in fact equal to
\ben
\CC[\pa_x^i (\xi)\, (\ \xi\in \Xi,\, i\in \ZZ_{\geq 0}), e^{\pm
  \alpha} ][\![\epsilon]\!].
\een

\subsection{Projections}

\begin{lemma}\label{le:od}
Let $\mathcal{A}$ be one of the rings $\mathcal{E}_{(\pm)}$,
$\mathcal{E}_{(2)}$, or $\mathcal{E}_{(3)}$. Then we have an
orthogonal decomposition
\ben
\mathcal{A} = \mathcal{A}^0 \oplus \mathcal{A}\, H. 
\een
\end{lemma}
\proof
Let us give the argument for the case $\A=\E_{(+)}$. The remaining
cases are similar. If $P\in \E_{(+)}$ then it can be written as
\ben
\sum_{j=j_0}^\infty \Lambda^{-j} \, p_j(\pa_2,\pa_3),
\een
where the coefficients $p_j$ are differential operators in $\pa_2$ and
$\pa_3$. Since $H_2=(\Lambda-1) \pa_2 -q_2 (\Lambda+1)$ we have 
\ben
\pa_2 = \iota_{\Lambda^{-1}}
\left( (\Lambda-1)^{-1} q_2 (\Lambda+1) + (\Lambda-1)^{-1} H_2
  \right).
\een
Similarly
\ben
\pa_3 = \iota_{\Lambda^{-1}}
\left( (\Lambda+1)^{-1} q_3 (\Lambda-1) + (\Lambda+1)^{-1} H_3
  \right). 
\een
Using these identities and that $[\pa_a,(\Lambda\pm 1)^{-1}H_b]\in B[\![\Lambda^{-1}]\!]$
we get that each differential operator $p_j(\pa_2,\pa_3)$ can be
decomposed as
$
p_j^{(0)}+
p_j^{(2)} (\Lambda-1)^{-1}H_2 +
p_j^{(3)}(\Lambda+1)^{-1} H_3,
$
where  
$
p_j^{(0)} \in
B[\![\Lambda^{-1}]\!]
$
and
$
\ p_j^{(2)},\ p_j^{(3)}\in B[\pa_2,\pa_3][\![\Lambda^{-1}]\!].
$
This proves that $\E_{(+)}=\E^0_{(+)}+\E_{(+)}H$. The sum must be
direct because if $P(x,\mathbf{t},\Lambda)\in \E^0_{(+)}$ anihilates the wave function
$\Psi_1(x,\mathbf{t},z)$, then $P(x,\mathbf{t},\Lambda)
S_1(x,\mathbf{t},\Lambda)=0$. However, the ring $\E^0_{(+)}$ does not
have zero divisors and $S_1\neq 0$, so $P=0$.
\qed

Let us denot by $\pi_\alpha: \E_{(\alpha)}\to \E^0_{(\alpha)}$
$(\alpha=\pm,2,3)$  the projection defined via the orthogonal
decomposition from Lemma \ref{le:od}. Note that if we have an operator
$P=\sum_{j,k,l} p_{j,k,l}(x,\mathbf{t}) \Lambda^{-j} \pa_2^k\pa_3^l\in
\E_{(\alpha)}$, then the projection $\pi_\alpha(P)$ is a
pseudo-difference (if $\alpha=\pm$) or
a pseudo-differential (if $\alpha=2,3$) operator whose coefficients are polynomials on
the shifted derivatives $\pa_2^b\pa_3^c(p_{j,k,l}[a])$, where $b,c\in \ZZ_{\geq 0}$ and $a\in
\ZZ$.

\subsection{The $t_{1,k}$-flows}
Let us define 
\ben
B^+_{1,k} :=
\left(
-\left( L_1^k \, \sum_{m=0}^\infty
  \Lambda^{-2m-1} \right)_{1,< 0} +
  \left(\sum_{m=0}^\infty
  \Lambda^{2m+1}\,  (L_1^\#)^k \right)_{1,<0} \right)\,
(\Lambda-\Lambda^{-1}).
\een
Note that $B^+_{1,k}\in B(\!(\Lambda^{-1})\!)=-L_1^k+B_{1,k}$ is a 
pseudo-difference operator. We claim that the coefficients of $B^+_{1,k}$ 
are differential polynomials in $\mathcal{R}$. Indeed, it is sufficient 
to prove that the coefficients of $L_1$ belong to $\mathcal{R}$. 
Recalling Corollary \ref{cor:proj-L} we get $\pi_{+}(\L) = 
\tfrac{L_1^{n-2}}{n-2}$. The coefficients of $\pi_+(\L)$ belong 
to $\mathcal{R}$. The rest of the proof is the same as the proof 
of Lemma \ref{le:L1}.

\begin{lemma}\label{le:S_1-flow1}
  The operator series $S_1$ satisfies the differential
  equations
  \ben
  \pa_{1,k} S_1 = B_{1,k}^+\, S_1,\quad k\geq 1.
  \een
\end{lemma}
\proof
Let us differentiate the HBEs in Proposition \ref{HQE-W} with respect
to $\pa_{1,k}$, substitute $\mathbf{t}'=\mathbf{t}$, and compare the
coefficients in front of the negative powers of $\Lambda$. We get
\ben
\pa_{1,k} S_1^+ \Lambda^{-1} S_1^- +
\left(
S_1^+ \Lambda^{k-1} S_1^-
-(S_1^-)^{\#} \Lambda^{1-k} (S_1^+)^\#
\right)_{1,<0}= 0. 
\een
It remains only to recall Proposition \ref{prop:dres}, a).
\qed
\begin{proposition}\label{prop:L-flow1}
  We have
  \ben
  \pa_{1,k} \L & = &  \pi_+([B^+_{1,k},\L]) \\
  \pa_{1,k} H_i & = & -\pi_+(H_i B^+_{1,k}),\quad (1\leq i\leq 3).
  \een
\end{proposition}
\proof
Recalling Corollary \ref{cor:anul} we have
$
H_2= (\Lambda-1) S_1 \cdot \pa_2\cdot S_1^{-1} ,
$
$
H_3= (\Lambda+1) S_1 \cdot \pa_3\cdot S_1^{-1} ,
$
and
\ben
H_1=\pa_2(\Lambda+1)^{-1} H_3 + \pa_3 (\Lambda-1)^{-1} H_2 -
(\Lambda-1)^{-1} H_2 
(\Lambda+1)^{-1} H_3. 
\een
The differential equations for $H_i$ follow from the above formulas
and Lemma \ref{le:S_1-flow1}.

It remain only to prove the differential equation for $\L$. Using that
\ben
\L = \frac{1}{n-2}\, S_1\Lambda^{n-2}S_1^{-1} + A_2 H_2 + A_3 H_3
\een
for some $A_2,A_3\in \E_{(+)}$ we get that modulo $\ker(\pi_+)$ the
derivative $\pa_{1,k}\L $ is given by
\ben
\frac{1}{n-2} [B_{1,k}, S_1\Lambda^{n-2}S_1^{-1}] - A_2 H_2 B_{1,k} -
A_3 H_3 B_{1,k},
\een
where we used Lemma \ref{le:S_1-flow1} and the differential equations
for $H_i$ which we have established already. The above expression
coincides with $[B_{1,k},\L]$ modulo terms in $\ker(\pi_+)$.
\qed

\subsection{The $t_{a,2l+1}$-flows}
Let us introduce the pseudo-difference operators
\ben
B^+_{2,k} = \frac{1}{2}\,
\iota_{\Lambda^{-1}} \,
\left(
  S_2 \pa_2^{k} (\Lambda-1)^{-1} S_2^{-1}
\right)_{2,[0]}
(\Lambda-\Lambda^{-1})
\een
and
\ben
B^+_{3,k} = \frac{1}{2}\,
\iota_{\Lambda^{-1}} \,
\left(
  S_3 \pa_3^{k} (\Lambda+1)^{-1} S_3^{-1}
\right)_{3,[0]}
(\Lambda-\Lambda^{-1})
\een
We claim that $B_{2,k}^+$ and $B_{3,k}^+$ are formal Laurent series in $\Lambda^{-1}$
whose coefficients belong to $\R$. Let us prove this for
$B_{2,k}^+$. The argument for $B_{3,k}^+$ is similar. Using Corollary
\ref{cor:anul}, b) we have
\ben
B_{2,k}^+ = L_2^k \, H_2^{-1} \, (\pa_2+ q_2) (\Lambda-\Lambda^{-1})=
L_2^k \,
\left(
  (1+\pa_2^{-1} q_2)^{-1} (1-\pa_2^{-1} q_2) \Lambda -1
\right)^{-1}
(\Lambda-\Lambda^{-1}).
\een
It is enough to prove that the coefficients of the pseudo-differential
operator $L_2$ belong to $\R$. On the other hand, $\pi_2(\L)=
\tfrac{L_2^2}{2}$. Solving for $L_2$ and using that the coefficients of 
$\pi_2(\L)$ belong to $\R$, we get that $L_2$ has coefficients in $\R$
(see the proof of Lemma \ref{le:La} for more details).
\begin{lemma}\label{le:S_1-flow2}
  We have
  \ben
  \pa_{a,2l+1} (S_1)  =  B_{a,2l+1}^+ \, S_1,\quad l\geq 0,\quad a=2,3.
  \een
\end{lemma}
\proof
Again we give the argument only for $a=2$the first differential equation
when $a=2$,
because the argument for the second one is identical.
Let us differentiate the HBEs in Proposition \ref{HQE-W} corresponding
to $k=0$ with respect
to $\pa_{2,2l+1}$, set $\mathbf{t}'=\mathbf{t}$, and compare the
coefficients in front of the negative powers of $\Lambda$. If $l>0$
then we get
\ben
\pa_{2,2l+1} (S_1) S_1^{-1} \iota_{\Lambda^{-1}}
(\Lambda-\Lambda^{-1})^{-1} =
\frac{1}{2} \sum_{m=1}^\infty (S_2 \pa_2 ^{2l+1} \Lambda^{-m} S_2^{-1})_{2,[0]}.
\een
This is equivalent to what we have to prove. If $l=0$ we get
\ben
\pa_{2} (S_1) S_1^{-1} \iota_{\Lambda^{-1}}
(\Lambda-\Lambda^{-1})^{-1} =
\frac{1}{2} \sum_{m=1}^\infty (\pa_2\cdot S_2 \Lambda^{-m} S_2^{-1})_{2,[0]}.
\een
It remains only to use that $L_2=S_2\pa_2S_2^{-1} =
\pa_2+O(\pa_2^{-1})$, i.e., the coefficients in front of $\pa_2^0$ of
$L_2$ is 0.
\qed

The same argument that we used to prove Proposition \ref{prop:L-flow1}
yields the following Proposition.
\begin{proposition}\label{prop:L-flow2}
  We have
  \ben
  \pa_{a,2l+1} \L & = &  \pi_+([B^+_{a,2l+1},\L]), \\
  \pa_{a,2l+1} H_i & = & -\pi_+(H_i B^+_{a,2l+1}),\quad (1\leq i\leq 3),
  \een
  where $a=2,3.$
\end{proposition}

\subsection{The derivations $\pa_2$ and $\pa_3$}\label{sec:der-23}
Let $\xi$ be one of the coordinate functions \eqref{parameters}.
We would like to prove that the derivative $\pa_a(\xi)$ ($a=2,3$) is a
formal power series in $\epsilon$ whose coefficients are differential
polynomials involving only the derivative $\pa_x$. To begin with let
us point out that by comparing the formulas for $q_a$ and $c_a$
($a=2,3$) in terms of the tau-function it follows that 
\beq\label{pa_a(q_b)}
\pa_2(q_3)=\pa_3(q_2) = \frac{e^{\epsilon \pa_x}-1}{e^{\epsilon
    \pa_x}+1} (q_2q_3)
\eeq
and that
\beq\label{pa_a(q_a)}
\pa_a(q_a) = \frac{1}{2}\,\left(1-e^{\epsilon\pa_x}\right)(c_a),\quad a=2,3.
\eeq
Let us check that the derivatives $\pa_2(a_i)$ ($1\leq i\leq
n-3$) and $\pa_2(c_a)$ $(a=2,3)$ can
be expressed in terms of $\pa_x$-differential polynomials. The
argument for $\pa_3(a_i)$ and $\pa_3(c_a)$ is identical. Let us
recall the differential equation
\beq\label{pa_2(L)}
\pa_2(\L) =  \pi_+([B_{2,1}^+,\L]).
\eeq
We have $B_{2,1}^+ = Q_2$ and the operator $\L$ has the form
\ben
&&
a_{n-3} \Lambda^{n-2} + a_{n-4} \Lambda^{n-3} + \sum_{k=2}^{n-4}
(a_{k-1} - a_{k+1}) \Lambda^{k} +
\Big(-a_2+\frac{1}{4}(c_2-c_3)\Big)\Lambda +\\
&&
\Big(-a_1-a_1[-1] + \frac{1}{2} (c_2+c_3)\Big) \Lambda^0+ \cdots,
\een
where the dots stand for terms involving only negative powers of
$\Lambda$. Let us split $[Q_2,\L]$ into sum
of two commutators $[Q_2,\L-\tfrac{1}{2}(\pa_2^2+\pa_3^2)]$ and
\ben
\frac{1}{2}
[Q_2, \pa_2^2+\pa_3^2] = -\frac{1}{2}(\pa_2^2+\pa_3^2)(Q_2) -
\pa_2(Q_2)\pa_2 -\pa_3(Q_2) \pa_3.
\een
The first commutator is already in $\E^{0}_{(+)}$ and it is a Laurent
series in $\Lambda^{-1}$ whose coefficients are differential
polynomials involving only $\pa_x$-derivatives. The projection $\pi_+$
of the second commutator is
\ben
-\frac{1}{2}(\pa_2^2+\pa_3^2)(Q_2) -
\pa_2(Q_2)Q_2 -\pa_3(Q_2) Q_3 .
\een
A straightforward computation, using formulas \eqref{pa_a(q_b)} and
\eqref{pa_a(q_a)}, shows that the above expression has 
leading order term of the type
\ben
&&\Big(\frac{1}{4}\pa_2 (c_2-c_2[-1]) -\frac{1}{2}\frac{1-e^{-\epsilon\pa_x}}{1+e^{\epsilon\pa_x}}\Big(\frac{e^{\epsilon\pa_x}-1}{e^{\epsilon\pa_x}+1}(q_2q_3)\cdot q_3+\frac{1}{2}(c_3-c_3[1])q_2\Big)+\\
&&+\frac{1}{2}(c_2-c_2[-1])q_2[-1]+
\frac{1}{2}(c_3-c_3[-1])q_3[-1]\Big) \Lambda^0 + O(\Lambda^{-1}). 
\een
Comparing the coefficients in front of $\Lambda^k$ in \eqref{pa_2(L)} for
$1\leq k\leq n-2$ we get that $\pa_2(a_i)$ $(1\leq i\leq n-3)$ and
$\pa_2(c_2-c_3)$ can be expressed as 
differential polynomials that involve only $\pa_x$-derivatives.  
Comparing the coefficients in front of $\Lambda^0$ in \eqref{pa_2(L)} we get that
\ben
\frac{1}{2}\pa_2(c_2+c_3) - \frac{1}{4}\pa_2 (c_2-c_2[-1])
= \frac{1}{8}(3+e^{-\epsilon\pa_x}) \pa_2(c_2+c_3)-\frac{1}{8}(1-e^{-\epsilon\pa_x})\pa_2(c_2-c_3)
\een
can be expressed in terms of differential polynomials that involve
only $\pa_x$-derivatives. The operator $3+e^{-\epsilon\pa_x}$ is
invertible, so the derivative $\pa_2(c_2+c_3)$ is also a differential
polynomial involving only $\pa_x$-derivatives.

\subsection{The extended flows}
Let us define the following operator series 
\ben
B_{0,l}^+:= (B_{0,l,1}^++B_{0,l,2}^++B_{0,l,3}^+)(\Lambda-\Lambda^{-1}),
\een
where 
\ben
B_{0,l,1}^+:= 
-\left(S_1\, \Big(\frac{ \Lambda^{(n-2)l}}{(n-2)^l l!}\left(
\epsilon \pa_x -h_l\right)
\Big)S_1^{-1} 
\sum_{m=0}^\infty\Lambda^{-2m-1} - 
\sum_{m=0}^\infty\Lambda^{2m+1}
\Big(
S_1\, \Big(\frac{ \Lambda^{(n-2)l}}{(n-2)^l l!}\left(
\epsilon \pa_x -h_l\right)
\Big)S_1^{-1} 
\Big)^{\#}
\right)_{1,<0},
\een
\ben
B_{0,l,2}^+:= \frac{1}{2}\iota_{\Lambda^{-1}}
\left(
S_2 \, \frac{\pa_2^{2l}}{2^l l!}\, \frac{\epsilon \pa_x}{\Lambda-1} \, S_2^{-1}
\right)_{2,[0]},
\een
and 
\ben
B_{0,l,3}^+:= \frac{1}{2}\iota_{\Lambda^{-1}}
\left(
S_3 \, \frac{\pa_3^{2l}}{2^l l!}\, \frac{\epsilon \pa_x}{\Lambda+1} \, S_3^{-1}
\right)_{3,[0]}.
\een
Note that $B_{0,l}^+=-A_{1,k}^++\pi_{+}(B_{0,l})$. Next we will prove that 
$B_{0,l}^+$ are Laurent series in $\Lambda^{-1}$ whose coefficients 
belong to $\R$. Note that the coefficients are apriori differential 
operators in $\pa_x$ of order 1. However, recalling Proposition 
\ref{HQE-W} with $k=l$ and $\mathbf{t}'=\mathbf{t}$ we see that 
the coefficients in front of
$\pa_x$ is $0$, i.e., the coefficients of $B_{0,l}^+$ viewed as a
Laurent series in $\Lambda^{-1}$ are scalar functions. Since we already
know that the coefficients of $L_1=S_1\Lambda S_1^{-1}$ and $L_a=S_a\pa_a
S_a^{-1}$ ($a=2,3$) are in $\R$ it is enough to prove the following
Lemma. 
\begin{lemma}\label{le:ell_i}
The coefficients of the operator series 
\ben
\ell_1:= \epsilon\pa_x (S_1)\, S_1^{-1},\quad
\ell_a:=\epsilon \pa_x(S_a) \, S_a^{-1}\ (a=2,3) 
\een
are in $\R$.
\end{lemma}
The proof of Lemma \ref{le:ell_i} is the same as the proof of
Proposition \ref{commuteonR1}, b) and Proposition \ref{commuteonR23},
b). 
\begin{lemma}
  The wave operator $S_1$ satisfies the differential equations
  \ben
  \pa_{0,l} (S_1) = B_{0,l}^+\, S_1,\quad \ell \geq 1.
  \een
\end{lemma}
\proof
Let us differentiate the HBEs corresponding to $k=0$ in Proposition \ref{HQE-W} with respect
to $t_{0,l}$, put $\mathbf{t}'=\mathbf{t}$, and compare the
coefficients in front of the negative powers of $\Lambda$.
We get exactly the identity stated in the lemma.
\qed

Just like before, this lemma alows us to prove the following
proposition.
\begin{proposition}\label{prop:L-flow0}
  We have
  \ben
  \pa_{0,l} \L & = &  \pi_+([B^+_{0,l},\L]), \\
  \pa_{0,l} H_i & = & -\pi_+(H_i B^+_{0,l}),\quad (l\geq 1).
  \een
\end{proposition}
Let us point out that $B_{0,0}^+ = \ell_1$. Therefore, if we put
$\pa_{0,0}=\epsilon \pa_x$, then the equations in Proposition
\ref{prop:L-flow0} will hold for $l=0$ as well.

\bibliographystyle{amsalpha}

\section{Appendix}
In this section, the examples of the projections and flows equations are given. 

\subsection{Example of projections} $\bullet$
The projection $\pi_+$ of $\pa_2^k$,
\begin{align*}
&\pi_{+}(\pa_2)=q_2[-1]+\sum_{l=1}^{+\infty}(q_2[-l-1]+q_2[-l])\Lambda^{-l},\\
&\pi_{+}(\pa_2^2)=\Big(\pa_2(q_2)+q_2^2\Big)[-1]+\Big(\pa_2(q_2[-1]+q_2)+(q_2[-1]+q_2)^2\Big)[-1]\Lambda^{-1}\\
&\quad\quad\quad+\Big(\pa_2(q_2[-3]+q_2[-2])+(q_2[-3]+q_2[-2])(2q_2[-1]+q_2[-2]+q_2[-3])\Big)\Lambda^{-2}+O(\Lambda^{-3}),\\
&\pi_{+}(\pa_2^3)=\Big(\pa_2^2(q_2)+3\pa_2(q_2)\cdot q_2+q_2^3\Big)[-1]+O(\Lambda^{-1}),\\
&\pi_{+}(\pa_2^4)=\Big(\pa_2^3(q_2)+4\pa_2^2(q_2)\cdot q_2+3(\pa_2q_2)^2+6\pa_2(q_2)q_2^2+q_2^4\Big)[-1]+O(\Lambda^{-1}).
\end{align*}
$\bullet$ The projection $\pi_+$ of $\pa_3^k$,
\begin{align*}
&\pi_{+}(\pa_3)=q_3[-1]+\sum_{k=1}^{+\infty}(-1)^k(q_3[-k-1]+q_3[-k])\Lambda^{-k},\\
&\pi_{+}(\pa_3^2)=\Big(\pa_3(q_3)+q_3^2\Big)[-1]-\Big(\pa_3(q_3[-1]+q_3)+(q_3[-1]+q_3)^2\Big)[-1]\Lambda^{-1}\\
&\quad\quad\quad+\Big(\pa_2(q_3[-3]+q_3[-2])+(q_3[-3]+q_3[-2])(2q_3[-1]+q_3[-2]+q_3[-3])\Big)\Lambda^{-2}+O(\Lambda^{-2}),\\
&\pi_{+}(\pa_3^3)=\Big(\pa_3^2(q_3)+3\pa_3(q_3)\cdot q_3+q_3^3\Big)[-1]+O(\Lambda^{-1}),\\
&\pi_{+}(\pa_3^4)=\Big(\pa_3^3(q_3)+4\pa_3^2(q_3)\cdot q_3+3(\pa_3q_3)^2+6\pa_3(q_3)q_3^2+q_3^4\Big)[-1]+O(\Lambda^{-1}).
\end{align*}
$\bullet$ The projection $\pi_2$ of $\Lambda^{\pm k}$,
\begin{align*}
&\pi_{2}(\Lambda)=1+2q_2\pa_2^{-1}+2(q_2^2-\pa_2(q_2))\pa_2^{-2}+O(\pa_2^{-3}),\\
&\pi_{2}(\Lambda^{-1})=1-2q_2[-1]\pa_2^{-1}+2(q_2^2-\pa_2(q_2))[-1]\pa_2^{-2}+O(\pa_2^{-3}),\\
&\pi_{2}(\Lambda^2)=1+2(q_2+q_2[1])\pa_2^{-1}+2\Big((q_2+q_2[1])^2-\pa_2(q_2+q_2[1])\Big)\pa_2^{-2}+O(\pa_2^{-3}),\\
&\pi_{2}(\Lambda^{-2})=1-2(q_2+q_2[-1])[-1]\cdot\pa_2^{-1}+2\Big((q_2+q_2[-1])^2-\pa_2(q_2+q_2[-1])\Big)[-1]\cdot\pa_2^{-2}+O(\pa_2^{-3}).
\end{align*}
$\bullet$ The projection $\pi_2$ of $\pa_3^k$,
\begin{align*}
&\pi_{2}(\pa_3)=-\pa_2^{-1}q_1,\\
&\pi_{2}(\pa_3^2)=-\pa_3(q_1)\pa_2^{-1}+(\pa_2\pa_3(q_1)+q_1^2)\pa_2^{-2}+O(\pa_2^{-2}).
\end{align*}
$\bullet$ The projection $\pi_3$ of $\Lambda^{\pm k}$,
\begin{align*}
&\pi_{3}(\Lambda)=-1-2q_3\pa_3^{-1}+2(\pa_3(q_3)-q_3^2)\pa_3^{-2}+O(\pa_3^{-3}),\\
&\pi_{3}(\Lambda^{-1})=-1+2q_3[-1]\pa_3^{-1}-2(q_3^2+\pa_3(q_3))[-1]\pa_2^{-2}+O(\pa_3^{-3}),\\
&\pi_{3}(\Lambda^2)=1+2(q_3+q_3[1])\pa_3^{-1}+2\Big((q_3+q_3[1])^2-\pa_3(q_3+q_3[1])\Big)\pa_3^{-2}+O(\pa_3^{-3}),\\
&\pi_{2}(\Lambda^{-2})=1-2(q_3+q_3[-1])[-1]\cdot\pa_3^{-1}+2\Big((q_3+q_3[-1])^2+\pa_3(q_3+q_3[-1])\Big)[-1]\cdot\pa_3^{-2}+O(\pa_3^{-3}).
\end{align*}
$\bullet$ The projection $\pi_3$ of $\pa_2^k$,
\begin{align*}
\pi_{3}(\pa_2)&=-\pa_3^{-1}q_1,\\
\pi_{3}(\pa_2^2)&=-\pa_2(q_1)\pa_3^{-1}+(\pa_2\pa_3(q_1)+q_1^2)\pa_3^{-2}+O(\pa_3^{-2}).
\end{align*}
\subsection{Examples of Lax operator}
Take $n=4$ as an example. The Lax operator will be
\begin{align*}
\mathcal{L}=a\Lambda^{2}+\frac{1}{4}(c_2-c_3)\Lambda+\frac{1}{2}(c_2+c_3)-a-a[-1]+\frac{1}{4}(c_2-c_3)\Lambda^{-1}+a[-1]\Lambda^{-2}+\frac{1}{2}\pa_2^2+\frac{1}{2}\pa_3^2,
\end{align*}
where $a=\frac{1}{2}e^{\frac{\alpha}{2}}$. Then we will obtain
\begin{align*}
\pi_+(\mathcal{L})=&a\Lambda^{2}+\frac{1}{4}(c_2-c_3)\Lambda+v_{1,0}+v_{1,1}\Lambda^{-1}+v_{1,2}\Lambda^{-2}+O(\Lambda^{-3}),\\
\pi_i(\mathcal{L})=&\frac{1}{2}\pa_i^2+c_i+v_{i,1}\pa_i^{-1}+v_{i,2}\pa_{i}^{-2}+O(\pa_i^{-3}),
\end{align*}
where
\begin{align*}
v_{1,0}=&\frac{1}{4}(c_2+c_3)+\frac{1}{4}(c_2+c_3)[1]-a-a[-1]+\frac{1}{2}\Big(q_2^{2}+q_3^2\Big)[-1],\\
v_{1,1}=&\frac{1}{4}(c_2-c_3)[-2]+\frac{1}{2}\sum_{l=2}^3(-1)^l(q_l[-2]+q_l[-1])^2,\\
v_{1,2}=&a[-1]+\frac{1}{4}(c_2+c_3)[-3]-\frac{1}{4}(c_2+c_3)[-1]\\
&+\frac{1}{2}\sum_{l=2}^3(q_l[-3]+q_l[-2])(2q_l[-1]+q_l[-2]+q_l[-3]),\\
v_{i,1}=&2a(q_i+q_i[1])+\frac{(-1)^i}{2}(c_2-c_3)(q_i+q_i[-1])-2a[-1](q_i[-1]+q_i[-2]),\\
v_{i,2}=&2a((q_i+q_i[1])^2-\pa_i(q_i+q_i[1]))+\frac{(-1)^i}{2}(c_2-c_3)(q_i^2-q_i[-1]^2\\
&-\pa_i(q_i-q_i[-1]))+2a[-1]((q_i[-1]+q_i[-2])^2-\pa_i(q_i[-1]+q_i[-2])).
\end{align*}

The expressions of $B_{i,k}$ are listed as follows
\begin{align*}
&B_{0,1}=\mathcal{L}\cdot\epsilon\pa_x-\frac{1}{2}(a_{2,1}\pa_2+a_{3,1}\pa_3)+\sum_{l=-2}^2b_{0,l}\Lambda^{l}\\
&B_{1,1}=e^\beta (\Lambda-\Lambda^{-1}),\quad B_{2,1}=\pa_2,\quad B_{3,1}=\pa_3\\
&B_{1,2}=\Big(2a\Lambda+\frac{1}{2}(c_2-c_3)-2a[-1]\Lambda^{-1}\Big)\cdot(\Lambda-\Lambda^{-1})\\
&B_{2,3}=\pa_2^3+3c_2\pa_2,\quad B_{3,3}=\pa_3^2+3c_3\pa_3,
\end{align*}
where
\begin{align*}
&\beta=\frac{1}{1+e^{\epsilon\pa_x}}\Big(\frac{\alpha}{2}\Big),\quad b_{0,2}=-a\Big(\frac{1}{2}+a_{1,0}[2]\Big),\\
&b_{0,1}=-\frac{1}{4}(c_2-c_3)(\frac{1}{2}+a_{1,0}[1])-a a_{1,1}[2],\\
&b_{0,0}=-(b_{0,2}+b_{0,-2})+\epsilon\pa_x a[-1]-a_{1,0}[1]\cdot a[-1]-a_{2,1}-\frac{1}{2}a_{2,2},\\
&b_{0,-1}=-b_{0,1}+\frac{1}{4}\epsilon\pa_x(c_2-c_3),\\
&b_{0,-2}=\Big(\frac{1}{2}+a_{1,0}[2]\Big)\cdot a[-1]+\epsilon\pa_x a[-1].
\end{align*}
Here $a_{k,l}$ comes from $\ell_i=\epsilon\pa_x(S_i)\cdot S_i^{-1}$ with $i=1,2,3$,
\begin{align*}
&\ell_1=a_{1,0}+a_{1,1}\Lambda^{-1}+O(\Lambda^{-2})\\
&\ell_2=a_{2,1}\pa_2^{-1}+a_{2,2}\pa_2^{-2}+O(\pa_2^{-3}),\quad \ell_3=a_{3,1}\pa_3^{-1}+a_{3,2}\pa_3^{-2}+O(\pa_3^{-3}),
\end{align*}
where
\begin{align*}
&a_{1,0}=\frac{\epsilon\pa_x}{1-e^{\epsilon\pa_x}}(\beta),\quad a_{1,1}=\frac{1}{2}e^{\beta}\frac{\epsilon\pa_x}{1-e^{2\epsilon\pa_x}}\Big((c_2-c_3)e^{-\beta}\Big),\quad a_{i,1}=\frac{2\epsilon\pa_x}{e^{\epsilon\pa_x}-1}(q_i),\\ &a_{i,2}=\frac{2\epsilon\pa_x}{e^{\epsilon\pa_x}-1}\Big(-\pa_i(q_i)-q_i^2+\sum_{m=0}^{+\infty}\sum_{i=1}^{m+1}\sum_{j=0}^{i-1}\binom{i-1}{j}(\epsilon\pa_x)^j(a_{i,1})(\epsilon\pa_x)^{m-j}(a_{i,1})\Big),\quad i=2,3.
\end{align*}
\subsection{Examples of flows of $t_{i,1}$ for $i=1,2,3$} Flows of $t_{1,1}$ are
\begin{align*}
&\pa_{1,1}(a)=\frac{1}{4}\Big(e^{\beta}\cdot (c_2-c_3)[1]-e^{\beta[1]}\cdot(c_2-c_3)\Big),\\
&\pa_{1,1}(q_2)=e^{\beta[1]}(q_2+q_2[1])-e^{\beta}(q_2+q_2[-1]),\\ &\pa_{1,1}(q_3)=e^{\beta}(q_3+q_3[-1])-e^{\beta[1]}(q_3+q_3[1]),\\
&\pa_{1,1}(c_2)=c_2[1]-c_2[-1],\quad \pa_{1,1}(c_3)=c_3[-1]-c[1].
\end{align*}
Flows of $t_{2,1}$ are
\begin{align*}
\pa_2(a)=&a(q_2[-1]-q_2[1]), \quad \pa_2(q_2)=\frac{1}{2}(c_2-c_2[1]),
\quad \pa_2(q_3)=\frac{e^{\epsilon\pa_x}-1}{e^{\epsilon\pa_x}+1}(q_2q_3),\\
\pa_2(c_2)=&\frac{1}{1+e^{-\epsilon\pa_x}}\Big(-(c_2-c_3)q_2+(c_2-c_3)[-1]\cdot q_2[-2]-4a(q_2+q_2[1])\\
&+4a[-1](q_2[-2]-q_2)+4a[-2](q_2[-2]+q_2[-3])-q_2[-1](c_2[-1]-c_2)\\
&-2q_3[-1]\frac{1-e^{-\epsilon\pa_x}}{1+e^{\epsilon\pa_x}}(q_2q_3)\Big),\\
\pa_2(c_3)=&\frac{1}{1+e^{-\epsilon\pa_x}}\Big(q_2[-1](c_3-c_3[-1])-2q_3[-1]\frac{1-e^{-\epsilon\pa_x}}{1+e^{\epsilon\pa_x}}(q_2q_3)\Big).
\end{align*}
And flows of $t_{3,1}$ are given by
\begin{align*}
\pa_3(a)=&a(q_3[-1]-q_3[1]), \quad \pa_3(q_3)=\frac{1}{2}(c_3-c_3[1]),
\quad \pa_3(q_3)=\frac{e^{\epsilon\pa_x}-1}{e^{\epsilon\pa_x}+1}(q_2q_3),\\
\pa_3(c_2)=&\frac{1}{1+e^{-\epsilon\pa_x}}\Big(q_3[-1](c_2-c_2[-1])-2q_2[-1]\frac{1-e^{-\epsilon\pa_x}}{1+e^{\epsilon\pa_x}}(q_2q_3)\Big),\\
\pa_3(c_3)=&\frac{1}{1+e^{-\epsilon\pa_x}}\Big((c_2-c_3)q_3-(c_2-c_3)[-1]\cdot q_3[-2]-4a(q_3+q_3[1])\\
&+4a[-1](q_3[-2]-q_3)+4a[-2](q_3[-2]+q_3[-3])-q_3[-1](c_3[-1]-c_3)\\
&-2q_2[-1]\frac{1-e^{-\epsilon\pa_x}}{1+e^{\epsilon\pa_x}}(q_2q_3)\Big).
\end{align*}
\subsection{Examples of flows of $t_{1,2}$, $t_{2,3}$ and $t_{3,2}$}
Flows of $t_{1,2}$ are
\begin{align*}
\pa_{1,2}(a)=&\frac{1}{2}a\Big((c_2+c_3)[2]+(c_2+c_3)[1]-(c_2+c_3)-(c_2+c_3)[-1]\\
&+2q_2[1]^2+2q_3[1]^2-2q_2[-1]^2-2q_3[-1]^2-8a[1]+8a[-1]\Big),\\
\pa_{1,2}(q_2)=& 2a[1](q_2[1]+q_2[2])+2a(q_2[-1]-q_2[1])-2a[-1](q_2[-1]+q_2[-2])\\
&-\frac{1}{2}(c_2-c_3)(q_2[-1]+q_2)+\frac{1}{2}(c_2-c_3)[1]\cdot (q_2+q_2[1]),\\
\pa_{1,2}(q_3)=& 2a[1](q_3[1]+q_3[2])+2a(q_3[-1]-q_3[1])-2a[-1](q_3[-1]+q_3[-2])\\
&+\frac{1}{2}(c_2-c_3)(q_3[-1]+q_3)-\frac{1}{2}(c_2-c_3)[1]\cdot (q_3+q_3[1]),\\
\pa_{1,2}(c_2)=&-2a\Big(2(q_2+q_2[1])^2+c_2-c_2[2]\Big)+2a[-1]\Big(2(q_2[-1]+q_2[-2])^2+c_2-c_2[-2]\Big)\\
&-\frac{1}{2}(c_2-c_3)\Big(2q_2[-1]^2-2q_2^2+c_2[-1]-c_2[1]\Big),\\
\pa_{1,2}(c_3)=&2a\Big(2(q_3+q_3[1])^2+c_3[2]-c_3\Big)-2a[-1]\Big(2(q_3[-1]+q_3[-2])^2+c_3[-2]-c_3\Big)\\
&+\frac{1}{2}(c_2-c_3)\Big(2q_3[-1]^2-2q_3^2+c_3[-1]-c_3[1]\Big).
\end{align*}
The flows of $t_{2,3}$ are given by
\begin{align*}
\pa_{2,3}(a)=&a\Big(\pa_2^2(q_2[-1]-q_2[1])+3\pa_2(q_2[-1])\cdot q_2[-1]-3\pa_2(q_2[1])\cdot q_2[1]\\
&+q_2[-1]^3-q_2[1]^3+3c_2\cdot q_2[-1]-3c_2[2]\cdot q_2[1]\Big),\\
\pa_{2,3}(c_2)=&\pa_2^3(c_2)+3\pa_2^2(v_{2,1})+3\pa_2(v_{2,2})+3c_2\pa_2(c_2),\\
\pa_{2,3}(c_3)=&\pa_2^2\pa_3(q_1)+3\pa_3(c_2q_1),\\
\pa_{2,3}(q_2)=&\pa_2^3(q_2)+3q_2\cdot\pa_2^2(q_2)+3\pa_2(q_2)\cdot q_2^2+3\pa_2(c_2[1]q_2)+3(\pa_2(q_2))^2\\
\pa_{2,3}(q_3)=&-(3c_2[1]+q_2^2)\cdot (q_2q_3+q_1[1])+3\pa_3(c_2[1])q_2-\pa_2(q_2)\cdot(3q_2q_3+q_1[1])\\
&-2\pa_2(q_1[1])\cdot q_2-\pa_2^2(q_1[1])-q_3\cdot\pa_2^2(q_2).
\end{align*}
The flows of $t_{3,3}$ are given by
\begin{align*}
\pa_{3,3}(a)=&a\Big(\pa_3^2(q_3[-1]-q_3[1])+3\pa_3(q_3[-1])\cdot q_3[-1]-3\pa_3(q_3[1])\cdot q_3[1]\\
&+q_3[-1]^3-q_3[1]^3+3c_3\cdot q_3[-1]-3c_3[2]\cdot q_3[1]\Big),\\
\pa_{3,3}(c_2)=&\pa_2\Big(\pa_3^2(q_1)+3c_3q_1\Big),\\
\pa_{3,3}(c_3)=&\pa_3^3(c_3)+3\pa_3^2(v_{3,1})+3\pa_3(v_{3,2})+3c_3\pa_3(c_3),\\
\pa_{3,3}(q_2)=&-(3c_3[1]+q_3^2)\cdot (q_2q_3+q_1[1])+3\pa_2(c_3[1])q_3-\pa_3(q_3)\cdot(3q_2q_3+q_1[1])\\
&-2\pa_3(q_1[1])\cdot q_3-\pa_3^2(q_1[1])-q_2\cdot\pa_3^2(q_3),\\
\pa_{3,3}(q_3)=&\pa_3^3(q_3)+3q_3\cdot\pa_3^2(q_3)+3\pa_3(q_3)\cdot q_3^2+3\pa_3(c_3[1]q_3)+3(\pa_3(q_3))^2.
\end{align*}
\subsection{Examples of flows of $t_{0,1}$}
Flows of $t_{0,1}$ are
\begin{align*}
\pa_{0,1}(a)=&a\cdot\epsilon\pa_x v_{1,0}[2]
+\frac{1}{16}(c_2-c_3)\cdot\epsilon\pa_x(c_2-c_3)[1]
+v_{1,0}\cdot\epsilon\pa_x(a)+b_{0,2}(v_{1,0}[2]-v_{1,0})\\
&+\frac{1}{4}\Big(b_{0,1}\cdot(c_2-c_3)[1]-(c_2-c_3)\cdot b_{1,0}[1]\Big)
+a\Big(b_{0,0}-b_{0,0}[2]-\frac{1}{2}(a_{2,1}\cdot q_2[-1]\\
&+a_{3,1}\cdot q_3[-1])+\frac{1}{2}(a_{2,1}[2]q_2[1]+a_{3,1}[2]\cdot q_3[1])\Big),
\\
\pa_{0,1}(c_2)=&\frac{1}{2}\epsilon\pa_x\Big(\pa_2^2(c_2)+c_2^2
+2\pa_2(v_{2,1})+v_{2,2}\Big)
-\frac{1}{2}a_{2,1}\pa_2(c_2)-\frac{1}{2}\pa_2^2(b_{0,0,2})-\pa_2(b_{0,1,2}),\\
\pa_{0,1}(c_3)=&\frac{1}{2}\epsilon\pa_x\Big(\pa_3^2(c_3)+c_2^2
+2\pa_3(v_{3,1})+v_{3,2}\Big)
-\frac{1}{2}a_{3,1}\pa_3(c_3)-\frac{1}{2}\pa_3^2(b_{0,0,3})-\pa_3(b_{0,1,3}),\\
\pa_{0,1}(q_2)=&\pa_2\Big(b_{0,0}[1]-b_{0,1}-\frac{1}{2}a_{2,1}[1]q_2-\frac{1}{2}a_{3,1}[1]q_3
\Big)+q_2\Big(b_{0,2}[1]-b_{0,2}+b_{0,1}[1]-2b_{0,1}\Big)\\
&+q_2[1]\Big(b_{0,2}[1]-b_{0,2}+b_{0,1}[1]\Big)-\epsilon\pa_xA_2,\\
\pa_{0,1}(q_3)=&\pa_3\Big(b_{0,0}[1]+b_{0,1}-\frac{1}{2}a_{2,1}[1]q_2-\frac{1}{2}a_{3,1}[1]q_3
\Big)+q_2\Big(b_{0,2}[1]+b_{0,2}+b_{0,1}[1]+2b_{0,1}\Big)\\
&+q_2[1]\Big(b_{0,2}[1]+b_{0,2}+b_{0,1}[1]\Big)-\epsilon\pa_xA_3.
\end{align*}
Here
\begin{align*}
b_{0,0,2}=&\epsilon\pa_x\Big(a[-1]+\frac{1}{4}(c_2-c_3)\Big)-a_{1,0}[1]\cdot a[-1]-a_{2,1}-\frac{1}{2}a_{2,2},\\
b_{0,0,3}=&\epsilon\pa_x\Big(a[-1]-\frac{1}{4}(c_2-c_3)\Big)-a_{1,0}[1]\cdot a[-1]-a_{2,1}-\frac{1}{2}a_{2,2},\\
b_{0,1,2}=&2b_{0,2}\Big(q_2+q_2[1]\Big)+2b_{0,1} q_2+2b_{0,-1}q_2[-1]-2b_{0,-2}\Big(q_2[-1]+q_2[-2]\Big)+\frac{1}{2}a_{3,1}q_1,\\
b_{0,1,3}=&2b_{0,2}\Big(q_3+q_3[1]\Big)-2b_{0,1} q_3-2b_{0,-1}q_3[-1]-2b_{0,-2}\Big(q_3[-1]+q_3[-2]\Big)+\frac{1}{2}a_{2,1}q_1,\\
A_2=&\pa_2\Big(\frac{1}{2}(c_2+c_3)[1]-a-a[1]-\frac{1}{4}(c_2-c_3)\Big)+q_2\Big(a[1]-a
+\frac{1}{4}(c_2-c_3)[1]-\frac{1}{2}(c_2-c_3)\\
&-\frac{1}{2}(\pa_2(q_2)+\pa_3(q_3)+q_2^2+q_3^2)\Big)+
q_2[1]\cdot\Big(a[1]-a+\frac{1}{4}(c_2-c_3)[1]\Big),\\
A_3=&\pa_3\Big(\frac{1}{2}(c_2+c_3)[1]-a-a[1]-\frac{1}{4}(c_2-c_3)\Big)+q_3\Big(a[1]-a
-\frac{1}{4}(c_2-c_3)[1]+\frac{1}{2}(c_2-c_3)\\
&+\frac{1}{2}(\pa_2(q_2)+\pa_3(q_3)+q_2^2+q_3^2)\Big)+
q_3[1]\cdot\Big(a[1]-a-\frac{1}{4}(c_2-c_3)[1]\Big).
\end{align*}

\end{document}